\documentclass[a4paper,11pt]{article}
\usepackage{amssymb,amscd,amsfonts,amsbsy}
\usepackage{enumerate}
\usepackage{epsfig}
\input amssymb.sty
\input euscript.sty
\input diagrams.tex

\def\ring{\mathaccent"0017 }

\newcommand{\RR}{{\mathbb{R}}}
\newcommand{\NN}{{\mathbb{N}}}
\newcommand{\ZZ}{{\mathbb{Z}}}
\newcommand{\CC}{{\mathbb{C}}}
\newcommand{\meanint}{{\int{\mkern-19mu}-}}

\newtheorem{proposition}{Proposition}[section]
\newtheorem{theorem}[proposition]{Theorem}
\newtheorem{lemma}[proposition]{Lemma}
\newtheorem{corollary}[proposition]{Corollary}

\newtheorem{definition}{Definition}[section]

\begin{document}

\title{{The Dirichlet problem in Lipschitz domains with boundary data 
in Besov spaces for higher order elliptic systems with rough coefficients}
\thanks{2000 {\it Math Subject Classification.} Primary: 35G15, 35J55, 35J40
Secondary 35J67, 35E05, 46E39.
\newline
{\it Key words}: higher order elliptic systems, Besov spaces, 
weighted Sobolev spaces, mean oscillations, BMO, VMO, Lipschitz domains, 
Dirichlet problem
\newline
The work of authors was supported in part by from NSF DMS and FRG grants
as well as from the Swedish National Science Research Council}} 
 
\author{V.\, Maz'ya, M.\, Mitrea and T.\, Shaposhnikova}

\date{~}

\maketitle

\begin{abstract}
We settle the issue of well-posedness for the Dirichlet problem for 
a higher order elliptic system ${\mathcal L}(x,D_x)$ with 
complex-valued, bounded, measurable coefficients in a Lipschitz domain 
$\Omega$, with boundary data in Besov spaces. 

The main hypothesis under which our principal result is established is 
in the nature of best possible and requires that, at small scales, the mean 
oscillations of the unit normal to $\partial\Omega$ and of the coefficients 
of the differential operator ${\mathcal L}(x,D_x)$ are not too large. 
\end{abstract}

\section{Introduction}
\setcounter{equation}{0}
A fundamental theme in the theory of partial differential equations, which 
has profound and intriguing connections with many other subareas of analysis, 
is the well-posedness of various classes of boundary value problems under 
sharp smoothness assumptions on the boundary of the domain and on the 
coefficients of the corresponding differential operator. 
In this paper we initiate a program broadly aimed at extending the scope 
of the agenda set forth by Agmon, Douglis, Nirenberg and Solonnikov
(cf. {\bf\cite{ADN}}, {\bf\cite{Sol1}}, {\bf\cite{Sol2}}) in connection with
general elliptic boundary value problems on Sobolev-Besov scales, 
as to allow minimal smoothness assumptions (on the underlying domain 
and on the coefficients of the differential operator).
Our main result is the solvability of the Dirichlet problem for general higher 
order elliptic systems in divergence form, with complex-valued, bounded, 
measurable coefficients in Lipschitz domains, and for boundary data in 
Besov spaces. In order to be more specific we need to introduce some notation. 

Let $m,l\in\NN$ be two fixed integers and, for a bounded Lipschitz domain 
$\Omega$ in $\mathbb{R}^n$ (a formal definition is given in \S{6.1})
with outward unit normal $\nu=(\nu_1,...,\nu_n)$ consider the Dirichlet problem
for the operator 

\begin{equation}\label{LOL}
{\mathcal L}(X,D_X)\,{\mathcal U}
:=\sum_{|\alpha|=|\beta|=m}D^\alpha(A_{\alpha\beta}(X)D^\beta{\mathcal U})
\end{equation}

\noindent i.e., 

\begin{equation}\label{e0}
\left\{
\begin{array}{l}
\displaystyle{\sum_{|\alpha|=|\beta|=m}
D^\alpha(A_{\alpha\beta}(X)\,D^\beta\,{\mathcal U})}=0
\qquad\mbox{for}\,\,X\in\Omega,
\\[28pt] 
{\displaystyle\frac{\partial^k{\mathcal U}}{\partial\nu^k}}=g_k
\,\,\quad\mbox{on}\,\,\partial\Omega,\qquad 0\leq k\leq m-1.
\end{array}
\right.
\end{equation}
 
\noindent Here and elsewhere, 
$D^\alpha=(-i\partial/\partial x_1)^{\alpha_1}\cdots
(-i\partial/\partial x_n)^{\alpha_n}$ if $\alpha=(\alpha_1,...,\alpha_n)$. 
The coefficients $A_{\alpha\beta}$ are $l\times l$ matrix-valued functions 
with complex entries satisfying 

\begin{equation}\label{A-bdd}
\sum_{|\alpha|=|\beta|=m}\|A_{\alpha\beta}\|_{L_\infty(\Omega)}\leq\kappa_1
\end{equation}

\noindent for some finite constant $\kappa_1$, and such that 
the coercivity condition  

\begin{equation}\label{coercive}
\Re\,\int_\Omega\sum_{|\alpha|=|\beta|=m}\langle A_{\alpha\beta}(X) 
D^\beta\,{\mathcal U}(X),\,D^\alpha\,{\mathcal U}(X)\rangle\,dX 
\geq\kappa_0\sum_{|\alpha|=m}\|D^\alpha\,{\mathcal U}\|^2_{L_2(\Omega)}
\end{equation}

\noindent with $\kappa_0=const>0$ holds for all $\CC^l$-valued 
functions ${\mathcal U}\in C^\infty_0(\Omega)$. Throughout the paper, 
$\Re\,z$ denotes the real part of $z\in\CC$ and $\langle\cdot,\cdot\rangle$ 
stands for the canonical inner product in $\mathbb{C}^l$. 

Since, generally speaking, $\nu$ is merely bounded and measurable,
care should be exercised when defining iterated normal derivatives. For the 
setting we have in mind it is natural to take
$\partial^k/\partial\nu^k:=(\sum_{j=1}^n\xi_j\partial/\partial x_j)^k
\mid_{\xi=\nu}$ or, more precisely, 

\begin{equation}\label{nuk}
\frac{\partial^k{\mathcal U}}{\partial\nu^k}
:=i^k\sum_{|\alpha|=k}\frac{k!}{\alpha!}\,
\nu^\alpha\,{\rm Tr}\,[D^\alpha{\mathcal U}],\qquad 0\leq k\leq m-1,
\end{equation}
 
\noindent where ${\rm Tr}$ is the boundary trace operator and
$\nu^\alpha:=\nu_1^{\alpha_1}\cdots\nu_n^{\alpha_n}$ if 
$\alpha=(\alpha_1,...,\alpha_n)$.
With $\rho(X):={\rm dist}\,(X,\partial\Omega)$ and $p\in(1,\infty)$, 
$a\in(-1/p,1-1/p)$ fixed, a solution for (\ref{e0}) is sought in 
$W^{m,a}_p(\Omega)$, defined as the space of vector-valued functions for which 

\begin{equation}\label{W-Nr}
\Bigl(\sum_{0\leq|\alpha|\leq m}\int_\Omega|D^\alpha{\mathcal U}(X)|^p 
\rho(X)^{pa}\,dX\Bigr)^{1/p}<\infty.
\end{equation}

\noindent In particular, as explained later on, the traces in (\ref{nuk}) 
exist in the Besov space $B_p^{s}(\partial\Omega)$, where $s:=1-a-1/p\in(0,1)$,
for any ${\mathcal U}\in W^{m,a}_p(\Omega)$. 
Recall that, with $d\sigma$ denoting the area element on $\partial\Omega$,  

\begin{equation}\label{Bes-xxx}
g\in B_p^s(\partial\Omega)\Leftrightarrow 
\|g\|_{B_p^s(\partial\Omega)}:=\|f\|_{L_p(\partial\Omega)}
+\Bigl(\int_{\partial\Omega}\int_{\partial\Omega}
\frac{|g(X)-g(Y)|^p}{|X-Y|^{n-1+sp}}\,d\sigma_Xd\sigma_Y\Bigr)^{1/p}<\infty.
\end{equation}

\noindent The above definition takes advantage of the Lipschitz manifold 
structure of $\partial\Omega$. On such manifolds, smoothness spaces of index 
$s\in(0,1)$ can be defined in an intrinsic, invariant fashion by lifting 
their Euclidean counterparts onto the manifold itself via local charts. 
We shall, nonetheless, find it useful to consider higher order smoothness
spaces on $\partial\Omega$ in which case the above approach is no longer
effective. An alternative point of view has been developed 
by H.\,Whitney in {\bf\cite{Wh}} where he considered what amounts 
to higher order Lipschitz spaces on arbitrary closed sets. 
A far-reaching extension of this circle of ideas pertaining to the full 
scale of Besov and Sobolev spaces on irregular subsets of $\RR^n$ 
can be found in the book {\bf\cite{JW}} by A.\,Jonsson and H.\,Wallin.  

For the purpose of this introduction we note that one possible description of  
these higher order Besov spaces on the boundary of a Lipschitz domain 
$\Omega\subset\RR^n$ and for $m\in\NN$, $p\in(1,\infty)$, $s\in(0,1)$, 
reads

\begin{equation}\label{Bes-X}
\dot{B}^{m-1+s}_p(\partial\Omega)=\,\mbox{the closure of}\,\,
\Bigl\{(D^\alpha\,{\mathcal V}|_{\partial\Omega})_{|\alpha|\leq m-1}:\,
{\mathcal V}\in C^\infty_0(\RR^n)\Bigr\}\mbox{ in }
B_p^s(\partial\Omega)
\end{equation}

\noindent (we shall often make no notational distinction 
between a Banach space ${\mathfrak X}$ and 
${\mathfrak X}^N={\mathfrak X}\oplus\cdots\oplus{\mathfrak X}$ for a finite, 
positive integer $N$). A formal definition along with other equivalent 
characterizations of $\dot{B}^{m-1+s}_p(\partial\Omega)$ can be found 
in \S{6.4}. 

Given (\ref{nuk})-(\ref{W-Nr}), a necessary condition for the boundary data 
$\{g_k\}_{0\leq k\leq m-1}$ in (\ref{e0}) is that 

\begin{equation}\label{data-B}
\begin{array}{l}
\displaystyle{
\mbox{there exists 
$\dot{f}=\{f_\alpha\}_{|\alpha|\leq m-1}\in\dot{B}^{m-1+s}_p(\partial\Omega)$
such that}}
\\[10pt]
\displaystyle{
g_k=i^k\sum_{|\alpha|=k}\frac{k!}{\alpha!}\,\nu^\alpha\,f_\alpha,
\qquad\mbox{for each}\,\,\,\,0\leq k\leq m-1.}
\end{array}
\end{equation}

To state the (analytical and geometrical) conditions under which 
the problem (\ref{e0}), formulated as above, is well-posed, we need 
one final piece of terminology. By the {\it infinitesimal mean oscillation} 
of a function $F\in L_1(\Omega)$ we shall understand the quantity

\begin{equation}\label{e60}
\{F\}_{{\rm Osc}(\Omega)}:=\mathop{\hbox{lim\,sup}}_{\varepsilon\to 0} 
\left(\mathop{\hbox{sup}}_{{\{B_\varepsilon\}}_\Omega}
\meanint_{\!\!\!B_\varepsilon\cap\Omega}\,\,
\meanint_{\!\!\!B_\varepsilon\cap\Omega}\,
\Bigl|\,F(X)-F(Y)\,\Bigr|\,dXdY\right), 
\end{equation}

\noindent where $\{B_\varepsilon\}_\Omega$ stands for the set of arbitrary 
balls centered at points of $\Omega$ and of radius $\varepsilon$, and the 
barred integral is the mean value. In a similar fashion, the infinitesimal 
mean oscillation of a function $f\in L_1(\partial\Omega)$ is defined by

\begin{equation}\label{e61}
\{f\}_{{\rm Osc}(\partial\Omega)}
:=\mathop{\hbox{lim\,sup}}_{\varepsilon\to 0}
\left(\mathop{\hbox{sup}}_{\{B_\varepsilon\}_{\partial\Omega}}
\meanint_{\!\!\!B_\varepsilon\cap\partial\Omega}\,\,
\meanint_{\!\!\!B_\varepsilon\cap\partial\Omega}\,
\Bigl|\,f(X)-f(Y)\,\Bigr|\,d\sigma_Xd\sigma_Y\right), 
\end{equation}

\noindent where $\{B_\varepsilon\}_{\partial\Omega}$ is the collection 
of $n$-dimensional balls with centers on $\partial\Omega$ and of radius 
$\varepsilon$. 

\bigskip

Our main result reads as follows; see also Theorem~\ref{Theorem1}
for a more general version.

\vskip 0.08in
\begin{theorem}\label{Theorem}
In the above setting, for each $p\in (1,\infty)$ and $s\in (0,1)$, the 
problem {\rm (\ref{e0})} with boundary data as in (\ref{data-B}) has a 
unique solution ${\mathcal U}$ for which 
(\ref{W-Nr}) holds with $a=1-s-1/p$ provided the coefficient matrices 
$A_{\alpha\beta}$ and the exterior normal vector $\nu$ to $\partial\Omega$ 
satisfy

\begin{equation}\label{a0}
\{\nu\}_{{\rm Osc}(\partial\Omega)}
+\sum_{|\alpha|=|\beta|=m}\{ A_{\alpha\beta}\}_{{\rm Osc}(\Omega)}
\leq\,C\,s(1-s)\Bigl(pp'+s^{-1}(1-s)^{-1}\Bigr)^{-1}
\end{equation}

\noindent where $p'=p/(p-1)$ is the conjugate exponent of $p$. Above, $C$ 
is a sufficiently small constant which depends on $\kappa_0$, $\kappa_1$
and the Lipschitz constant of $\Omega$, and is independent of $p$ and $s$. 
Furthermore, the bound {\rm (\ref{a0})} can be improved for 
second order operators, i.e. when $m=1$, when the factor $s(1-s)$ 
in {\rm (\ref{a0})} can be removed.
\end{theorem}

\vskip 0.08in

Let ${\rm BMO}$ and ${\rm VMO}$ stand, respectively, for the John-Nirenberg 
space of functions of bounded mean oscillations and the Sarason space 
of functions of vanishing mean oscillations (considered either on $\Omega$ 
or on $\partial\Omega$). Since for an arbitrary function $F$ we have
(with the dependence on the domain dropped) 
$\{F\}_{{\rm Osc}}\leq 2\,{\rm dist}\,(F,{\rm VMO})$ where the 
distance is taken in ${\rm BMO}$, the smallness condition (\ref{a0}) in 
Theorem~\ref{Theorem} is satisfied if  

\begin{equation}\label{axxx}
{\rm dist}\,(\nu,{\rm VMO}\,(\partial\Omega))
+\sum_{|\alpha|=|\beta|=m}{\rm dist}\,(A_{\alpha\beta},{\rm VMO}\,(\Omega))
\leq\,C\,s(1-s)\Bigl(pp'+s^{-1}(1-s)^{-1}\Bigr)^{-1}. 
\end{equation}

\noindent In particular, this is trivially the case when 
$\nu\in {\rm VMO}(\partial\Omega)$ and the $A_{\alpha\beta}$'s belong 
to ${\rm VMO}(\Omega)$, irrespective of $p$, $s$, $\kappa_0$, $\kappa_1$ 
and the Lipschitz constant of $\Omega$. 

While the Lipschitz character of a domain $\Omega$ controls
the infinitesimal mean oscillation of its unit normal,
the inequality in the opposite direction is false in general, as seen 
by considering $\Omega:=\{(x,y)\in\RR^2:\,y>\varphi_\varepsilon(x)\}$ 
with $\varphi_\varepsilon(x):=x\,\sin\,(\varepsilon\log |x|^{-1})$. 
Indeed, a simple calculation gives 
$\|\varphi'_\varepsilon\|_{{\rm BMO}(\RR)}\leq C\varepsilon$, yet  
$\|\varphi'_\varepsilon\|_{L_\infty(\RR)}\sim 1$ uniformly 
for $\varepsilon\in(0,1/2)$. 

\medskip

An essentially equivalent reformulation of (\ref{e0}) is 

\begin{equation}\label{e0-bis}
\left\{
\begin{array}{l}
\displaystyle{\sum_{|\alpha|=|\beta|=m}
D^\alpha(A_{\alpha\beta}(X)\,D^\beta\,{\mathcal U})}=0
\qquad\mbox{in}\,\,\,\Omega,
\\[24pt] 
{\rm Tr}\,[D^\gamma\,{\mathcal U}]=g_\gamma
\,\,\quad\mbox{on}\,\,\partial\Omega,\qquad |\gamma|\leq m-1, 
\end{array}
\right.
\end{equation}
 
\noindent where ${\mathcal U}$ satisfies (\ref{W-Nr}) and 

\begin{equation}\label{data-G}
\dot{g}:=\{g_\gamma\}_{|\gamma|\leq m-1}\in\dot{B}^{m-1+s}_p(\partial\Omega),
\end{equation}

\noindent though an advantage of the classical formulation (\ref{e0}) 
is that the number of the data is minimal. 
For a domain $\Omega\subset\RR^2$ of class $C^r$, 
$r>{\textstyle\frac 1{2}}$, and for constant coefficient operators, 
the Dirichlet problem (\ref{e0-bis}) has been considered by S.\,Agmon 
in {\bf\cite{Ag1}} where he proved that there exists a unique solution 
${\mathcal U}\in C^{m-1+s}(\bar{\Omega})$, $0<s<r$, whenever 
$g_\gamma=D^\gamma{\mathcal V}|_{\partial\Omega}$ , $|\gamma|\leq m-1$,
for some function ${\mathcal V}\in C^{m-1+s}(\bar{\Omega})$.
See also {\bf\cite{Ag2}} for a related version. 

The innovation that allows us to consider, for the first time, 
boundary data in Besov spaces as in (\ref{data-B}) and (\ref{data-G}), is the 
systematic use of {\it weighted Sobolev spaces} such as those associated with 
the norm in (\ref{W-Nr}). In relation to the standard Besov scale in $\RR^n$,
we would like to point out that, thanks to Theorem~4.1 in {\bf\cite{JK}}
on the one hand, and Theorem~1.4.2.4 and Theorem~1.4.4.4 in {\bf\cite{Gr}}
on the other, we have

\begin{equation}\label{incls}
\begin{array}{l}
a=1-s-\frac{1}{p}\in (0,1-1/p)\Longrightarrow
W^{m,a}_p(\Omega)\hookrightarrow B^{m-1+s+1/p}_p(\Omega),
\\[10pt]
a=1-s-\frac{1}{p}\in (-1/p,0)\Longrightarrow
B^{m-1+s+1/p}_p(\Omega)\hookrightarrow W^{m,a}_p(\Omega).
\end{array}
\end{equation}

\noindent Of course, $W^{m,a}_p(\Omega)$ is just a classical Sobolev 
space when $a=0$. 

Remarkably, the classical trace theory for unweighted
Sobolev spaces turns out to have a most satisfactory analogue in this 
weighted context; for the upper half-space this has been worked out by 
S.V.\,Uspenski\u{\i} in {\bf\cite{Usp}}, a paper preceded by the significant 
work of E.\,Gagliardo in {\bf\cite{Ga}} in the unweighted case.  
As a consequence, we note that under the assumptions made in 
Theorem~\ref{Theorem}, 

\begin{equation}\label{Trace}
\sum_{|\alpha|\leq m-1}
\|{\rm Tr}\,[D^\alpha\,{\mathcal U}]\|_{B_p^{s}(\partial\Omega)}
\sim \left(\sum_{0\leq|\alpha|\leq m}\int_{\Omega}\rho(X)^{p(1-s)-1}\,
|D^\alpha{\mathcal U}(X)|^p\,dX\right)^{1/p},
\end{equation}

\noindent uniformly in ${\mathcal U}$ satisfying 
${\mathcal L}(X,D_X)\,{\mathcal U}=0$ in $\Omega$. The estimate (\ref{Trace}) 
can be viewed as a far-reaching generalization of a well-known 
characterization of the membership of a function to a Besov space in 
$\RR^{n-1}$ in terms of weighted Sobolev norm estimates for 
its harmonic extension to $\RR^n_+$ (see, e.g., Proposition~$7'$ 
on p.\,151 of {\bf\cite{St}}). 

Theorem~1.1 is new even in the case when $m=1$ and 
$A_{\alpha\beta}\in\CC^{l\times l}$ (i.e., for second order, constant 
coefficient systems) and provides a complete answer to the issue of 
well-posedness of the problem (\ref{e0}), (\ref{data-B}), (\ref{W-Nr}) in 
the sense that the small mean oscillation condition, depending on $p$ and $s$, 
is in the nature of best possible if one insists of allowing arbitrary 
indices $p$ and $s$ in (\ref{data-B}). This can be seen by considering 
the following Dirichlet problem for the Laplacian in a domain 
$\Omega\subset\RR^n$:

\begin{equation}\label{LapJK}
\left\{
\begin{array}{l}
\Delta\,{\mathcal U}=0\mbox{ in }\Omega,
\\[6pt]
{\rm Tr}\,{\mathcal U}=g\in B^s_p(\partial\Omega),
\\[6pt]
D^\alpha{\mathcal U}\in L_p(\Omega,\,\rho(X)^{p(1-s)-1}\,dX),\qquad
\forall\,\alpha\,:\,|\alpha|\leq 1. 
\end{array}
\right.
\end{equation}

\noindent It has long been known that, already in the case when 
$\partial\Omega$ exhibits one cone-like singularity, the well-posedness 
of (\ref{LapJK}) prevents the indices $(s,1/p)$ from taking arbitrary 
values in $(0,1)\times(0,1)$. At a more sophisticated level, the work 
of D.\,Jerison and C.\,Kenig in {\bf\cite{JK}} shows that (\ref{LapJK}) is 
well-posed in an arbitrary, given Lipschitz domain $\Omega$ if and only if 
the point $(s,1/p)$ belongs to a certain open region 
${\mathcal R}_{\Omega}\subseteq(0,1)\times(0,1)$, determined exclusively 
by the geometry of the domain $\Omega$ (cf. {\bf\cite{JK}} for more details). 
Let us also mention here that, even when $\partial\Omega$ is smooth and 
$m=l=1$, a well-known example due to N.\,Meyers (cf. {\bf\cite{Mey}}) shows 
that the well-posedness of (\ref{e0}) in the class of operators with bounded, 
measurable coefficients confines $p$ to a small neighborhood of $2$. 

\medskip 

Broadly speaking, there are two types of questions pertaining to the
well-posedness of the Dirichlet problem in a Lipschitz domain 
$\Omega$ for a divergence form, elliptic system (\ref{LOL}) of order $2m$
with boundary data in Besov spaces.  

\vskip 0.08in
\noindent {\it Question I.} Granted that the coefficients of ${\mathcal L}$ 
exhibit a certain amount of smoothness, identifying the Besov spaces 
for which this boundary value problem is well-posed. 

\vskip 0.08in
\noindent {\it Question II.} Alternatively, for a given Besov space 
characterize the class of Lipschitz domains $\Omega$ and elliptic 
operators ${\mathcal L}$ for which the aforementioned boundary value problem 
is well-posed. 

\vskip 0.08in 
\noindent These, as well as other related issues, have been a driving 
force behind many exciting, recent developments in partial differential 
equations and allied fields. Ample evidence of their impact can be found 
in C.\,Kenig's excellent account {\bf\cite{Ke}} which describes the state of 
the art in this field of research up to mid 1990's, with a particular 
emphasis on the role played by harmonic analysis techniques. 

One generic problem which falls under the scope of {\it Question I} 
is to determine the optimal scale of spaces on which the Dirichlet 
problem for an elliptic system of order $2m$ is solvable in an 
{\it arbitrary Lipschitz domain} $\Omega$ in $\RR^n$. 
The most basic case, that of the constant coefficient Laplacian in 
arbitrary Lipschitz domains in $\RR^n$, is now well-understood thanks 
to the work of B.\,Dahlberg and C.\,Kenig {\bf\cite{DK}}, in the case 
of $L_p$-data,
and D.\,Jerison and C.\,Kenig {\bf\cite{JK}}, in the case of Besov data. 
The case of (\ref{LapJK}) for boundary data exhibiting higher regularity 
(i.e., $s>1$) has been recently dealt with by V.\,Maz'ya and T.\,Shaposhnikova
in {\bf\cite{MS2}} where nearly optimal smoothness conditions for 
$\partial\Omega$
are found in terms of the properties of $\nu$ as a Sobolev space multiplier.  
Generalizations of (\ref{LapJK}) to the case of variable-coefficient, single, 
second order elliptic equations have been obtained 
by M.\,Mitrea and M.\,Taylor in {\bf\cite{MT1}}, {\bf\cite{MT2}}, 
{\bf\cite{MT3}}. 

In spite of substantial progress in recent years, there remain many basic 
open questions, particularly for $l>1$ and/or $m>1$ (corresponding to 
genuine systems and/or higher order equations), even in the case of 
{\it constant coefficient} operators in Lipschitz domains. In this context, 
one significant problem (as mentioned in, e.g., {\bf\cite{Fa}}) is to 
determine the sharp range of $p$'s for which the Dirichlet problem for 
elliptic systems with $L_p$-boundary data is well-posed. 
In {\bf\cite{PV}}, J.\,Pipher and G.\,Verchota have developed a $L_p$-theory 
for real, constant coefficient, higher order systems 
$L=\sum_{|\alpha|=2m}A_\alpha D^\alpha$ when $p$ is near $2$, i.e. 
$2-\varepsilon<p<2+\varepsilon$ with $\varepsilon>0$ depending on the 
Lipschitz character of $\Omega$ but this range is not optimal.
Recently, more progress for the biharmonic equation and for general constant 
coefficient, second order systems with real coefficients, which are elliptic 
in the sense of Legendre-Hadamard was made by Z.\,Shen in {\bf\cite{Sh}}, 
where he further extended the range of $p$'s from 
$(2-\varepsilon,2+\varepsilon)$ to 
$(2-\varepsilon,\frac{2(n-1)}{n-3}+\varepsilon)$ for a general Lipschitz 
domain $\Omega$ in $\RR^n$, $n\geq 4$, where as before 
$\varepsilon=\varepsilon(\partial\Omega)>0$. Let us also mention
here the work {\bf\cite{AP}} of V.\,Adolfsson and J.\,Pipher who have dealt 
with the Dirichlet problem for the biharmonic operator in arbitrary 
Lipschitz domains and with data in Besov spaces, {\bf\cite{Ve}} 
where G.\,Verchota formulates and solves a Neumann-type problem for the
bi-Laplacian in Lipschitz domains and with boundary data in $L_2$, 
{\bf\cite{MMT}} where the authors treat the Dirichlet 
problem for variable coefficient symmetric, real, elliptic systems of 
second order in an arbitrary Lipschitz domain $\Omega$ and with boundary 
data in $B^s_p(\partial\Omega)$, when $2-\varepsilon<p<2+\varepsilon$ and 
$0<s<1$, as well as the paper {\bf\cite{KM}} by V.\,Kozlov and V.\,Maz'ya,
which contains an explicit description of the asymptotic behavior
of null-solutions of constant coefficient, higher order, elliptic operators
near points on the boundary of a domain with a sufficiently small Lipschitz
constant. 

A successful strategy for dealing with {\it Question II} consists of  
formulating and solving the analogue of the original problem in a 
standard case, typically when $\Omega=\RR^n_+$ and ${\mathcal L}$ has 
constant coefficients, and then deviating from this most standard setting 
by allowing perturbations of a certain magnitude. 
A paradigm result in this regard, going back to the work of Agmon, Douglis, 
Nirenberg and Solonnikov in the 50's and 60's is that the Dirichlet problem 
is solvable in the context of Sobolev-Besov spaces if $\partial\Omega$ 
is sufficiently smooth and if ${\mathcal L}$ has continuous coefficients. 
The latter requirement is an artifact of the method of proof (based 
on Korn's trick of freezing the coefficients) which requires measuring 
the size of the oscillations of the coefficients in a {\it pointwise sense} 
(as opposed to integral sense, as in (\ref{e60})). 
For a version of {\it Question II}, corresponding to boundary data of 
higher regularity, optimal results have been obtained by V.\,Maz'ya and 
T.\,Shaposhnikova in {\bf\cite{MS}}. In this context, the natural language 
for describing the smoothness of the domain $\Omega$ is that of Sobolev 
space multipliers. 

While the study of boundary value problems in a domain $\Omega\subset\RR^n$ 
for elliptic differential operators with discontinuous coefficients 
goes a long way back (for instance, C.\,Miranda has considered in 
{\bf\cite{Mir}} operators with coefficients in the Sobolev space $W^1_n$), 
a lot of attention has been devoted lately to the class of operators with 
coefficients in ${\rm VMO}$ (it is worth pointing out here that 
$W^1_n\hookrightarrow{\rm VMO}$ on Lipschitz subdomains of $\RR^n$). 
Much of the impetus for the recent surge of interest in this particular 
line of work stems from an observation made by F.\,Chiarenza, M.\,Frasca 
and P.\,Longo in the early 1990's. 
More specifically, while investigating interior estimates for 
the solution of a scalar, second-order elliptic differential equation 
of the form ${\mathcal L}\,{\mathcal U}=F$, these authors have noticed 
in {\bf\cite{CFL1}} that ${\mathcal U}$ can be related to $F$ via a potential 
theoretic representation formula in which the residual terms 
are commutators between operators of Calder\'on-Zygmund type, on the one 
hand, and operators of multiplication by the coefficients of ${\mathcal L}$, 
on the other hand. This made it possible to control these terms by invoking 
the commutator estimate of Coifman-Rochberg-Weiss ({\bf\cite{CRW}}).  
Various partial extensions of this result can be found in 
{\bf\cite{AQ}}, {\bf\cite{By1}}, {\bf\cite{CaPe}}, 
{\bf\cite{CFL2}}, {\bf\cite{Faz}}, {\bf\cite{Gu}}, {\bf\cite{Ra}}, 
and the references therein. Here we would just like to mention that,  
in the whole Euclidean space, a different approach (based on estimates for 
the Riesz transforms) has been devised by T.\,Iwaniec and C.\,Sbordone in 
{\bf\cite{IS}}. 

Compared to the aforementioned works, our approach is more akin to that  
of F.\,Chiarenza and collaborators ({\bf\cite{CFL1}}, {\bf\cite{CFL2}}), 
though there are fundamental differences between solving boundary problems 
for higher order and for second order operators. 
One difficulty inherently linked with the case $m>1$ arises from the 
way the norm in (\ref{W-Nr}) behaves under a change of variables
$\varkappa:\Omega=\{(X',X_n):\,X_n>\varphi(X')\}\to\RR^n_+$ designed to 
flatten the Lipschitz surface $\partial\Omega$. When $m=1$, a simple 
bi-Lipschitz changes of variables such as the inverse of the map 
$\RR^n_+\ni (X',X_n)\mapsto(X',\varphi(X')+X_n)\in\Omega$ will do, but 
matters are considerable more subtle in the case $m>1$. In this latter 
situation, we employ a special global flattening map first introduced by 
J.\,Ne\v{c}as (in a different context; cf. p.\,188 in {\bf\cite{Nec}}) 
and then independently rediscovered and/or further adapted to new 
 settings by several authors, including V.\,Maz'ya and 
T.\,Shaposhnikova in {\bf\cite{MS}}, 
B.\,Dahlberg, C.\,Kenig J.\,Pipher, E.\,Stein and G.\,Verchota 
(cf. {\bf\cite{Dah}} and the discussion in {\bf\cite{DKPV}}), 
and S.\,Hofmann and J.\,Lewis in {\bf\cite{HL}}. Our main novel contribution 
in this regard is adapting this circle of ideas to the context when one 
seeks pointwise estimates for higher order derivatives of $\varkappa$ 
and $\lambda:=\varkappa^{-1}$ in terms of 
$[\nabla\varphi]_{{\rm BMO}(\RR^{n-1})}$. 

Another ingredient of independent interest is deriving estimates for 
$D_x^\alpha D_y^\beta G(x,y)$ where $G$ is the Green function associated 
with a constant (complex) coefficient system $L(D)$ of order $2m$ in the 
upper half space, which are sufficiently well-suited for deriving commutator 
estimates in the spirit of {\bf\cite{CRW}}. The methods employed in earlier 
work are largely based on explicit representation formulas for $G(x,y)$ and, 
hence, cannot be adapted easily to the case of non-symmetric, complex 
coefficient, higher order systems. By way of contrast, our approach 
consists of proving directly that the residual part 
$R(x,y):=G(x,y)-\Phi(x-y)$, where $\Phi$ is a fundamental solution 
for $L(D)$, has the property that $D_x^\alpha D_y^\beta R(x,y)$ 
is a Hardy-type kernel whenever $|\alpha|=|\beta|=m$. 

The layout of the paper is as follows. Section~2 contains estimates 
for the Green function in the upper-half space. Section~3 deals with 
integral operators (of Calder\'on-Zygmund and Hardy type) as well as
commutator estimates on weighted Lebesgue spaces. In the last part 
of this section we also revisit Gagliardo's extension operator and
establish estimates in the context of ${\rm BMO}$. 
Section~4 contains a discussion of the Dirichlet problem for higher 
order, variable coefficient elliptic systems in the upper-half space.
Then the adjustments necessary to treat the case of an unbounded domain 
lying above the graph of a Lipschitz function are presented in Section~5.
Finally, in Section~6, we explain how to handle the case of a bounded
Lipschitz domain, and state and prove Theorem~\ref{Theorem1}
(from which Theorem~\ref{Theorem} follows). This section also contains 
further complements and extensions of our main result.

\bigskip

\section{Green's matrix estimates in the half-space}
\setcounter{equation}{0}

In this section we prove a key estimate for derivatives of Green's matrix 
associated with the Dirichlet problem for homogeneous, higher-order 
constant coefficient elliptic systems in the half-space $\mathbb{R}^n_+$.

\subsection{Statement of the main result}

Let ${L}(D_x)$ be a matrix-valued differential operator

\begin{equation}\label{eq1.1}
{L}(D_x)=\sum_{|\alpha|=2m}A_{\alpha} D^{\alpha}_x,
\end{equation}
\noindent

\noindent where the $A_\alpha$'s are constant $l\times l$ matrices with
complex entries. Throughout the paper, $D^\alpha_x:= i^{-|\alpha|}
\partial_{x_1}^{\alpha_1}\partial_{x_2}^{\alpha_2}\cdots
\partial_{x_n}^{\alpha_n}$ if
$\alpha=(\alpha_1,\alpha_2,...,\alpha_n)\in\NN_0^n$.
Here and elsewhere, $\NN$ stands for the collection of all
positive integers and $\NN_0:=\NN\cup\{0\}$. 

We assume that ${L}$ is strongly elliptic, i.e. there exists $\kappa>0$ such 
that $\sum_{|\alpha|=m}\|A_\alpha\|_{\CC^{l\times l}}\leq\kappa^{-1}$ and

\begin{equation}\label{eq1.2}
\Re\,\langle{L}(\xi)\eta,\eta\rangle_{\mathbb{C}^l}\geq\kappa\,
|\xi|^{2m}\,\|\eta\|^2_{\mathbb{C}^l},\qquad 
\forall\,\xi\in\RR^n,\,\,\,\forall\,\eta\in\mathbb{C}^l.
\end{equation}

\noindent In what follows, in order to simplify notations, we shall denote the 
norms in different finite-dimensional real Euclidean spaces by $|\cdot|$ 
irrespective of their dimensions. Also, quite frequently, we shall make 
no notational distinction between a space of scalar functions, call 
it ${\mathfrak X}$, and the space of vector-valued functions 
(of a fixed, finite dimension) whose components are in ${\mathfrak X}$.

We denote by $F(x)$ a  fundamental matrix of the  operator
${L}(D_x)$, i.e. a $l\times l$ matrix solution of the system

\begin{equation}\label{eq1.3}
{L}(D_x)F(x)=\delta(x)I_l\quad\mbox{in}\,\,\mathbb{R}^n,
\end{equation}

\noindent where $I_l$ is the $l\times l$ identity matrix and $\delta$ is the
Dirac function. 

Consider the Dirichlet problem

\begin{equation}\label{eq1.5}
\left\{
\begin{array}{l}
{L}(D_x)u=f\qquad\qquad\qquad \mbox{in}\,\,\mathbb{R}^n_+,
\\[6pt] 
{\rm Tr}\,[\partial^j u/\partial x_n^j]=f_j\quad j=0,1,\ldots,m-1,   
\qquad\quad\,\,\mbox{on}\,\,\mathbb{R}^{n-1}
\end{array}
\right.
\end{equation}

\noindent where 
$\mathbb{R}^n_+:=\{x=(x',x_n):\,x'\in\mathbb{R}^{n-1},\,x_n>0\}$
and ${\rm Tr}$ is the boundary trace operator. 
Hereafter, we shall identify $\partial\RR^n_+$ with $\RR^{n-1}$ in 
a canonical fashion.

For each $y'\in \mathbb{R}^{n-1}$ we introduce the Poisson matrices 
$P_0,\ldots, P_{m-1}$ for problem (\ref{eq1.5}), i.e. the solutions of 
the boundary-value problems

\begin{equation}\label{eq1.6}
\left\{
\begin{array}{l}
{ L}(D_x)P_j(x,y')= 0\,I_l
\qquad\qquad\qquad\mbox{on}\,\,\mathbb{R}^n_+,
\\[10pt]
\displaystyle{\left(\frac{\partial^k}{\partial x_n^k}P_j\right)
(\,(x',0),y'\,)}
=\delta_{jk}\,\delta(x'-y')I_l\,\,\,{\rm for}\,\,\,x'\in\mathbb{R}^{n-1},\,\,
0\leq k\leq m-1,
\end{array}
\right.
\end{equation}

\noindent where $\delta_{jk}$ is the usual Kronecker symbol and 
$0\leq j\leq m-1$. The matrix-valued function $P_j(x,0')$ is positive 
homogeneous of degree $j+1-n$, i.e.

\begin{equation}\label{eq1.7}
P_j(x,0')=|x|^{j+1-n}\,P_j(x/|x|,0'),\qquad x\in\mathbb{R}^n,
\end{equation}

\noindent where $0'$ denotes the origin of $\mathbb{R}^{n-1}$.
The restriction of $P_j(\cdot,0')$ to the upper half-sphere
$S^{n-1}_+$ is smooth and vanishes on the equator along with all of
its derivatives up to order $m-1$ (see for example, \S{10.3} in 
{\bf\cite{KMR2}}). Hence,

\begin{equation}\label{eq1.8}
\|P_j(x,0')\|_{\mathbb{C}^{l\times l}}
\leq C\,\frac{x_n^m}{|x|^{n+m-1-j}},\qquad x\in\RR^n_+,
\end{equation}

\noindent and, consequently,

\begin{equation}\label{eq1.9}
\|P_j(x,y')\|_{\mathbb{C}^{l\times l}}
\leq C\,\frac{x_n^m}{|x-(y',0)|^{n+m-1-j}},\qquad
x\in\RR^n_+,\,\,\,\,y'\in\RR^{n-1}.
\end{equation}

By $G(x,y)$ we shall denote the Green's matrix of the problem (\ref{eq1.5}), 
i.e. the unique solution of the boundary-value problem

\begin{equation}\label{eq1.10}
\left\{
\begin{array}{l}
{L}(D_x)G(x,y)=\delta(x-y)I_l\quad\mbox{for}\,\,x\in\mathbb{R}^n,
\\[6pt]
\displaystyle{\left(\frac{\partial ^j}{\partial x_n^j}G\right)((x',0),y)}
=0\,I_l\qquad
\mbox{for}\,\,x'\in\mathbb{R}^{n-1}, \,\,\, 0\leq j\leq m-1,
\end{array}
\right.
\end{equation}

\noindent where $y\in\mathbb{R}^n_+$ is regarded as a parameter.

We now introduce the matrix 

\begin{equation}\label{defRRR}
R(x,y):=F(x-y)-G(x,y),\qquad x,y\in\RR^n_+,
\end{equation}

\noindent so that, for each fixed $y\in\mathbb{R}^n_+$,

\begin{equation}\label{eq1.12}
\left\{
\begin{array}{l}
{L}(D_x)\,R(x,y)=0\qquad\qquad\qquad\qquad\qquad\quad\,\,\,\,
\mbox{for}\,\,x\in\mathbb{R}^n,
\\[6pt]
\displaystyle{\left(\frac{\partial^j}{\partial x_n^j}R\right)((x',0),y)}
=\left(\frac{\partial^j}{\partial x_n^j}F\right)((x',0)-y)
\quad\mbox{for}\,\,x'\in\mathbb{R}^{n-1},\,\, 0\leq j\leq m-1.
\end{array}
\right.
\end{equation}

\noindent Our goal is to prove the following result.

\begin{theorem}\label{th1}
For all multi-indices $\alpha,\beta$ of length $m$

\begin{equation}\label{mainest}
\|D^\alpha_xD^\beta_y R(x,y)\|_{\mathbb{C}^{l\times l}}
\leq C\,|x-\bar{y}|^{-n},
\end{equation}

\noindent for $x,y\in\RR^n_+$, where $\bar{y}:=(y',-y_n)$ is
the reflection of the point $y\in\RR^n_+$ with respect to $\partial\RR^n_+$.
\end{theorem}

In the proof of Theorem~\ref{th1} we distinguish two cases, $n>2m$ and 
$n\leq 2m$, which we shall treat separately. Our argument pertaining to
the situation when $n>2m$ is based on a lemma to be proved in the 
subsection below.

\subsection{Estimate for a parameter dependent integral}

As a preamble to the proof of Theorem~\ref{th1}, here we dispense with the
following technical result. 

\begin{lemma}\label{lem1}
Let $a$ and $b$ be two non-negative numbers and assume that 
$\zeta\in \mathbb{R}^N$. Then for every $\varepsilon>0$ and $0<\delta<N$
there exists a constant $c(N,\varepsilon,\delta)>0$ such that 

\begin{equation}\label{E1}
\int_{\mathbb{R}^N}
\frac{d\eta}{(|\eta|+a)^{N+\varepsilon}(|\eta-\zeta|+ b)^{N-\delta}} 
\leq\frac{c(N,\varepsilon,\delta)}{a^\varepsilon(|\zeta|+a+b)^{N-\delta}}.
\end{equation}
\end{lemma}

\noindent{\bf Proof.} Write ${\mathcal J}={\mathcal J}_1+{\mathcal J}_2$ 
where ${\mathcal J}$ stands for the integral in the left side of (\ref{E1}), 
whereas  ${\mathcal J}_1$ and ${\mathcal J}_2$ denote the integrals
obtained by splitting the domain of integration in ${\mathcal J}$ into 
$B_a=\{\eta\in \mathbb{R}^N:\,|\eta|<a\}$ and $\RR^n\setminus B_a$, 
respectively. If $|\zeta|<2a$, then

\begin{equation}\label{I1}
{\mathcal J}_1\leq a^{-N-\varepsilon}
\int_{B_a}\frac{d\eta}{(|\eta-\zeta|+b)^{N-\delta}}
\leq c\,a^{-N-\varepsilon}\int_{B_{4a}}\frac{d\xi}{(|\xi|+b)^{N-\delta}}.
\end{equation}

\noindent Hence

\begin{equation}\label{II1}
{\mathcal J}_1\leq 
\left\{
\begin{array}{l}
c\,a^{-N-\varepsilon}\,a^{N}/b^{N-\delta}\qquad{\rm if}\,\,a<b,
\\[6pt]
c\,a^{-N-\varepsilon+\delta}\qquad\quad{\rm if}\,\,a>b,
\end{array}
\right.
\end{equation}

\noindent so that, in particular, 

\begin{equation}\label{I1est}
|\zeta|<2a\Longrightarrow 
{\mathcal J}_1\leq c\,a^{-\varepsilon}(|\zeta|+a+b)^{\delta-N}.
\end{equation}

Let us now assume that $|\zeta|>2a$. Then

\begin{equation}\label{I1est2}
{\mathcal J}_1\leq\int_{B_a}\frac{d\eta}{(|\eta|+a)^{N+\varepsilon}}\, 
\frac{c}{(|\zeta|+b)^{N-\delta}}
\leq c\,a^{-\varepsilon}(|\zeta|+a+b)^{\delta-N},
\end{equation}

\noindent which is of the right order. As for ${\mathcal J}_2$, we write 

\begin{equation}\label{I2}
{\mathcal J}_2\leq\int_{\mathbb{R}^n\backslash B_a}
\frac{d\eta}{|\eta|^{N+\varepsilon}(|\eta-\zeta|+b)^{N-\delta}} 
={\mathcal J}_{2,1}+{\mathcal J}_{2,2}.
\end{equation}

\noindent where ${\mathcal J}_{2,1}$, ${\mathcal J}_{2,2}$ are obtained by 
splitting the domain of integration in the above integral into the set 
$\{\eta:|\eta|>\max\{a,2|\zeta|\}\}$ and its complement in 
$\mathbb{R}^n\backslash B_a$. We have

\begin{eqnarray}\label{Eq4}
{\mathcal J}_{2,1} & \leq & \int_{|\eta|>\max\{a,b,2|\zeta|\}}
\frac{d\eta}{|\eta|^{N+\varepsilon}(|\eta|+b)^{N-\delta}} 
+\int_{b>|\eta|>\max\{a,b,2|\zeta|\}}
\frac{d\eta}{|\eta|^{N+\varepsilon}(|\eta|+b)^{N-\delta}}
\nonumber\\[6pt]
& \leq & c\Biggl(\int_{|\eta|>\max\{a,b,2|\zeta|\}}
\frac{d\eta}{|\eta|^{2N+\varepsilon-\delta}}
+\frac{1}{b^{N-\delta}}\int_{b>|\eta|>\max\{a,b,2|\zeta|\}}
\frac{d\eta}{|\eta|^{N+\varepsilon}}\Biggr)
\nonumber\\[6pt]
& \leq & \frac{c}{(a+b+|\zeta|)^{N+\varepsilon-\delta}}
+\frac{c}{a^\varepsilon(a+b+|\zeta|)^{N-\delta}}
\nonumber\\[6pt]
& \leq & \frac{c}{a^\varepsilon(a+b+|\zeta|)^{N-\delta}}.
\end{eqnarray}

There remains to estimate the integral

\begin{equation}\label{Eq5}
{\mathcal J}_{2,2}=\int_{B_{2|\zeta|}\backslash B_a}
\frac{d\eta}{|\eta|^{N+\varepsilon}(|\eta-\zeta|+b)^{N-\delta}} 
={\mathcal J}_{2,2}^{(1)}+{\mathcal J}_{2,2}^{(2)},
\end{equation}

\noindent where ${\mathcal J}_{2,2}^{(1)}$ and ${\mathcal J}_{2,2}^{(2)}$ are 
obtained by splitting the domain of integration in ${\mathcal J}_{2,2}$ into 
$B_{|\zeta|/2}\backslash B_a$ and its complement (relative to 
$B_{2|\zeta|}\backslash B_a$). On the one hand, 

\begin{equation}\label{Eq6}
{\mathcal J}_{2,2}^{(1)}\leq\frac{c}{(|\zeta|+b)^{N-\delta}}
\int_{B_{|\zeta|/2}\backslash B_a}\frac{d\eta}{|\eta|^{N+\varepsilon}} 
\leq\frac{c}{a^\varepsilon(|\zeta|+a+b)^{N-\delta}}.
\end{equation}

\noindent On the other hand, whenever $|\zeta|>a/2$, the integral 
${\mathcal J}_{2,2}^{(2)}$, which extends over all $\eta$'s such 
that $|\eta|>a$, $2|\zeta|>|\eta|>|\zeta|/2$, can be estimated as

\begin{eqnarray*}
{\mathcal J}_{2,2}^{(2)} & \leq & \frac{c}{|\zeta|^{N+\varepsilon}}
\int_{B_{2|\zeta|}\backslash B_a}\frac{d\eta}{(|\eta-\zeta|+b)^{N-\delta}}
\leq\frac{c}{|\zeta|^{N+\varepsilon}}
\int_{B_{4|\zeta|}}\frac{d\xi}{(|\xi|+b)^{N-\delta}}
\\[6pt]
& \leq & \frac{c}{|\zeta|^{N+\varepsilon}}
\Biggl(\int_{{|\xi|<4|\zeta|}\atop{|\xi|<b}}
\frac{d\xi}{(|\xi|+b)^{N-\delta}}
+\int_{{|\xi|<4|\zeta|}\atop{|\xi|>b}}\frac{d\xi}{(|\xi|+b)^{N-\delta}}\Biggr).
\end{eqnarray*}

\noindent Consequently, 

\begin{equation}\label{J22}
{\mathcal J}_{2,2}^{(2)}\leq
\frac{c\,\min\{|\zeta|,b\}^N}{|\zeta|^{N+\varepsilon} b^{N-\delta}}.
\end{equation}

\noindent Using $|\zeta|>a/2$ and the obvious inequality

\begin{equation}\label{trivial}
\min\{|\zeta|,b\}^N\,\max\{|\zeta|,b^{N-\delta}\}\leq|\zeta|^N\,b^{N-\delta}
\end{equation}

\noindent we arrive at

\begin{equation}\label{Eq7}
{\mathcal J}_{2,2}^{(2)}\leq c\,a^{-\varepsilon}(|\zeta|+a+b)^{\delta-N}.
\end{equation}

\noindent The estimate (\ref{Eq7}), along with (\ref{Eq6}) and (\ref{Eq5}), 
gives the upper bound $c\,a^{-\varepsilon}(|\zeta|+ a+b)^{\delta-N}$ for 
${\mathcal J}_{2,2}$. Combining this with (\ref{Eq4}) we obtain the same 
majorant for ${\mathcal J}_{2}$ which, together with a similar result 
for ${\mathcal J}_{1}$ already obtained leads to (\ref{E1}). 
The proof of the lemma is therefore complete.
\hfill$\Box$
\vskip 0.08in

\subsection{Proof of Theorem~\ref{th1} for $n>2m$}

In the case when $n>2m$ there exists a unique fundamental matrix $F(x)$ 
for the the operator (\ref{eq1.1}) which is positive homogeneous of 
degree $2m-n$. We shall use the integral representation formula 

\begin{equation}\label{IntRRR}
R(x,y)=R_0(x,y)+\ldots+ R_{m-1}(x,y),\qquad x,y\in\RR^n_+,
\end{equation}

\noindent where $R(x,y)$ has been introduced in (\ref{defRRR}) and, 
with $P_j$ as in (\ref{eq1.6}), we set  

\begin{equation}\label{eq1.13}
R_j(x,y):=\int_{\mathbb{R}^{n-1}}P_j(x,\xi')\,
\left(\frac{\partial^j}{\partial x_n^j}F\right)((\xi',0)-y)\,d\xi',
\qquad 0\leq j\leq m-1.  
\end{equation}

\noindent Then, thanks to (\ref{eq1.8}) we have

\begin{equation}\label{eq1.14}
\|R_j(x,y)\|_{\mathbb{C}^{l\times l}}\leq C\,\int_{\mathbb{R}^{n-1}}
\frac{x_n^m}{|x-(\xi',0)|^{n+m-1-j}}\cdot\frac{d\xi'}{|(\xi',0)-y|^{n-2m+j}}.
\end{equation}

Next, putting

\begin{eqnarray*}
N=n-1 &, & \quad a=x_n,\\
\varepsilon=m-j &, & \quad b=y_n,\\
\delta=2m-j-1 &, & \quad \zeta=y'-x',
\end{eqnarray*}

\noindent in the formulation of Lemma~\ref{lem1}, we obtain from 
(\ref{eq1.14})

\begin{equation}\label{Rj}
\|R_j(x,y)\|_{\mathbb{C}^{l\times l}}
\leq\frac{C\,x_n^j}{(|y'-x'|+x_n+y_n)^{n-2m+j}},\qquad 0\leq j\leq m-1.
\end{equation}

\noindent Summing up over $j=0,\ldots,m-1$ gives, by virtue of (\ref{IntRRR}), 
the estimate

\begin{equation}\label{eq1.16}
\|R(x,y)\|_{\mathbb{C}^{l\times l}}\leq C\,|x-{\bar y}|^{2m-n},
\qquad x,y\in\RR^n_+.
\end{equation}

In order to obtain pointwise estimates for derivatives of $R(x,y)$, we make 
use of the following local estimate for a solution of problem (\ref{eq1.5}) 
with $f=0$. Recall that $W^s_p$ stands for the classical $L_p$-based Sobolev 
space of order $s$. The particle {\it loc} is used to brand the local versions 
of these (and other similar) spaces.

\begin{lemma}\label{l1.1}{\rm [see {\bf\cite{ADN}}]}
Let $\zeta$ and $\zeta_0$ be functions in $C_0^\infty(\mathbb{R}^n)$ 
such that $\zeta_0 =1$ in a neighborhood of $\mbox{supp}\,\zeta$. 
Then the solution $u\in W_2^m(\mathbb{R}_+^n,loc)$ of problem (\ref{eq1.5}) 
with $f=0$ and $\varphi_j\in W_p^{k+1-j-1/p}(\mathbb{R}^{n-1},loc)$, 
where $k\geq m$ and $p\in(1,\infty)$, belongs to 
$W_p^{k+1}(\mathbb{R}^n_+,loc)$ and satisfies the estimate

\begin{equation}\label{eq1.17}
\|\zeta u\|_{W_p^{k+1}(\mathbb{R}^n_+)}\leq C\,
\Bigl(\sum_{j=0}^{m-1}\|\zeta_0 \varphi_j\|_{W_p^{k+1-j-1/p}(\mathbb{R}^{n-1})}
+\|\zeta_0 u\|_{L_p(\mathbb{R}^n_+)}\Bigr),
\end{equation}

\noindent where $C$ is independent of $u$ and $\varphi_j$.
\end{lemma}

\noindent Let $B(x,r)$ denote the ball of radius $r>0$ centered at $x$.

\begin{corollary}\label{c1.2}
Assume that $u\in W_2^m(\mathbb{R}_+^n,loc)$ is a solution of 
problem (\ref{eq1.5}) 
with $f=0$ and $\varphi_j\in C^{k+1-j}(\mathbb{R}^{n-1},loc)$. Then for any
$z\in\overline{\mathbb{R}^n_+}$ and $\rho>0$

\begin{equation}\label{eq1.18}
\sup_{\mathbb{R}^n_+\cap B(z,\rho)}|\nabla _k u|\leq C\,
\Bigl(\,\rho^{-k}\sup_{\mathbb{R}^n_+\cap B(z,2\rho)}|u|
+\sum_{j=0}^{m-1}\sum_{s=0}^{k+1-j}
\rho^{s+j-k}\sup_{\mathbb{R}^{n-1}\cap B(z,2\rho)}|\nabla'_{s}\varphi_j|\Bigr),
\end{equation}

\noindent where $\nabla'_s$ is the gradient of order $s$ in $\mathbb{R}^{n-1}$.
Here $C$ is a  constant independent of $\rho$, $z$, $u$ and $f_j$.
\end{corollary}

\noindent{\bf Proof.} Given the dilation invariant nature of the estimate
we seek, it suffices to assume that $\rho=1$. 
Given $\phi\in C^{k+1-j}(\RR^{n-1})$ supported in $\RR^{n-1}\cap B(z,2)$,
we observe that, for a suitable $\theta\in (0,1)$,

\begin{eqnarray}\label{simple}
\|\phi\|_{W_p^{k+1-j-1/p}(\mathbb{R}^{n-1})}
& \leq C & \|\phi\|_{L_p(\mathbb{R}^{n-1})}^\theta
\|\phi\|_{W_p^{k+1-j}(\mathbb{R}^{n-1})}^{1-\theta}
\nonumber\\[6pt]
& \leq C &
\sum_{s=0}^{k+1-j}\sup_{\mathbb{R}^{n-1}\cap B(z,2)}|\nabla'_{s}\phi|.
\end{eqnarray}

\noindent Also, if $p>n$, 

\begin{equation}\label{eq1.19}
\sup\limits_{\mathbb{R}^n_+}|\nabla_k v|\leq C\,
\|v\|_{W_p^{k+1}(\mathbb{R}^n_+)},
\end{equation}

\noindent by virtue of the classical Sobolev inequality. 
Combining (\ref{simple}), (\ref{eq1.19}) with Lemma~\ref{l1.1} now readily
gives (\ref{eq1.18}).
\hfill$\Box$
\vskip 0.08in

Given $x,y\in\RR^n_+$, set $\rho:=|x-\bar{y}|/5$ and pick
$z\in\partial\RR^n_+$ such that $|x-z|=\rho/2$. It follows that for any
$w\in\RR^n_+\cap B(z,2\rho)$ we have
$|x-\bar{y}|\leq |x-z|+|z-w|+|w-\bar{y}|\leq \rho/2+2\rho+|w-\bar{y}|
\leq |x-\bar{y}|/2+|w-\bar{y}|$. Consequently,
$|x-\bar{y}|/2\leq |w-\bar{y}|$  for every $w\in\RR^n_+\cap B(z,2\rho)$,
so that, ultimately,

\begin{equation}\label{nablaPhi}
\rho^{\nu-k}\sup_{w\in\mathbb{R}^{n-1}\cap B(z,2\rho)}
\|\nabla'_{\nu}F(w-y)\|_{\CC^{l\times l}}
\leq \frac{C}{|x-\bar{y}|^{n-2m+k}},
\end{equation}

\noindent for each $\nu\in\NN_0$.
Granted (\ref{eq1.16}) and our choice of $\rho$, we altogether obtain that

\begin{equation}\label{eq1.20}
\|D^\alpha_x R(x,y)\|_{\mathbb{C}^{l\times l}}
\leq C_{k}\,|x-\bar{y}|^{2m-n-k},\qquad x,y\in\RR^n_+,
\end{equation}

\noindent for each multi-index $\alpha\in\NN_0^n$ of length $k$.

In the following two formulas, it will be convenient to use the notation 
$R_{\mathcal L}$ for the matrix $R$ associated with the operator 
${\mathcal L}(D_x)$ as in (\ref{defRRR}). By Green's formula

\begin{equation}\label{eq1.11}
R_{\mathcal L}(y,x)=\Bigl[R_{{\mathcal L}^*}(x,y)\Bigr]^*,\qquad x,y\in\RR^n_+,
\end{equation}

\noindent where the superscript star indicates adjunction.

In order to estimate {\it mixed} partial derivatives, we observe 
that (\ref{eq1.11}) entails

\begin{equation}\label{eq1.21}
(D^\beta_y R_{\mathcal L})(x,y)
=\Bigl[(D^\beta_x R_{{\mathcal L}^*})(y,x)\Bigr]^*
\end{equation}

\noindent and remark that ${\mathcal L}^*$ has properties similar 
to ${\mathcal L}$. This, in concert with (\ref{eq1.20}) and 

\begin{equation}\label{reflect}
|x-\bar{y}|=|\bar{x}-y|,\qquad x, y\in\mathbb{R}^n_+.
\end{equation}
 
\noindent yields

\begin{equation}\label{eq1.22}
\|D^\beta_y R(x,y)\|_{\mathbb{C}^{l\times l}}
\leq C_{\beta}\,|x-\bar{y}|^{2m-n-|\beta|}.
\end{equation}

\noindent Let us also point out that by formally differentiating
(\ref{eq1.12}) with respect to $y$ we obtain

\begin{equation}\label{eq1.23}
\left\{
\begin{array}{l}
{\mathcal L}(D_x)\,[D^\beta_yR_{\mathcal L}(x,y)]=0
\qquad\qquad\qquad\qquad\qquad 
\mbox{for}\,\,x\in\mathbb{R}^n,
\\[10pt]
\displaystyle{\left(\frac{\partial^j}{\partial x_n^j}
D^\beta_y R\right)((x',0),y)
=\left(\frac{\partial ^j}{\partial x_n^j}(-D)^\beta F\right)((x',0)-y)},
\,\,x'\in\mathbb{R}^{n-1},\,\,0\leq j\leq m-1.
\end{array}
\right.
\end{equation}

\noindent With (\ref{eq1.22}) and (\ref{eq1.23}) 
in place of (\ref{eq1.16}) and (\ref{eq1.12}), respectively, we can now 
run the same program as above and obtain the estimate 

\begin{equation}\label{Eq3}
\|D^\alpha_xD^\beta_y R(x,y)\|_{\mathbb{C}^{l\times l}}
\leq C_{\alpha\beta}\,|x-\bar{y}|^{2m-n-|\alpha|-|\beta|},
\end{equation}
 
\noindent for all multi-indices $\alpha$ and $\beta$.

\subsection{Proof of Theorem~\ref{th1} for $n\leq 2m$}

When $n\leq 2m$ we shall use the method of descent. To get started, fix an 
integer $N$ such that $N>2m$ and let $(x,z)\mapsto {\mathcal G}(x,y,z-\zeta)$ 
denote the Green matrix with singularity at $(y,\zeta)\in\RR^n\times\RR^{N-n}$
of the Dirichlet problem for the operator ${\mathcal L}(D_x)+(-\Delta_z)^m$ 
in the $N$-dimensional half-space 

\begin{equation}\label{RN}
\mathbb{R}^N_+:=
\{(x,z):\,z\in\mathbb{R}^{N-n},\,x=(x',x_n),\,x'\in\mathbb{R}^{n-1},\,x_n>0\}.
\end{equation}

\noindent Also, recall that $G(x,y)$ stands for the Green matrix of 
the problem (\ref{eq1.5}). Our immediate goal is to establish the following. 

\begin{lemma}\label{lem3}
For all multi-indices $\alpha$ and $\beta$ of order $m$ and for all 
$x$ and $y$ in $\mathbb{R}^n_+$

\begin{equation}\label{Eq12}
D^\alpha_x D^\beta_y G(x',y)
=\int_{\mathbb{R}^{N-n}}D^\alpha_x D^\beta_y{\mathcal G}(x,y,-\zeta)\,d\zeta.
\end{equation}
\end{lemma}

\noindent{\bf Proof.} The strategy is to show that 

\begin{equation}\label{pairGf}
\int_{\mathbb{R}^n_+}D^\alpha_xD_y^\beta G(x,y)\,f_\beta(y)\,dy
=\int_{\mathbb{R}^n_+}\int_{\mathbb{R}^{N-n}}D^\alpha_xD_y^\beta 
{\mathcal G}(x, y,-\zeta)\,d\zeta\,f_\beta(y)\,dy
\end{equation}

\noindent for each $f_\beta\in C^\infty_0(\mathbb{R}^n_+)$, from which 
(\ref{Eq12}) clearly follows. To justify (\ref{pairGf}) for a fixed,
arbitrary $f_\beta\in C^\infty_0(\mathbb{R}^n_+)$, we let $u$ be the unique 
vector-valued function satisfying $D^\alpha u\in L^2(\mathbb{R}^n_+)$ for 
all $\alpha$ with $|\alpha|=m$, and such that 

\begin{equation}\label{Eq9}
\left\{
\begin{array}{l}
{\mathcal L}(D_x)u=D^\beta_x f_\beta \qquad{\rm in}\,\,\mathbb{R}^n_+,
\\[6pt]
\displaystyle{\left(\frac{\partial^j u}{\partial x_n^j}\right)(x',0)=0
\qquad {\rm on}\,\,\mathbb{R}^{n-1}},\,\,0\leq j\leq m-1. 
\end{array}
\right.
\end{equation}

\noindent It is well-known that for each $\gamma\in\NN_0^n$ 

\begin{equation}\label{Eq10}
|D^\gamma u(x)|\leq C_\gamma\,|x|^{m-n-|\gamma|}\qquad{\rm for}\,\,|x|>1.
\end{equation}

\noindent This follows, for instance, from Theorem~6.1.4 {\bf\cite{KMR1}} 
combined with Theorem~10.3.2 {\bf\cite{KMR2}}. Also, as a consequence of 
Green's formula, the solution of the problem (\ref{Eq9}) satisfies

\begin{equation}\label{Eq14}
D^\alpha_x u(x)
=\int_{\mathbb{R}^n_+}D^\alpha_x(-D_y)^\beta G(x,y)\,f_\beta(y)\,dy.
\end{equation}

We shall now derive yet another integral representation formula for 
$D^\alpha_x u$ in terms of (derivatives of) ${\mathcal G}$ which is 
similar in spirit to (\ref{Eq14}). To get started, we note that 
since $N>2m$ the estimate (\ref{Eq3}) implies
 
\begin{equation}\label{Eq13}
\|D^\alpha_x D^\beta_y {\mathcal G}(x,y,-\zeta)\|_{\mathbb{C}^{l\times l}} 
\leq c\,(|x-y|+|\zeta|)^{-N}.
\end{equation}

\noindent Let us now fix $x\in\mathbb{R}^n_+$, $\rho>0$ and introduce 
a cut-off function $H\in C^\infty(\mathbb{R}^{N-n})$ which satisfies
$H(z)=1$ for $|z|\leq 1$ and $H(z)=0$ for $|z|\geq 2$. We may then write

\begin{equation}\label{u=G}
u(x)=\int_{\mathbb{R}^N} {\mathcal G}(x,y,-\zeta) 
\Bigl[H\bigl(\zeta/\rho\bigr)D^\beta f_\beta(y)+(-\Delta_\zeta)^m 
\bigl(H\bigl(\zeta/\rho\bigr)\,u(y)\bigr)\Bigr]\,dy\,d\zeta,
\end{equation}

\noindent which further implies

\begin{eqnarray}\label{est1}
&& \Bigl|D^\alpha_x u(x)-\int_{\mathbb{R}^N} D^\alpha_x(-D_y)^\beta\,
{\mathcal G}(x,y,-\zeta)\,H\bigl(\zeta/\rho\bigr)\,f_\beta(y)\,dy\,d\zeta\Bigr|
\nonumber\\[6pt]
&& \qquad\leq c\,\sum_{|\gamma|=m}\int_{\mathbb{R}^N_+}
\|D^\alpha_x D_\zeta^\gamma\,{\mathcal G}(x,y,-\zeta)\|
_{\mathbb{C}^{l\times l}}\,
\bigl|u(y)\,D^\gamma_\zeta\bigl(H\bigl(\zeta/\rho\bigr)\bigr)\bigr|\,d\zeta.
\end{eqnarray}

\noindent By (\ref{Eq10}) and (\ref{Eq13}), the expression in the 
right-hand side of (\ref{est1}) does not exceed

\begin{eqnarray*}
&& c\,\rho^{-m}\int_{\rho<|\zeta|<2\rho}d\zeta
\int_{\mathbb{R}^{n-1}}(|x-y|+|\zeta|)^{-N}\,|y|^{m-n}\,dy
\\[6pt]
&& \qquad\qquad
\leq c\,\rho^{N-n-m}\int_{\mathbb{R}^{n-1}}(|y|+\rho)^{-N}|y|^{m-n}\,dy
=c\,\rho^{-n}.
\end{eqnarray*}

\noindent This estimate, in concert with (\ref{Eq13}), allows us 
to obtain, after making $\rho\to\infty$, that 

\begin{equation}\label{DalphaU}
D^\alpha_x u(x)
=\int_{\mathbb{R}^n_+}\int_{\mathbb{R}^{N-n}}D^\alpha_x(-D_y)^\beta 
{\mathcal G}(x, y,-\zeta)\,d\zeta\,f_\beta(y)\,dy.
\end{equation}

\noindent Now (\ref{pairGf}) follows readily from this and (\ref{Eq14}). 
\hfill$\Box$
\vskip 0.08in

Having disposed of Lemma~\ref{lem3}, we are ready to present the

\vskip 0.08in
\noindent{\bf End of Proof of Theorem~\ref{th1}.} Assume that $2m\geq n$ 
and let $N$ be  again an integer such that $N>2m$. Denote by 
${\mathcal F}(x,z)$ the fundamental solution of the operator 
${\mathcal L}(D_x)+(-\Delta_z)^m$, which is positive homogeneous of 
degree $2m-N$ and is singular at $(0,0)\in\RR^n\times\RR^{N-n}$. 
Then the identity 

\begin{equation}\label{Eq17}
D^{\alpha+\beta}_x F(x)=\int_{\mathbb{R}^{N-n}}D^{\alpha+\beta}_x
{\mathcal F}(x,-\zeta)\,d\zeta
\end{equation}

\noindent can be established by proceeding as in the proof of Lemma~\ref{lem3}.
Combining (\ref{Eq17}) with Lemma~\ref{lem3}, we arrive at

\begin{equation}\label{Eq18}
D^\alpha_x D^\beta_y R(x',y)=\int_{\mathbb{R}^{N-n}}D^\alpha_x D^\beta_y 
{\mathcal R}(x,y,-\zeta)\,d\zeta,
\end{equation}

\noindent where ${\mathcal R}(x,y,z):={\mathcal G}(x,y,z)-{\mathcal F}(x-y,z)$.
Consequently, 

\begin{equation}\label{DalDbetU}
\|D^\alpha_x D^\beta_y {\mathcal R}(x,y,-\zeta)\|_{\mathbb{C}^{l\times l}} 
\leq C(|x-\bar{y}|+|\zeta|)^{-N}.
\end{equation}

\noindent by (\ref{eq1.20})  with $k=0$ and $N$ in place of $n$. 
This estimate, together with (\ref{Eq18}), yields (\ref{mainest}) and  
the proof of Theorem~\ref{th1} is therefore complete.
\hfill$\Box$
\vskip 0.08in

\section{Properties of integral operators in a half-space}
\setcounter{equation}{0}

\noindent In \S{3.1} and \S{3.2} we prove estimates for commutators 
(and certain commutator-like operators) between 
integral operators in $\mathbb{R}^n_+$ and multiplication operators with 
functions of bounded mean oscillations, in weighted Lebesgue spaces on 
$\mathbb{R}^n_+$. Subsection~3.3 contains ${\rm BMO}$ and pointwise 
estimates for extension operators from $\mathbb{R}^{n-1}$ onto 
$\mathbb{R}^n_+$.

\smallskip 

Throughout this section, given two Banach spaces $E,F$, we let 
${\mathfrak L}(E,F)$ stand for the space of bounded linear operators 
from $E$ into $F$, and abbreviate ${\mathfrak L}(E):={\mathfrak L}(E,E)$.
Also, given $p\in[1,\infty]$, an open set ${\mathcal O}\subset\RR^n$ 
and a measurable nonnegative function $w$ on ${\mathcal O}$, we let 
$L_p({\mathcal O},w(x)\,dx)$ denote the usual Lebesgue space of (classes of) 
functions which are $p$-th power integrable with respect to the weighted 
measure $w(x)\,dx$ on ${\mathcal O}$. Finally, following a well-established 
tradition, $A(r)\sim B(r)$ will mean that each quantity is dominated by 
a fixed multiple of the other, uniformly in the parameter $r$.

\subsection{Kernels with singularity at $\partial\mathbb{R}^{n}_+$}

Recall $L_p(\RR^n_+,\,x_n^{ap}\,dx)$ stands for the weighted Lebesgue space 
of $p$-th power integrable functions in $\RR^n_+$ corresponding to the 
weight $w(x):=x_n^{ap}$, $x=(x',x_n)\in\RR^n_+$. 

\begin{proposition}\label{tp3}
Let $a\in\RR$, $1<p<\infty$, and assume that ${\mathcal Q}$ is a 
non-negative measurable function on 
$\{\zeta=(\zeta',\zeta_n)\in\RR^{n-1}\times\RR:\,\zeta_n>-1\}$, 
which also satisfies

\begin{equation}\label{CC-1}
\int_{\mathbb{R}^n_+}
{\mathcal Q}(\zeta',\zeta_n-1)\,\zeta_n^{-a-1/p}\,d\zeta<\infty.
\end{equation}

\noindent Then the operator

\begin{equation}\label{CC-2}
Qf(x):=x_n^{-n}\int_{\mathbb{R}^n_+}{\mathcal Q}
\Bigl(\frac{y-x}{x_n}\Bigr)f(y)\,dy,\qquad x=(x',x_n)\in\RR^n_+,
\end{equation}

\noindent initially defined on functions $f\in L_p(\mathbb{R}^n_+)$ with 
compact support in $\mathbb{R}^n_+$, can be extended by continuity to an 
operator acting from $L_p(\RR^n_+,\,x_n^{a p}\,dx)$ into itself, with the 
norm satisfying

\begin{equation}\label{CC-3}
\|Q\|_{{\mathfrak L}(L_p(\RR^n_+,\,x_n^{ap}dx))}\leq\int_{\RR^n_+}
{\mathcal Q}(\zeta',\zeta_n-1)\,\zeta_n^{-a-1/p}\,d\zeta.
\end{equation}
\end{proposition}

\noindent{\bf Proof.} Introducing the new variable 
$\zeta:=(x_n^{-1}(y'-x'), x_n^{-1}y_n)\in\RR^n_+$, we may write 

\begin{equation}\label{CC-4}
|Qf(x)|\leq\int_{\RR^n_+}{\mathcal Q}(\zeta',\zeta_n-1) 
|f(x' +x_n\zeta', x_n\zeta_n)|d\zeta,\qquad\forall\,x\in\RR^n_+.
\end{equation}

\noindent Then, by Minkowski's inequality,

\begin{eqnarray}\label{CC-5}
\|Qf\|_{L_p(\RR^n_+,x_n^{a p}\,dx)} & \leq & \int_{\RR^n_+}
{\mathcal Q}(\zeta',\zeta_n-1)\Bigl(\int_{\RR^n_+} 
x_n^{a p}|f(x'+x_n\zeta',x_n\zeta_n)|^p\,dx\Bigr)^{1/p}d\zeta
\nonumber\\[6pt]
& = & \Bigl(\int_{\RR^n_+}{\mathcal Q}(\zeta',\zeta_n-1)\,
\zeta_n^{-a-1/p}\,d\zeta\Bigr)\|f\|_{L_p(\RR^n_+,x_n^{a p}\,dx)},
\end{eqnarray}

\noindent as desired. 
\hfill$\Box$
\vskip 0.08in

Recall that $\bar{y}:=(y',-y_n)$ if $y=(y',y_n)\in\RR^{n-1}\times\RR$. 

\begin{corollary}\label{Cor1}
Consider 

\begin{equation}\label{1a}
Rf(x):=\int_{\mathbb{R}^n_+}\frac{\log\,\bigl(\frac{|x-y|}{x_n}+2\bigr)}
{|x-\bar{y}|^n} f(y)\,dy,\qquad x=(x',x_n)\in\RR^n_+.
\end{equation}
 
\noindent Then for each $1<p<\infty$ and each $a\in (-1/p,1-1/p)$ the 
operator $R$ is bounded from $L_p(\RR^n_+,\,x_n^{a p}\,dx)$ into itself.
Moreover, 

\begin{equation}\label{70a}
\|R\|_{{\mathfrak L}(L_p(\RR^n_+,\,x_n^{a p}\,dx))}
\leq\frac{c(n)\,p^2}{(pa+1)(p(1-a)-1)}=\frac{c(n)}{s(1-s)},
\end{equation}

\noindent where $s=1-a-1/p$ and $c(n)$ is independent of $p$ and $a$.
\end{corollary}

\noindent{\bf Proof.} The result follows from Proposition~\ref{tp3} with

\begin{equation}\label{CC-6}
{\mathcal Q}(\zeta):=\frac{\log\,(|\zeta|+2)}{(|\zeta|^2+1)^{n/2}},
\end{equation}
 
\noindent and from the obvious inequality $2|x-\bar{y}|^2 \geq |x-y|^2 +x_n^2$.
\hfill$\Box$
\vskip 0.08in

Let us note here that Corollary~\ref{Cor1} immediately yields the following. 

\begin{corollary}\label{Cor2}
Consider 

\begin{equation}\label{2a}
Kf(x):=\int_{\RR^n_+}\frac{f(y)}{|x-\bar{y}|^n}\,dy,\qquad x\in\RR^n_+.
\end{equation}
 
\noindent Then for each $1<p<\infty$ and $a\in (-1/p,1-1/p)$ the 
operator $K$ is bounded from $L_p(\RR^n_+,\,x_n^{a p}\,dx)$ into itself. 
Moreover, 

\begin{equation}\label{71a}
\|K\|_{{\mathfrak L}(L_p(\RR^n_+,\,x_n^{a p}\,dx))}
\leq \frac{c(n)\,p^2}{(pa+1)(p(1-a)-1)}=\frac{c(n)}{s(1-s)},
\end{equation}

\noindent where $s=1-a-1/p$ and $c(n)$ is independent of $p$ and $a$.
\end{corollary}

Recall that the barred integral stands for the mean-value (taken in the 
integral sense). 

\begin{lemma}\label{Lem1}
Assume that $1<p<\infty$, $a\in (-1/p,1-1/p)$, and recall the operator 
$K$ introduced in {\rm (\ref{2a})}. Further, consider a non-negative, 
measurable function $w$ defined on $\RR^n_+$ and fix a family of 
balls ${\mathcal F}$ which form a Whitney covering of $\RR^n_+$. 
Then the norm of $wK$ as an operator from $L_p(\RR^n_+,\,x_n^{a p}\,dx)$ 
into itself is equivalent to

\begin{equation}\label{CC-7}
\sup\limits_{B\in{\mathcal F}}\meanint_{\!\!\!\!B}w(y)^p\,dy.
\end{equation}

\noindent Furthermore,

\begin{equation}\label{CC-70}
\|w\,K\|_{{\mathfrak L}(L_p(\RR^n_+,\,x_n^{a p}\,dx))}\leq\frac{c(n)}{s(1-s)}
\sup\limits_{B\in{\mathcal F}}\Bigl(\meanint_{\!\!\!\!B}w(y)^p\,dy\Bigr)^{1/p},
\end{equation}

\noindent where $c(n)$ is independent of $w$, $p$, and $\alpha$.
\end{lemma}

\noindent{\bf Proof.} Fix $f\geq 0$ and denote by $|B|$ the Euclidean 
volume of $B$. Sobolev's embedding theorem allows us to write

\begin{equation}\label{CC-8}
\|Kf\|^p_{L_\infty(B)}\leq c(n)\,|B|^{-1}\sum_{j=0}^n |B|^{jp/n}
\|\nabla_j Kf\|^p_{L_p(B)},\qquad\forall\,B\in{\mathcal F}. 
\end{equation}

\noindent Hence,

\begin{equation}\label{N1-bis}
\int_{\RR^n_+}|x_n^{a}w(x)(Kf)(x)|^p\,dx\leq c(n)\,
\sup\limits_{B\in{\mathcal F}}\meanint_{\!\!\!\!B}w(y)^p\,dy
\int_{\RR^n_+}x_n^{pa}\sum_{0\leq j\leq l}x_n^{jp}|\nabla_j Kf|^p\,dx. 
\end{equation}

\noindent Observing that $x_n^j|\nabla_j\, Kf|\leq c(n)\,Kf$ and referring to 
Corollary~\ref{Cor2}, we arrive at the required upper estimate for the norm 
of $wK$. The lower estimate is obvious.
\hfill$\Box$
\vskip 0.08in

We momentarily pause in order to collect some definitions and set up 
basic notation pertaining to functions with bounded mean oscillations.
Let $f$ be a locally integrable function defined on $\mathbb{R}^n$ and 
define the seminorm

\begin{equation}\label{semi1}
[f]_{{\rm BMO}(\mathbb{R}^n)}:=\sup_{B}
\meanint_{\!\!\!B}\,\Bigl|f(x)-\meanint_{\!\!\!B}f(y)\,dy\Bigr|\,dx, 
\end{equation}

\noindent where the supremum is taken over all balls $B$ in ${\mathbb{R}^n}$. 
Next, if $f$ is a locally integrable function defined on $\mathbb{R}^n_+$, 
we set 

\begin{equation}\label{semi2}
[f]_{{\rm BMO}(\mathbb{R}^n_+)}:=\mathop{\hbox{sup}}_{\{B\}}
\meanint_{\!\!\!B\cap\mathbb{R}^n_+}\,
\Bigl|f(x)-\meanint_{\!\!\!B\cap\mathbb{R}^n_+}f(y)\,dy\Bigr|\,dx, 
\end{equation}
 
\noindent where, this time, the supremum is taken over the collection 
$\{B\}$ of all balls $B$ with centers in $\overline{\mathbb{R}^n_+}$.
Then the following inequalities are straightforward

\begin{equation}\label{N4-bis}
[f]_{{\rm BMO}(\mathbb{R}^n_+)}\leq\mathop{\hbox{sup}}_{\{B\}}
\meanint_{\!\!\!B\cap\mathbb{R}^n_+}\,\meanint_{\!\!\!B\cap\mathbb{R}^n_+}
\Bigl|f(x)-f(y)\,\Bigr|\,dxdy\leq 2\,[f]_{{\rm BMO}(\mathbb{R}^n_+)}.
\end{equation}

\noindent We also record here the equivalence relation 

\begin{equation}\label{semi}
[f]_{{\rm BMO}(\mathbb{R}^n_+)}\sim [{\rm Ext}\,f]_{{\rm BMO}(\mathbb{R}^n)},
\end{equation}

\noindent where ${\rm Ext}\,f$ is the extension of $f$ onto 
$\mathbb{R}^n$ as an even function in $x_n$.

Finally, by ${\rm BMO}({\mathbb{R}^n_+})$ we denote the collection of 
equivalence classes, modulo constants, of functions $f$ on $\mathbb{R}^n_+$ 
for which $[f]_{{\rm BMO}({\mathbb{R}^n_+})}<\infty$.

\begin{proposition}\label{tp3-bis}
Let $b\in{\rm BMO}(\RR^n_+)$ and consider the operator 

\begin{equation}\label{eqp8-bis}
Tf(x):=\int_{\RR^n_+}\frac{|b(x)-b(y)|}{|x-\bar{y}|^n}f(y)\,dy,
\qquad x\in\RR^n_+.
\end{equation}

\noindent Then for each $p\in(1,\infty)$ and $a\in (-1/p,1-1/p)$

\begin{equation}\label{eqp9bis}
T:L_p(\RR^n_+,\,x_n^{a p}\,dx)\longrightarrow
L_p(\RR^n_+,\,x_n^{a p}\,dx)
\end{equation}

\noindent is a well-defined, bounded operator with 
\begin{equation}\label{CC-71}
\|T\|_{{\mathfrak L}(L_p(\RR^n_+,\,x_n^{a p}\,dx))}
\leq\frac{c(n)}{s(1-s)}\,[b]_{{\rm BMO}(\RR^n_+)},
\end{equation}

\noindent where $c(n)$ is  a constant which depends only on $n$.
\end{proposition}

\noindent{\bf Proof.} Given $x\in\RR^n_+$ and $r>0$, we shall use 
the abbreviations

\begin{equation}\label{CC-9}
\bar{b}_r(x):=\meanint_{\!\!\!B(x,r)\cap\RR^n_+}b(y)\,dy,
\qquad\quad D_r(x):=|b(x)-\bar{b}_{r}(x)|,
\end{equation}

\noindent and make use of the integral operator

\begin{equation}\label{CC-10}
Sf(x):=\int_{\mathbb{R}^n_+}\frac{D_{|x-\bar{y}|}(x)}{|x-\bar{y}|^n}
\,f(y)\,dy,\qquad x\in\RR^n_+,
\end{equation}

\noindent as well as its adjoint $S^*$. Clearly, for each nonnegative, 
measurable function $f$ on $\RR^n_+$ and each $x\in\RR^n_+$, 

\begin{eqnarray}\label{CC-11}
Tf(x) & \leq & Sf(x)+S^*f(x)+\int_{\RR^n_+}
\frac{|\bar{b}_{|x-\bar{y}|}(x)-\bar{b}_{|x-\bar{y}|}(y)|}{|x-\bar{y}|^n}\, 
f(y)dy
\nonumber\\[6pt]
& \leq & Sf(x)+S^*f(x)+c(n)\,[b]_{{\rm BMO}(\RR^n_+)}Kf(x),
\end{eqnarray}

\noindent where $K$ has been introduced in (\ref{2a}). Making use of 
Corollary~\ref{Cor2}, we need to estimate only the norm of $S$. Obviously,

\begin{equation}\label{CC-12}
Sf(x)\leq D_{x_n}(x)Kf(x)+\int_{\RR^n_+}
\frac{|\bar{b}_{x_n}(x)-\bar{b}_{|x-\bar{y}|}(x)|}{|x-\bar{y}|^n}\,f(y)\,dy.
\end{equation}

\noindent Setting $r=|x-\bar{y}|$ and $\rho=x_n$ in the standard inequality

\begin{equation}\label{CC-13}
|\bar{b}_\rho(x)-\bar{b}_r(x)|\leq c(n)\,\log\,\Bigl(\frac{r}{\rho}+1\Bigr)
[b]_{{\rm BMO}(\RR^n_+)},
\end{equation}

\noindent where $r>\rho$ (cf., e.g., p.\,176 in {\bf\cite{MS}}, or p.\,206 in 
{\bf\cite{Tor}}), we arrive at

\begin{equation}\label{CC-14}
Sf(x)\leq D_{x_n}(x)Kf(x)+ c(n)\,[b]_{BMO(\mathbb{R}^n_+)}\, Rf(x),
\end{equation}

\noindent where $R$ is defined in (\ref{1a}). Let ${\mathcal F}$ 
be a Whitney covering of $\RR^n_+$ with open balls. For an arbitrary
$B\in{\mathcal F}$, denote by $\delta$ the radius of $B$. 
By Lemma~\ref{Lem1} with $w(x):=D_{x_n}(x)$, the norm of the operator 
$D_{x_n}(x)K$ does not exceed

\begin{eqnarray}\label{CC-15}
\sup\limits_{B\in{\mathcal F}}\Bigl(\meanint_{\!\!\!B}|D_{x_n}(x)|^p\,dx
\Bigr)^{1/p} 
& \leq & c(n)\,\sup\limits_{B\in{\mathcal F}}\Bigl(\meanint_{\!\!\!B}
|b(x)-\bar{b}_\delta(x)|^p\,dx\Bigr)^{1/p}
+c(n)\,[b]_{{\rm BMO}(\mathbb{R}^n_+)}
\nonumber\\[6pt]
&\leq & c(n)\,[b]_{{\rm BMO}(\mathbb{R}^n_+)},
\end{eqnarray}

\noindent by the John-Nirenberg inequality. Here we have also used the 
triangle inequality and the estimate (\ref{CC-13}) in order to replace 
$\bar{b}_{x_n}(x)$ in the definition of $D_{x_n}(x)$ by $\bar{b}_{\delta}(x)$.
The intervening logarithmic factor is bounded independently of $x$ since $x_n$
is comparable with $\delta$, uniformly for $x\in B$. With this estimate
in hand, a reference to Corollary~\ref{Cor1} gives that 

\begin{eqnarray}\label{CC-16}
&& S:L_p(\RR^n_+,\,x_n^{a p}\,dx)\to L_p(\RR^n_+,\,x_n^{a p}\,dx)
\mbox{ boundedly} 
\\[6pt]
&&\mbox{for each }p\in(1,\infty)\mbox{ and each }a\in (-1/p,1-1/p).
\nonumber
\end{eqnarray}

\noindent The corresponding estimate for the norm $S$ results as well. 
By duality, it follows that $S^*$ enjoys the same property and,
hence, the operator $T$ is bounded on $L_p(\RR^n_+,\,x_n^{a p}\,dx)$
for each $p\in(1,\infty)$ and $a\in(-1/p,1-1/p)$, thanks to (\ref{CC-11})
and Corollary~\ref{Cor2}. The fact that the operator norm of $T$ admits the 
desired estimate is implicit in the above reasoning and this finishes 
the proof of the proposition.
\hfill$\Box$
\vskip 0.08in

\subsection{Singular integral operators}

We need the analogue of Proposition~\ref{tp3-bis} for the class 
of Mikhlin-Calder\'on-Zygmund singular integral operators. Recall that

\begin{equation}\label{CZ-op}
{\mathcal S}f(x)=p.v.\int_{\RR^n}k(x,x-y)f(y)\,dy,\qquad x\in\RR^n,
\end{equation}

\noindent (where $p.v.$ indicates that the integral is taken in the
principal value sense, which means excluding balls centered at the
singularity and then passing to the limit as the radii shrink to zero),
is called a Mikhlin-Calder\'on-Zygmund operator (with a variable 
coefficient kernel) provided the function 
$k:\RR^n\times(\RR^n\setminus\{0\})\to\RR$ satisfies:

\begin{itemize}
\item[(i)] $k(x,\cdot)\in C^\infty(\RR^n\setminus\{0\})$
and, for almost each $x\in\RR^n$,

\begin{equation}\label{kk-est}
\max_{|\alpha|\leq 2n}\|D_z^\alpha k(x,z)\|_{L_\infty(\RR^n\times S^{n-1})}
<\infty,
\end{equation}

\noindent where $S^{n-1}$ is the unit sphere in $\RR^n$;

\item[(ii)]  $k(x,\lambda z)=\lambda ^{-n}k(x,z)$ for each $z\in\RR^n$
and each $\lambda\in\RR$, $\lambda>0$;
\item[(iii)] $\int_{S^{n-1}} k(x,\omega)\,d\omega=0$,
where $d\omega$ indicates integration with respect to $\omega\in S^{n-1}$.
\end{itemize}
\vskip 0.08in

It is well-known that the Mikhlin-Calder\'on-Zygmund operator ${\mathcal S}$ 
and its commutator $[{\mathcal S},b]$ with the operator of multiplication by 
a function $b\in{\rm BMO}(\mathbb{R}^n_+)$ are bounded operators in 
$L_p(\mathbb{R}^n_+)$ for each $1<p<\infty$. The norms of these operators 
admit the estimates

\begin{equation}\label{b1}
\|{\mathcal S}\|_{{\mathfrak L}(L_p(\mathbb{R}^n_+))}\leq c(n)\,p\,p',\qquad 
\|\,[{\mathcal S},b]\,\|_{{\mathfrak L}(L_p(\mathbb{R}^n_+))}
\leq c(n)\,p\,p'\,[b]_{{\rm BMO}(\RR^n_+)},
\end{equation}

\noindent where $c(n)$ depends only on $n$  and the quantity 
in (\ref{kk-est}). The first estimate in (\ref{b1}) goes back to the work 
of A.\,Calder\'on and A.\,Zygmund (cf., e.g., {\bf\cite{CaZy}}, 
{\bf\cite{CaZy2}}; see also the comment on p.\,22 of {\bf\cite{St}} 
regarding the dependence on the parameter $p$ of the constants involved).
The second estimate in (\ref{b1}) was originally proved for 
convolution type operators by R.\,Coifman, R.\,Rochberg and G.\,Weiss in 
{\bf\cite{CRW}} and a standard expansion in spherical harmonics allows
to extend this result to the case of operators with variable-kernels 
of the type considered above. 

We are interested in extending (\ref{b1}) to the weighted case, i.e. 
when the measure $dx$ on $\RR^n_+$ is replaced by $x_n^{ap}\,dx$,
where $1<p<\infty$ and $a\in(-1/p,1-1/p)$. Parenthetically, we wish to 
point out that $a\in(-1/p,1-1/p)$ corresponds precisely to the range 
of $a$'s for which $w(x):=x_n^{ap}$ is a weight in Muckenhoupt's $A_p$ 
class, and that while in principle this observation can help with the 
goal just stated, we prefer to give a direct, self-contained proof. 

\begin{proposition}\label{Prop4}
Retain the above conventions and hypotheses. Then the operator ${\mathcal S}$ 
and its commutator $[{\mathcal S},b]$ with a function 
$b\in{\rm BMO}(\mathbb{R}^n_+)$ are bounded when acting from 
$L_p(\RR^n_+,\,x_n^{ap}\,dx)$ into itself for each $p\in(1,\infty)$ 
and $a\in (-1/p, 1-1/p)$. The norms of these operators satisfy 

\begin{eqnarray}\label{b2}
&&\|{\mathcal S}\|_{{\mathfrak L}(L_p(\mathbb{R}^n_+,\,x_n^{ap}\,dx))}
\leq c(n)\Bigl(p\,p'+\frac{1}{s(1-s)}\Bigr),
\\[6pt]
&&\|\,[{\mathcal S},b]\,\|_{{\mathfrak L}(L_p(\mathbb{R}^n_+,\,x_n^{ap}\,dx))} 
\leq c(n)\Bigl(p\,p'+\frac{1}{s(1-s)}\Bigr)\,[b]_{{\rm BMO}(\RR^n_+)}.
\label{b2'}
\end{eqnarray}
\end{proposition}

\noindent{\bf Proof.} Let $\chi_j$ be the characteristic function of the 
layer $2^{j/2}<x_n\leq 2^{1+j/2}$, $j=0,\pm 1,\ldots$, so that
$\sum_{j\in\ZZ}\chi_j=2$. We then write ${\mathcal S}$ as the sum 
${\mathcal S}_1+{\mathcal S}_2$, where

\begin{equation}\label{CC-17}
{\mathcal S}_1:=\frac{1}{4}\sum_{|j-k|\leq 3}\chi _j{\mathcal S}\chi_k.
\end{equation}

\noindent The following chain of inequalities is evident

\begin{eqnarray}\label{CC-18}
\|{\mathcal S}_1\, f\|_{L_p(\mathbb{R}^n_+,\,x_n^{ap}\,dx)}
& \leq & \Bigl(\sum_j\int_{\RR^n_+}\chi_j(x)\,
\Bigl|{\mathcal S}\Bigl(\sum_{|k-j|\leq 3}\chi_k f\Bigr)(x)\Bigr|^p\,
x_n^{ap}\,dx\Bigr)^{1/p}
\nonumber\\[6pt]
& \leq & c(n)\Bigl(\sum_j\int_{\mathbb{R}^n_+}
\Bigl|{\mathcal S}\Bigl(\sum_{|k-j|\leq 3} 
\chi_k 2^{ja/2} f\Bigr)(x)\Bigr|^p\,dx\Bigr)^{1/p}.
\end{eqnarray}

\noindent In concert with the first estimate in (\ref{b1}), this entails

\begin{eqnarray}\label{CC-19}
\|{\mathcal S}_1\,f\|_{L_p(\mathbb{R}^n_+,\,x_n^{ap}\,dx)}
& \leq & c(n)\,p\, p'\Bigl(\sum_j 
\int_{\RR^n_+}\Bigl(\sum_{|k-j|\leq 3}
\chi_k 2^{ja/2}|f|\Bigr)^p\,dx\Bigr)^{1/p}
\nonumber\\[6pt]
& \leq & c(n)\,p\,p'\Bigl(\int_{\RR^n_+}|f(x)|^p\,x_n^{ap}\,dx\Bigr)^{1/p},
\end{eqnarray}

\noindent which is further equivalent to

\begin{equation}\label{b3}
\|{\mathcal S}_1\|_{{\mathfrak L}(L_p(\mathbb{R}^n_+,\,x_n^{ap}\,dx))}
\leq c(n)\, p\, p'.
\end{equation}
  
\noindent Applying the same  argument to $[{\mathcal S},b]$ and referring 
to (\ref{b1}), we arrive at
  
\begin{equation}\label{b4}
\|\,[{\mathcal S}_1,b]\,\|_{{\mathfrak L}(L_p(\mathbb{R}^n_+,\,x_n^{ap}\,dx))} 
\leq c(n)\,p\,p'\,[b]_{{\rm BMO}(\RR^n_+)}.
\end{equation}

It remains to obtain the analogues of (\ref{b3}) and ({\ref{b4}) with 
${\mathcal S}_2$ in place of ${\mathcal S}_1$. One can check directly that the 
modulus of the kernel of ${\mathcal S}_2$ does not exceed 
$c(n)\,|x-\bar{y}|^{-n}$ and that the modulus of the kernel of 
$[{\mathcal S}_2,b]$ is majorized by $c(n)\,|b(x)-b(y)|\,|x-\bar{y}|^{-n}$. 
Then the desired conclusions follow from 
Corollary~\ref{Cor2} and Proposition~\ref{tp3-bis}. 
\hfill$\Box$
\vskip 0.08in

\subsection{The Gagliardo extension operator}  

Here we shall revisit a certain operator $T$, extending functions 
defined on $\mathbb{R}^{n-1}$ into functions defined on $\mathbb{R}^n_+$,
first introduced by Gagliardo in {\bf\cite{Ga}}. Fix a smooth, radial, 
decreasing, even, non-negative function $\zeta$ in $\mathbb{R}^{n-1}$ 
such that $\zeta(t)=0$ for $|t|\geq 1$ and 

\begin{equation}\label{zeta-int}
\int\limits_{\mathbb{R}^{n-1}}\zeta(t)\,dt=1.
\end{equation}

\noindent (A standard choice is 
$\zeta(t):=c\,{\rm exp}\,(-1/(1-|t|^2)_+)$ for a suitable $c$.) 
Following {\bf\cite{Ga}} we then define 

\begin{equation}\label{10.1.20}
(T\varphi)(x',x_n)=\int\limits_{\mathbb{R}^{n-1}}\zeta(t)\varphi(x'+x_nt)\,dt,
\qquad(x',x_n)\in\RR^n_+,
\end{equation}

\noindent acting on functions $\varphi$ from $L_1(\mathbb{R}^{n-1},loc)$. 
To get started, we note that 

\begin{eqnarray}\label{10.1.29}
\nabla_{x'}(T\varphi)(x',x_n)
& = & \int\limits_{\mathbb{R}^{n-1}}\zeta(t)\nabla\varphi(x'+tx_n)\,dt,
\\[6pt]
{\partial\over\partial x_n}(T\varphi)(x',x_n)
& = & \int\limits_{\mathbb{R}^{n-1}}\zeta(t)\,t\,\nabla\varphi(x'+t x_n)\,dt,
\label{10.1.30}
\end{eqnarray}

\noindent and, hence, we have the estimate

\begin{equation}\label{10.1.21}
\|\nabla_x\,(T\varphi)\|_{L_\infty(\mathbb{R}^n_+)} 
\leq c\,\|\nabla_{x'}\,\varphi\|_{L_\infty(\mathbb{R}^{n-1})}.
\end{equation}

\noindent Refinements of (\ref{10.1.21}) are contained in the 
Lemmas~\ref{lem7}-\ref{lem8} below. 

\begin{lemma}\label{lem7}
(i) For all $x\in\mathbb{R}^n_+$ and for all multi-indices $\alpha$ 
with $|\alpha|>1$,

\begin{equation}\label{1.2}
\Bigl|D^\alpha_{x}(T\varphi)(x)\Bigr|
\leq c\,x_n ^{1-|\alpha|}[\nabla\varphi]_{\rm BMO(\mathbb{R}^{n-1})}.
\end{equation}

(ii) For all $x=(x',x_n)\in\mathbb{R}^n_+$, 

\begin{equation}\label{Tfi}
\Bigl|(T\varphi)(x)-\varphi(x')\Bigr|
\leq c\,x_n[\nabla\varphi]_{\rm BMO(\mathbb{R}^{n-1})}.
\end{equation}
\end{lemma}

\noindent{\bf Proof.} Rewriting (\ref{10.1.30}) as 

\begin{equation}\label{1.3}
{\partial\over\partial x_n}(T\varphi)(x',x_n)
=x_n^{1-n}\int\limits_{\mathbb{R}^{n-1}}\zeta\Bigl(\frac{\xi-x'}{x_n}\Bigr)
\frac{\xi-x'}{x_n}\Bigl(\nabla\varphi(\xi)-\meanint_{\!\!\!|z-x'|<x_n} 
\nabla\varphi(z)dz\Bigr)d\xi
\end{equation}

\noindent we obtain

\begin{equation}\label{1.12}
\Bigl|D^\gamma_{x}{\partial\over\partial x_n}(T\varphi)(x)\Bigr|
\leq c\,x_n^{-|\gamma|}[\nabla\varphi]_{\rm BMO(\mathbb{R}^{n-1})}
\end{equation}

\noindent for every non-zero multi-index $\gamma$. Furthermore, for 
$i=1,\ldots n-1$, by (\ref{10.1.29})

\begin{equation}\label{TT1}
\frac{\partial}{\partial x_i}\nabla_{x'}(T\varphi)(x) 
=x_n^{1-n}\int\limits_{\mathbb{R}^{n-1}}\partial_i\zeta
\Bigl(\frac{\xi-x'}{x_n}\Bigr)\Bigl(\nabla\varphi(\xi)
-\meanint_{\!\!\!|z-x'|<x_n}\nabla\varphi(z)dz\Bigr)d\xi,
\end{equation}

\noindent where $\partial_i$ is the differentiation with respect to 
the $i$-th component of the argument. Hence once again

\begin{equation}\label{TT2}
\Bigl|D^\gamma_x\frac{\partial}{\partial x_i}\nabla_{x'}(T\varphi)(x)\Bigr| 
\leq c\,x_n^{-|\gamma|-1}[\nabla\varphi]_{\rm BMO(\mathbb{R}^{n-1})},
\end{equation}

\noindent and the estimate claimed in (i) follows.

Finally, (ii) is a simple consequence of (i) and the fact that 
$(T\varphi)|_{\RR^{n-1}}=\varphi$. 
\hfill$\Box$
\vskip 0.08in

\noindent{\bf Remark.} In concert with Theorem~2 on p.\,62-63 in 
{\bf\cite{St}}, formula (\ref{10.1.29}) yields the pointwise estimate

\begin{equation}\label{HL-Max}
|\nabla\,(T\varphi)(x)|
\leq c\,{\mathcal M}(\nabla\varphi)(x'),\qquad x=(x',x_n)\in\RR^{n}_+,
\end{equation}

\noindent where ${\mathcal M}$ is the classical Hardy-Littlewood 
maximal function (cf., e.g., Chapter~I in {\bf\cite{St}}). 
As for higher order derivatives, an inspection of the above proof 
reveals that 

\begin{equation}\label{sharp}
\Bigl|D_x^\alpha(T\varphi)(x)\Bigr|
\leq c\,x_n^{1-|\alpha|}(\nabla\varphi)^\#(x'),\qquad (x',x_n)\in\RR^{n},
\end{equation}

\noindent holds for each multi-index $\alpha$ with $|\alpha|>1$, 
where $(\cdot)^\#$ is the Fefferman-Stein sharp maximal function 
(cf. {\bf\cite{FS}}). 

\smallskip 

\begin{lemma}\label{lem8}
If $\nabla_{x'}\varphi\in{\rm BMO}(\mathbb{R}^{n-1})$ then 
$\nabla(T\varphi)\in{\rm BMO}(\mathbb{R}^n_+)$ and 

\begin{equation}\label{1.8}
[\nabla(T\varphi)]_{{\rm BMO}(\mathbb{R}^n_+)}
\leq c\,[\nabla_{x'}\varphi]_{{\rm BMO}(\mathbb{R}^{n-1})}.
\end{equation}
\end{lemma}

\noindent{\bf Proof.} Since $(T\varphi)(x',x_n)$ is even with respect to $x_n$,
it suffices to estimate $[\nabla_x(T\varphi)]_{{\rm BMO}(\mathbb{R}^n)}$. 
Let $Q_r$ denote a cube with side-length $r$ centered at the point 
$\eta =(\eta',\eta_n)\in\RR^{n-1}\times\RR$. Also let $Q'_r$ be the 
projection of $Q_r$ on $\mathbb{R}^{n-1}$. Clearly,

\begin{equation}\label{xxx}
\nabla_{x'}(T\varphi)(x',x_n)-\nabla _{x'}\varphi(x')
=x_n^{1-n}\int\limits_{\mathbb{R}^{n-1}}\zeta\Bigl(\frac{\xi -x'}{x_n}\Bigr)
(\nabla\varphi(\xi)-\nabla\varphi(x'))\,d\xi.
\end{equation}

Suppose that $|\eta_n|<2r$ and write 

\begin{eqnarray}\label{1.10}
\int_{Q_r}\Bigl|\nabla_{x'}(T\varphi)(x',x_n)-\nabla_{x'}\varphi(x')\Bigr|\,dx 
& \leq & c\,r^{2-n}\int_{Q'_{4r}}
\int_{Q'_{4r}}|\nabla\varphi(\xi)-\nabla\varphi(z)|\,dz\,d\xi.
\nonumber\\[6pt]
& \leq & c\,r^n[\nabla\varphi]_{{\rm BMO}(\mathbb{R}^{n-1})}.
\end{eqnarray}

Therefore, for $|\eta_n|<2r$

\begin{eqnarray}\label{1.18}
\meanint_{\!\!\!Q_r}\meanint_{\!\!\!Q_r}|\nabla_{x'}T\varphi(x)-\nabla_{y'}
T\varphi(y)|\,dxdy
& \leq & 2\meanint_{\!\!\!Q_r}|\nabla_{x'}T\varphi(x)-\nabla\varphi(x')|\,dx
\nonumber\\[6pt]
&&+\meanint_{\!\!\!Q'_r}\meanint_{\!\!\!Q'_r}
|\nabla\varphi(x')-\nabla\varphi(y')|\,dx'dy' 
\nonumber\\[6pt]
& \leq & c\,[\nabla\varphi]_{{\rm BMO}(\mathbb{R}^{n-1})}.
\end{eqnarray}

Next, consider the case when $|\eta_n|\geq 2r$ and let $x$ and $y$ 
be arbitrary points in $Q_r(\eta)$. Then, using the generic abbreviation 
$\bar{f}_E:={\displaystyle{\meanint_{\!\!\!E}f}}$, we may write 

\begin{eqnarray}\label{arTT}
|\nabla_{x'}T\varphi(x)-\nabla_{y'}T\varphi(y)|
& \leq & 
\int\limits_{\mathbb{R}^{n-1}}\Bigl|x_n^{1-n}\zeta\Bigl(\frac{\xi-x'}{x_n}
\Bigr)-y_n^{1-n}\zeta\Bigl(\frac{\xi-y'}{y_n}\Bigr)
\Bigr|\Bigl|\nabla\varphi(\xi)- {\overline{\nabla\varphi}}_{Q'_{2|\eta_n|}}
\Bigr|\,d\xi 
\nonumber\\[6pt]
& \leq & \frac{c\,r}{|\eta_n|^n}\int_{Q'_{2|\eta_n|}}
\Bigl|\nabla\varphi(\xi)
-{\overline{\nabla\varphi}}_{Q'_{2|\eta_n|}}\Bigr|\,d\xi 
\nonumber\\[6pt]
& \leq & c\,[\nabla\varphi]_{{\rm BMO}(\mathbb{R}^{n-1})}. 
\end{eqnarray}

\noindent Consequently, for $|\eta_n|\geq 2r$, 

\begin{equation}\label{QrT3}
\meanint_{\!\!\!Q_r}\meanint_{\!\!\!Q_r}
|\nabla_{x'}T\varphi(x)-\nabla_{y'}T\varphi(y)|\,dxdy
\leq c\,[\nabla\varphi]_{{\rm BMO}(\mathbb{R}^{n-1})}
\end{equation}

\noindent which, together with (\ref{1.18}), gives

\begin{equation}\label{QrT4}
[\nabla_{x'}T\varphi]_{{\rm BMO}(\mathbb{R}^{n})}
\leq c[\nabla\varphi]_{\rm BMO(\mathbb{R}^{n-1})}.
\end{equation}

\noindent This inequality and (\ref{1.12}), where $|\gamma|=0$, 
imply (\ref{1.8}). 
\hfill$\Box$
\vskip 0.08in

\section{The Dirichlet problem in $\mathbb{R}^n_+$ for variable coefficient
systems}
\setcounter{equation}{0}

\subsection{Preliminaries}

For 

\begin{equation}\label{indices}
1<p<\infty,\quad -\frac{1}{p}<a<1-\frac{1}{p}\quad 
\mbox{and}\quad m\in\NN, 
\end{equation}

\noindent we let $V^{m,a}_p(\mathbb{R}^n_+)$ denote the weighted Sobolev space 
associated with the norm

\begin{equation}\label{defVVV}
\|u\|_{V_p^{m,a}(\mathbb{R}^n_+)}:=\Bigl(\sum_{0\leq|\beta|\leq m} 
\int_{\mathbb{R}^n_+}|x_n^{|\beta|-m} D^\beta u(x)|^p\,x_n^{pa}\,dx
\Bigr)^{1/p}.
\end{equation}

\noindent It is easily proved that $C_0^\infty(\mathbb{R}^n_+)$ is dense in 
$V^{m,a}_p(\mathbb{R}^n_+)$. Moreover, by the one-dimensional Hardy's 
inequality (see, for instance, {\bf\cite{Maz1}}, formula (1.3/1)), we have 

\begin{equation}\label{4.60}
\|u\|_{V_p^{m,a}(\mathbb{R}^n_+)}\leq cs^{-1}
\Bigl(\sum_{|\beta|=m}\int_{\mathbb{R}^n_+}
|D^\beta u(x)|^p\,x_n^{pa}\,dx\Bigr)^{1/p}
\,\,\,\,\mbox{ for }\,\,u\in C_0^\infty(\mathbb{R}^n_+). 
\end{equation}

\noindent The dual of $V^{m,-a}_{p'}(\mathbb{R}^n_+)$ will be denoted by 
$V^{m,a}_p(\mathbb{R}^n_+)$, where $1/p+1/p'=1$. 

Consider now the operator

\begin{equation}\label{LxDu}
{L}(x,D_x)u:=\sum_{0\leq|\alpha|,|\beta|\leq m}
D^\alpha_x({A}_{\alpha\beta}(x)\,x_n^{|\alpha|+|\beta|-2m} D_x^\beta\,u)
\end{equation}

\noindent where ${A}_{\alpha\beta}$ are $\CC^{l\times l}$-valued functions 
in $L_\infty(\mathbb{R}^n_+)$. We shall use the notation 
$\ring{L}(x,D_x)$ for the principal part of ${L}(x,D_x)$, i.e.

\begin{equation}\label{Lcirc}
\ring{L}(x,D_x)u:=\sum_{|\alpha|=|\beta|=m} 
D^\alpha_x({A}_{\alpha\beta}(x)\,D_x^\beta\,u).
\end{equation}

\subsection{Solvability and regularity result}

\begin{lemma}\label{lem5}
Assume that there exists $\kappa=const>0$ such that the coercivity condition

\begin{equation}\label{B5}
\Re\int_{\mathbb{R}^n_+}\sum_{|\alpha|=|\beta|=m}
\langle A_{\alpha\beta}(x)\,D^\beta u(x),\,D^\alpha u(x)
\rangle_{\mathbb{C}^l}\,dx 
\geq \kappa\sum_{|\gamma|=m}\|D^\gamma\,u\|^2_{L_2(\mathbb{R}^n_+)},
\end{equation}

\noindent holds for all $u\in C^\infty_0(\mathbb{R}^n_+)$, and that 

\begin{equation}\label{B6}
\sum_{|\alpha|=|\beta|= m}\|A_{\alpha\beta}\|_{L_\infty(\mathbb{R}^n_+)} 
\leq \kappa^{-1}.
\end{equation}

{\rm (i)} Let $p\in (1,\infty)$ and $-1/p<a< 1-1/p$ and suppose that 

\begin{equation}\label{E7}
\frac{1}{s(1-s)}\sum_{{|\alpha|+|\beta|<2m}\atop{0\leq |\alpha|,|\beta|\leq m}}
\|A_{\alpha\beta}\|_{{L_\infty}(\mathbb{R}^n_+)} 
+\sum_{|\alpha|=|\beta|=m}[A_{\alpha\beta}]_{{\rm BMO}(\mathbb{R}^n_+)}
\leq \delta,
\end{equation}

\noindent where $\delta$ satisfies

\begin{equation}\label{E8}
\Bigl(pp'+\frac{1}{s(1-s)}\Bigr)\,{\delta}<c(n,m,\kappa)
\end{equation}

\noindent with a sufficiently small constant $c(n,m,\kappa)$ and $s=1-a-1/p$. 
 
Then the operator 

\begin{equation}\label{Liso}
L=L(x,D_x):V_p^{m,a}(\mathbb{R}^n_+)\longrightarrow V_p^{-m,a}(\mathbb{R}^n_+)
\end{equation}

\noindent is an isomorphism.  

{\rm (ii)} Let $p_i\in (1,\infty)$ and $-1/p_i<a_i<1-1/p_i$, where $i=1,2$. 
Suppose that (\ref{E8}) holds with $p_i$ and $s_i=1-a_i-1/p_i$ in place 
of $p$ and $s$. If $u\in V_{p_1}^{m,a_1}(\mathbb{R}^n_+)$ is such that 
${L}u\in V_{p_1}^{-m,a_1}(\mathbb{R}^n_+)
\cap V_{p_2}^{-m,a_2}(\mathbb{R}^n_+)$, 
then $u\in V_{p_2}^{m,a_2}(\mathbb{R}^n_+)$.
\end{lemma}

\noindent{\bf Proof.} The fact that the operator (\ref{Liso}) is continuous 
is obvious. Also, the existence of a bounded inverse ${L}^{-1}$ for $p=2$ 
and $a=0$ follows from (\ref{B5}) and (\ref{E7})-(\ref{E8}) with $p=2$, $a=0$,
which allow us to implement the Lax-Milgram lemma. 

We shall use the notation $\ring{L}_y$ for the operator $\ring{L}(y,D_x)$, 
corresponding to (\ref{Lcirc}) in which the coefficients have been frozen at 
$y\in\mathbb{R}^n_+$, and the notation $G_y$ for the solution operator
for the Dirichlet problem for $\ring{L}_y$ in $\RR^n_+$ with homogeneous 
boundary conditions. Next, given $u\in V_p^{m,a}(\mathbb{R}^n_+)$, set 
$f:=Lu\in V_p^{-m,a}(\mathbb{R}^n_+)$ so that 

\begin{equation}\label{E10}
\left\{
\begin{array}{l}
{L}(x,D)u=f\qquad\qquad\mbox{in}\,\,\,\mathbb{R}^n_+,\\[6pt]
\displaystyle{\frac{\partial^j\,u}{\partial x_n^j}(x',0)=0}
\qquad\qquad {\rm on}\,\,\mathbb{R}^{n-1},\,\,0\leq j\leq m-1. 
\end{array}
\right.
\end{equation}

\noindent Applying the trick used for the first time in {\bf\cite{CFL1}}, 
we may write 

\begin{equation}\label{E12}
u(x)=(G_y f)(x)-(G_{y}(\ring{L}-\ring{L}_y)u)(x)-(G_{y}({L}-\ring{L})u)(x),
\qquad x\in\RR^n_+.
\end{equation}

\noindent We desire to use (\ref{E12}) in order to express $u$ in terms 
of $f$ (cf. (\ref{IntRepFor})-(\ref{SSS}) below) via integral operators
whose norms we can control. 

First, we claim that whenever $|\gamma|=m$, the norm of the operator 

\begin{equation}\label{ItalianTrick}
V_p^{m,a}(\mathbb{R}^n_+)\ni u
\mapsto D_x^\gamma(G_y(\ring{L}-\ring{L}_y)u)(x)\Bigl|_{x=y}\,
\in L_p(\mathbb{R}^n_+,y_n^{ap}\,dy)
\end{equation}

\noindent does not exceed 

\begin{equation}\label{smallCT}
C\,\Bigl(pp'+\frac{1}{s(1-s)}\Bigr)\sum_{|\alpha|=|\beta|=m}
[A_{\alpha\beta}]_{{\rm BMO}(\mathbb{R}^n_+)}.
\end{equation}

\noindent Given the hypotheses under which we operate, 
the expression (\ref{smallCT}) is therefore small if $\delta$ is small. 

In what follows,  
we denote by $G_y(x,z)$ the integral kernel of $G_y$ and integrate 
by parts in order to move derivatives of the form $D_z^\alpha$ with 
$|\alpha|=m$ from $(\ring{L}-\ring{L}_y)u$ onto $G_y(x,z)$ (the absence 
of boundary terms is due to the fact that $G_y(x,\cdot)$ satisfies
homogeneous Dirichlet boundary conditions). 
That (\ref{smallCT}) bounds the norm of (\ref{ItalianTrick}) can now be
seen by combining Theorem~\ref{th1} with (\ref{CC-71}) and 
Proposition~\ref{Prop4}. 

\smallskip

Let $\gamma$ and $\alpha$ be multi-indices with $|\gamma|=m$, 
$|\alpha|\leq m$ and consider the assignment

\begin{equation}\label{oppsi}
\begin{array}{l}
\displaystyle
C_0^\infty(\RR^n_+)\ni\Psi\mapsto\Bigl(D^\gamma_x\int_{\RR^n_+}
G_y(x,z) D^\alpha_z\frac{\Psi(z)}{z_n^{m-|\alpha|}}\,dz\Bigr)\Bigl|_{x=y}.
\end{array}
\end{equation}

\noindent After integrating by parts, the action of this operator can 
be rewritten in the form

\begin{equation}\label{DgammaInt} 
\Bigl(D^\gamma_x\int_{\RR^n_+}\Bigl[  
\Bigl(\frac{-1}{i}\frac{\partial}{\partial z_n}\Bigr)^{m-|\alpha|}
(-D_z)^\alpha G_y(x,z)\Bigr]\Gamma_\alpha(z)\,dz\Bigr)\Bigl|_{x=y},
\end{equation}

\noindent where 

\begin{equation}\label{Gamma-def}
\Gamma_\alpha(z):=
\left\{
\begin{array}{l}
\Psi(z),\qquad\mbox{if }\,\,\,|\alpha|=m,
\\[10pt]
{\displaystyle{\frac{(-1)^{m-|\alpha|}}{(m-|\alpha|-1)!}
\int_{z_n}^\infty (t-z_n)^{m-|\alpha|-1}\frac{\Psi(z',t)}{t^{m-|\alpha|}}\,dt,
\qquad\mbox{if }\,\,\,|\alpha|<m}}. 
\end{array}
\right.
\end{equation}

\noindent Using  Theorem~\ref{th1} along with (\ref{CC-71}) and 
Proposition~\ref{Prop4}, we may therefore conclude that 

\begin{eqnarray}\label{DZGamma}
&& \Bigl\|\Bigl(D^\gamma_x\int_{\RR^n_+} 
\Bigl[\Bigl(\frac{-1}{i}\frac{\partial}{\partial z_n}\Bigr)^{m-|\alpha|} 
(-D_z)^\alpha G_y(x,z)\Bigr]\Gamma_\alpha(z)\,dz\Bigr)\Bigl|_{x=y}
\Bigr\|_{L_p(\mathbb{R}^n_+,\,y_n^{ap}\,dy)}
\nonumber\\[6pt]
&&\qquad\qquad\qquad\qquad\qquad\qquad\qquad
\leq C\,\Bigl(pp'+\frac{1}{s(1-s)}\Bigr)
\|\Gamma_\alpha\|_{L_p(\mathbb{R}^n_+,\,x_n^{ap}\,dx)}.
\end{eqnarray}

\noindent On the other hand, Hardy's inequality gives 

\begin{equation}\label{Cs-1}
\|\Gamma_\alpha\|_{L_p(\mathbb{R}^n_+,\,x_n^{ap}\,dx)} 
\leq\frac{C}{1-s}\,\|\Psi\|_{L_p(\mathbb{R}^n_+,\,x_n^{ap}\,dx)}
\end{equation}

\noindent and, hence, the operator (\ref{oppsi}) can be extended from 
$C_0^\infty(\mathbb{R}^n_+)$ as a bounded operator in 
$L_p(\mathbb{R}^n_+,\,x_n^{ap}\,dx)$ and its norm  is majorized by 

\begin{equation}\label{Bigpp}
\frac{C}{1-s}\,\Bigl(pp'+\frac{1}{s(1-s)}\Bigr).
\end{equation}

Next, given an arbitrary $u\in V_p^{m,a}(\mathbb{R}^n_+)$, we let 
$\Psi=\Psi_{\alpha\beta}$ in (\ref{oppsi}) with 

\begin{equation}\label{Pizz}
\Psi_{\alpha\beta}(z):=z_n^{|\beta|-m}\,A_{\alpha\beta}\,D^\beta u(z),
\qquad |\alpha|+|\beta|<2m,
\end{equation}

\noindent and conclude that the norm of the operator 

\begin{equation}\label{altOp}
V_p^{m,a}(\mathbb{R}^n_+)\ni u\mapsto 
D_x^\gamma\,(G_y\,({L}-\ring{L})u)(x)\Bigl|_{x=y}
\in L_p(\mathbb{R}^n_+,y_n^{ap}\,dy)
\end{equation}

\noindent does not exceed

\begin{equation}\label{notExceed}
\frac{C}{1-s}\,\Bigl(pp'+\frac{1}{s(1-s)}\Bigr)\sum_{{|\alpha|+|\beta|<2m}
\atop{0\leq|\alpha|,|\beta|\leq m}} 
\|A_{\alpha\beta}\|_{{L_\infty}(\mathbb{R}^n_+)}.
\end{equation}

It is well-known (cf. (1.1.10/6) on p.\,22 of {\bf\cite{Maz1}}) 
that any $u\in V^{m,a}_p(\RR^n_+)$ can be represented in the form 

\begin{equation}\label{IntUU}
u=K\{D^\sigma u\}_{|\sigma|=m}
\end{equation}

\noindent where $K$ is a linear operator with the property that 

\begin{equation}\label{K-map}
D^\alpha K:L_p(\mathbb{R}^n_+,x_n^{ap}\,dx)\longrightarrow 
L_p(\mathbb{R}^n_+,x_n^{ap}\,dx)
\end{equation}

\noindent is bounded for every multi-index $\alpha$ with $|\alpha|=m$. 
In particular, by (\ref{4.60}), 

\begin{equation}\label{KDu}
\|K\{D^\sigma u\}_{|\sigma|=m}\|_{V_p^{m,a}(\mathbb{R}^n_+)}\leq C\,s^{-1}
\|\{D^\sigma u\}_{|\sigma|=m}\|_{L_p(\mathbb{R}^n_+,x_n^{ap}\,dx)}.
\end{equation}

At this stage, we transform the identity (\ref{E12}) in the following
fashion. First, we express the two $u$'s occurring inside the Green 
operator $G_y$ in the left-hand side of (\ref{E12}) as in (\ref{IntUU}).  
Second, for each multi-index $\gamma$ with $|\gamma|=m$, we apply $D^\gamma$ 
to both sides of (\ref{E12}) and, finally, set $x=y$. The resulting
identity reads 

\begin{equation}\label{IntRepFor}
\{D^\gamma u\}_{|\gamma|=m}+S\{D^\sigma u\}_{|\sigma|=m}=Q\,f
\end{equation}

\noindent where $Q$ is a bounded operator from $V_p^{-m,a}(\mathbb{R}^n_+)$ 
into $L_p(\mathbb{R}^n_+,\,x_n^{ap}\,dx)$ and $S$ is a linear operator
mapping $L_p(\mathbb{R}^n_+,\,x_n^{ap}\,dx)$ into itself. Furthermore, 
on account of (\ref{ItalianTrick})-(\ref{smallCT}), 
(\ref{altOp})-(\ref{notExceed}) and (\ref{KDu}), we can bound 
$\|S\|_{{\mathfrak L}(L_p(\mathbb{R}^n_+,\,x_n^{ap}\,dx))}$ by 

\begin{equation}\label{SSS}
C\,\Bigl(pp'+\frac{1}{s(1-s)}\Bigr) \Bigl(\sum_{|\alpha|=|\beta|=m}
[A_{\alpha\beta}]_{{\rm BMO}(\mathbb{R}^n_+)} 
+\frac{1}{s(1-s)}
\sum_{{|\alpha|+|\beta|<2m}\atop{0\leq |\alpha|,|\beta|\leq m}} 
\|A_{\alpha\beta}\|_{{L^\infty}(\mathbb{R}^n_+)}\Bigr).
\end{equation}

Owing to (\ref{E7})-(\ref{E8}) and with the integral representation formula 
(\ref{IntRepFor}) and the bound (\ref{SSS}) in hand, 
a Neumann series argument and standard functional analysis allow us to 
simultaneously settle the claims (i) and (ii) in the statement of the lemma.  
\hfill$\Box$
\vskip 0.08in

\section{The Dirichlet problem in a special Lipschitz domain}
\setcounter{equation}{0}

In this section as well as in subsequent ones, we shall work with 
an unbounded domain of the form 

\begin{equation}\label{10.1.26}
G=\{X=(X',X_n)\in\RR^n:\,X'\in\mathbb{R}^{n-1},\,\,X_n>\varphi(X')\},
\end{equation}

\noindent where $\varphi:\RR^{n-1}\to\RR$ is a Lipschitz function.

\subsection{The space ${\rm BMO}(G)$}

The space of functions of bounded mean oscillations in $G$ can be
introduced in a similar fashion to the case $G=\RR^n_+$. Specifically, 
a locally integrable function on $G$ belongs to the space ${\rm BMO}(G)$ if 

\begin{equation}\label{f-BMO}
[f]_{{\rm BMO}(G)}:=\sup\limits_{\{B\}}\meanint_{\!\!\!B\cap G}
\Bigl|f(X)-\meanint_{\!\!\!B\cap G}f(Y)\,dY\Bigr|\,dX<\infty,
\end{equation}

\noindent where the supremum is taken over all balls $B$ centered at 
points in ${\bar G}$. Much as before, 

\begin{equation}\label{10.26}
[f]_{{\rm BMO}(G)}\sim\sup\limits_{\{B\}}\meanint_{\!\!\!B\cap G}
\meanint_{\!\!\!B\cap G}\Bigl|f(X)-f(Y)\Bigr|\,dXdY.
\end{equation}

\noindent This implies the equivalence relation

\begin{equation}\label{1.32}
[f]_{{\rm BMO}(G)}\sim [f\circ\lambda]_{{\rm BMO}(\mathbb{R}^n_+)}
\end{equation}

\noindent for each bi-Lipschitz diffeomorphism $\lambda$ of $\mathbb{R}^n_+$ 
onto $G$. As direct consequences of definitions, we also have 

\begin{eqnarray}\label{1.33}
[\prod_{1\leq j\leq N} f_j]_{{\rm BMO}(G)}
& \leq & c\,\|f\|^{N-1}_{L_\infty(G)}
[f]_{{\rm BMO}(G)},\quad\mbox{where}\,\,f=(f_1,\ldots, f_N),
\\[6pt]
[f^{-1}]_{{\rm BMO}(G)} & \leq & 
c\,\|f^{-1}\|^2_{L_\infty(G)}[f]_{{\rm BMO}(G)}.
\label{1.34}
\end{eqnarray}

\subsection{A bi-Lipschitz map $\lambda:\mathbb{R}^n_+\to G$ and its inverse}

Let $\varphi:\RR^{n-1}\to\RR$ be the Lipschitz function whose graph is 
$\partial G$ and set $M:=\|\nabla\varphi\|_{L_\infty(\RR^{n-1})}$. 
Next, let $T$ be the extension operator defined as in (\ref{10.1.20}) and, 
for a fixed, sufficiently large constant $C>0$, consider the Lipschitz mapping

\begin{equation}\label{lambda}
\lambda:\,\mathbb{R}^n_+\ni(x',x_n)\mapsto (X',X_n)\in\,G
\end{equation}

\noindent defined by the equalities

\begin{equation}\label{10.1.27}
X':=x',\qquad X_n:=C\,M\,x_n+(T\varphi)(x',x_n)
\end{equation}

\noindent (see {\bf\cite{MS}}, \S{6.5.1} and an earlier, less accessible, 
reference {\bf\cite{MS1}}). The Jacobi matrix of $\lambda$ is given by

\begin{equation}\label{matrix}
\lambda'=
\left( 
\begin{array}{cc}
I & 0 \\
\nabla_{x'}(T\varphi) &  CM+\partial (T\varphi)/\partial x_n
\end{array} 
\right)
\end{equation}

\noindent where $I$ is the identity $(n-1)\times(n-1)$-matrix. Since 
$\vert\partial(T\varphi)/\partial x_n\vert\leq cM$ by (\ref{10.1.30}), 
it follows that ${\rm det}\,\lambda'>(C-c)\,M>0$.

Next, thanks to (\ref{Tfi}) and (\ref{lambda})-(\ref{10.1.27}) we have 

\begin{equation}\label{eq-phi}
X_n-\varphi(X')\sim x_n. 
\end{equation}

\noindent Also, based on (\ref{1.8}) we may write

\begin{equation}\label{1.35}
[\lambda']_{{\rm BMO}(\mathbb{R}^n_+)}
\leq c[\nabla\varphi]_{{\rm BMO}(\mathbb{R}^{n-1})}
\end{equation}

\noindent and further, by (\ref{10.1.21}) and (\ref{1.2}), 

\begin{equation}\label{1.4}
\|D^\alpha \lambda'(x)\|_{\mathbb{R}^{n\times n}} 
\leq c(M)\,x_n^{-|\alpha|}[\nabla\varphi]_{{\rm BMO}(\mathbb{R}^{n-1})},
\qquad\forall\,\alpha\,:\,|\alpha|\geq 1. 
\end{equation}

Next, by closely mimicking the proof of Proposition~2.6 from {\bf\cite{MS2}} 
it is possible to show the existence of the inverse Lipschitz  
mapping $\varkappa:=\lambda^{-1}:G\to\RR^n_+$. Owing to (\ref{1.32}), the 
inequality (\ref{1.35}) implies

\begin{equation}\label{1.40}
[\lambda'\circ\varkappa]_{{\rm BMO}(G)}
\leq c[\nabla\varphi]_{{\rm BMO}(\mathbb{R}^{n-1})}.
\end{equation}

\noindent Furthermore, (\ref{1.4}) is equivalent to

\begin{equation}\label{1.41}
\|(D^\alpha\lambda')(\varkappa(X))\|_{\mathbb{R}^{n\times n}} 
\leq c(M,\alpha)(X_n-\varphi(X'))^{-|\alpha|} 
[\nabla\varphi]_{{\rm BMO}(\mathbb{R}^{n-1})},
\end{equation}

\noindent whenever $|\alpha|>0$. Since 
$\varkappa'=(\lambda'\circ\varkappa)^{-1}$ we obtain from 
(\ref{1.34}) and (\ref{1.40})

\begin{equation}\label{1.36}
[\varkappa']_{{\rm BMO}(G)}\leq c[\nabla\varphi]_{{\rm BMO}(\mathbb{R}^{n-1})}.
\end{equation}

\noindent On the other hand, using $\varkappa'=(\lambda'\circ\varkappa)^{-1}$ 
and (\ref{1.41}) one can prove by induction on the order of differentiation 
that

\begin{equation}\label{1.5}
\|D^\alpha\varkappa'(X)\|_{\mathbb{R}^{n\times n}}
\leq c(M,\alpha)\,(X_n-\varphi(X'))^{-|\alpha|}
[\nabla\varphi]_{{\rm BMO}(\mathbb{R}^{n-1})}
\end{equation}

\noindent for all $X\in G$ if $|\alpha|>0$.

\subsection{The space $V_p^{m,a}(G)$}

Analogously to $V_p^{m,a}(\mathbb{R}^n_+)$, we define the weighted Sobolev 
space $V_p^{m,a}(G)$ naturally associated with the norm

\begin{equation}\label{VVV-space}
\|{\mathcal U}\|_{V_p^{m,a}(G)}:=\Bigl(\sum_{0\leq|\gamma|\leq m}\int_{G}
|(X_n-\varphi(X'))^{|\gamma|-m}D^\gamma{\mathcal U}(X)|^p
\,(X_n-\varphi(X'))^{pa}\,dX\Bigr)^{1/p}.
\end{equation}

\noindent Replacing the function $X_n-\varphi(X')$ by either 
$\rho(X):={\rm dist}\,(X,\partial G)$, or by the so-called regularized 
distance function $\rho_{\rm reg}(X)$ 
(defined as on pp.\,170-171 of {\bf\cite{St}}), 
yields equivalent norms on $V_p^{m,a}(G)$. Based on a standard localization 
argument involving a cut-off function vanishing near $\partial G$
(for example, take $\eta(\rho_{\rm reg}/\varepsilon)$ where 
$\eta\in C^\infty_0(\RR)$ satisfies $\eta(t)=0$ for $|t|<1$ and $\eta(t)=1$
for $|t|>2$) one can show that $C^\infty_0(G)$ is 
dense in $V_p^{m,a}(G)$. 

Next, we observe that for each ${\mathcal U}\in C^\infty_0(G)$, 

\begin{equation}\label{equivNr}
C\,s\,\|{\mathcal U}\|_{V_p^{m,a}(G)}\leq \Bigl(\sum_{|\gamma|=m}\int_G 
|D^\gamma {\mathcal U}(X)|^p\,(X_n-\varphi(X'))^{pa}\,dX\Bigr)^{1/p}
\leq\,\|{\mathcal U}\|_{V_p^{m,a}(G)}
\end{equation}

\noindent where, as before, $s=1-a-1/p$. 
Indeed, for each multi-index $\gamma$ with $0\leq|\gamma|\leq m$, 
the one-dimensional Hardy's inequality gives 

\begin{eqnarray}\label{Hardy}
&& \int_{G}|(X_n-\varphi(X'))^{|\gamma|-m}D^\gamma{\mathcal U}(X)|^p\,
(X_n-\varphi(X'))^{pa}\,dX 
\nonumber\\[6pt]
&&\qquad\qquad\qquad\qquad
\leq\bigl({C}/{s}\bigr)^p 
\sum_{|\alpha|=m}\int_G|D^\alpha{\mathcal U}(X)|^p\,(X_n-\varphi(X'))^{pa}\,dX,
\end{eqnarray}

\noindent and the first inequality in (\ref{equivNr}) follows readily from it. 
Also, the second inequality in (\ref{equivNr}) is a trivial consequence 
of (\ref{VVV-space}). 

Going further, we aim to establish that 

\begin{equation}\label{equivNr-bis}
c_1\,\|{u}\|_{V_p^{m,a}(\mathbb{R}^n_+)}
\leq\|{u}\circ\varkappa\|_{V_p^{m,a}(G)}
\leq c_2\,\|{u}\|_{V_p^{m,a}(\mathbb{R}^n_+)},
\end{equation}

\noindent where $c_1$ and $c_2$ do not depend on $p$ and $s$, whereas 
$\varkappa:G\longrightarrow\mathbb{R}^n_+$ is the map introduced in \S{5.2}. 
Clearly, it suffices to prove the upper estimate for 
$\|{u}\circ\varkappa\|_{V_p^{m,a}(G)}$ in (\ref{equivNr-bis}). 
As a preliminary matter, we remark that 

\begin{eqnarray}\label{1.49}
D^\gamma\bigl({u}(\varkappa(X))\bigr) 
& = & \bigl((\varkappa'^*(X)\xi)_{\xi=D}^\gamma\,{u}\bigr)(\varkappa(X))
\nonumber\\[6pt]
&& +\sum_{1\leq|\tau|<|\gamma|}
(D^\tau{u})(\varkappa(X))\sum_\sigma\,c_\sigma\prod_{i=1}^n\prod_j 
D^{\sigma_{ij}}\varkappa_i(X),
\end{eqnarray}

\noindent where

\begin{equation}\label{1.50}
\sum_{i,j}\sigma_{ij}=\gamma,\quad |\sigma_{ij}|\geq 1,\quad 
\sum_{i,j}(|\sigma_{ij}|-1)=|\gamma|-|\tau|.
\end{equation}

\noindent In turn, (\ref{1.49})-(\ref{1.50}) and (\ref{1.5}) 
allow us to conclude that 

\begin{equation}\label{DUU}
|D^\gamma\bigl({u}(\varkappa(X))\bigr)|\leq c\sum_{1\leq |\tau|\leq |\gamma|}
x_n^{|\tau|-|\gamma|}\,|D^\tau{u}(x)|, 
\end{equation}

\noindent which, in view of (\ref{eq-phi}), yields the desired conclusion. 

Finally, we set 

\begin{equation}\label{dual-VG}
V^{-m,a}_p(G):=\Bigl(V^{m,-a}_{p'}(G)\Bigr)^*.
\end{equation}

\noindent where, as usual, $p'=p/(p-1)$.

\bigskip

\subsection{Solvability and regularity result for the Dirichlet 
problem in the domain $G$}

Let us consider the differential operator 

\begin{equation}\label{E4}
{\mathcal L}\,{\mathcal U}
={\mathcal L}(X,D_X)\,{\mathcal U}=\sum_{|\alpha|=|\beta|=m} 
D^\alpha(\mathfrak{A}_{\alpha\beta}(X)\,D^\beta{\mathcal U}),\qquad X\in G,
\end{equation}

\noindent whose matrix-valued coefficients satisfy

\begin{equation}\label{E4a}
\sum_{|\alpha|=|\beta|=m}\|\mathfrak{A}_{\alpha\beta}\|_{L_\infty(G)} 
\leq \kappa^{-1}.
\end{equation}
 
\noindent This operator generates the sesquilinear form  
${\mathcal L}(\cdot,\cdot):V_p^{m,a}(G)\times V_{p'}^{m,-a}(G)\to\CC$
where $p'$ is the conjugate exponent of $p$, defined by 

\begin{equation}\label{LUV}
{\mathcal L}({\mathcal U},{\mathcal V}):=\sum_{|\alpha|=|\beta|=m}
\int_G\langle\mathfrak{A}_{\alpha\beta}(X)\,
D^\beta{\mathcal U}(X),\,D^\alpha{\mathcal V}(X)\rangle\,dX. 
\end{equation}

\noindent We assume that the inequality

\begin{equation}\label{B25}
\Re\,{\mathcal L}({\mathcal U},{\mathcal U})
\geq \kappa\sum_{|\gamma|=m}\|D^\gamma\,{\mathcal U}\|^2_{L_2(G)}
\end{equation}

\noindent holds for all ${\mathcal U}\in V_2^{m,0}(G)$. 

\begin{lemma}\label{lem5a}
{\rm (i)} Let $p\in (1,\infty)$, $-1/p<a<1-1/p$ and $s:=1-a-1/p$. 
Suppose that

\begin{equation}\label{E7a}
[\nabla\varphi]_{{\rm BMO}(\mathbb{R}^{n-1})}
+\sum_{|\alpha|=|\beta|=m}[{\mathfrak A}_{\alpha\beta}]_{{\rm BMO}(G)}
\leq \delta,
\end{equation}

\noindent where $\delta$ satisfies

\begin{equation}\label{E8b}
\Bigl(pp'+\frac{1}{s(1-s)}\Bigr)\,\frac{\delta}{s(1-s)}
<C(n,m,\kappa,\|\nabla\varphi\|_{L_\infty(\mathbb{R}^{n-1})})
\end{equation}

\noindent with a sufficiently small constant $C$, independent of $p$ and $s$. 
In the case $m=1$ the factor $\delta/s(1-s)$ in {\rm (\ref{E8b})} 
can be replaced by $\delta$. 
 
Then the operator 

\begin{equation}\label{Liso-2}
{\mathcal L}(X,D_X):V_p^{m,a}(G)\longrightarrow V_p^{-m,a}(G)
\end{equation}

\noindent is an isomorphism.  

{\rm (ii)} Let $p_i\in (1,\infty)$ and $-1/p_i<a_i<1-1/p_i$, where $i=1,2$. 
Suppose that {\rm (\ref{E8})} holds with $p_i$ and $s_i=1-a_i-1/p_i$ in place 
of $p$ and $s$. If ${\mathcal U}\in V_{p_1}^{m,a_1}(G)$ and 
${\mathcal L}\,{\mathcal U}\in V_{p_1}^{-m,a_1}(G)\cap V_{p_2}^{-m,a_2}(G)$, 
then ${\mathcal U}\in V_{p_2}^{m,a_2}(G)$.
\end{lemma}

\noindent{\bf Proof.} We shall extensively use the flattening mapping 
$\lambda$ and its inverse studied in \S{5.2}. The assertions (i) and (ii)  
will follow directly from Lemma~\ref{lem5} as soon as we show that 
the operator $L$ defined in $\mathbb{R}^n_+$ by

\begin{equation}\label{E20}
L({\mathcal U}\circ\lambda):=({\mathcal L}\,{\mathcal U})\circ\lambda
\end{equation}

\noindent satisfies all the hypotheses in that lemma. The sesquilinear form 
corresponding to the operator $L$ will be denoted by $L(u,v)$.

Set $u(x):={\mathcal U}(\lambda(x))$, $v(x):={\mathcal V}(\lambda(x))$ and
note that the identity (\ref{1.49}) implies

\begin{equation}\label{E40}
D^\beta {\mathcal U}(X)=\bigl((\varkappa'^*(\lambda(x))\xi)_{\xi=D}^\beta\,{u}
\bigr)(x)+\sum_{1\leq|\tau|<|\beta|}K_{\beta\tau}(x)\,x_n^{|\tau|-|\beta|} 
D^\tau u(x),
\end{equation}
\begin{equation}\label{E41}
D^\alpha {\mathcal V}(X)
=\bigl((\varkappa'^*(\lambda(x))\xi)_{\xi=D}^\alpha\,{v}
\bigr)(x)+\sum_{1\leq|\tau|<|\alpha|}K_{\alpha\tau}(x)\,x_n^{|\tau|-|\alpha|} 
D^\tau v(x),
\end{equation}

\noindent where, thanks to (\ref{1.5}), the coefficients $K_{\gamma\tau}$ 
satisfy

\begin{equation}\label{E42}
\|K_{\gamma\tau}\|_{L_\infty(\mathbb{R}^n_+)} 
\leq c[\nabla\varphi]_{{\rm BMO}(\mathbb{R}^{n-1})}.
\end{equation}

\noindent Plugging (\ref{E40}) and (\ref{E41}) into the definition of 
${\mathcal L}({\mathcal U},{\mathcal V})$, we arrive at

\begin{equation}\label{LUV2}
{\mathcal L}({\mathcal U},{\mathcal V})=L_0(u,v)
+\sum_{{1\leq|\alpha|,|\beta|\leq m}\atop{|\alpha|+|\beta|<2m}}
\int_{\mathbb{R}^n_+}\langle {A}_{\alpha\beta}(x)\,x_n^{|\alpha|+|\beta|-2m} 
D^\beta\,u(x),\,D^\alpha v(x)\rangle\,dx,
\end{equation}

\noindent where

\begin{equation}\label{Lzero}
L_0(u,v)=\sum_{|\alpha|=|\beta|=m}\int_{\mathbb{R}^n_+} 
\langle (\mathfrak{A}_{\alpha\beta}\circ\lambda)
((\varkappa'^*\circ\lambda)\xi)_{\xi=D}^\beta\,{u},\, 
((\varkappa'^*\circ\lambda)\xi)_{\xi=D}^\alpha\,{v}\rangle
\,{\rm det}\,\lambda'\,dx.
\end{equation}

\noindent It follows from (\ref{E40})-(\ref{E42}) that the coefficient 
matrices $A_{\alpha\beta}$ obey 

\begin{equation}\label{E43}
\sum_{{1\leq|\alpha|,|\beta|\leq m}\atop{|\alpha|+|\beta|<2m}}
\|A_{\alpha\beta}\|_{L_\infty(\mathbb{R}^n_+)} 
\leq {c}{\kappa}^{-1}\,[\nabla\varphi]_{{\rm BMO}(\mathbb{R}^{n-1})},
\end{equation}

\noindent where $c$ depends on $m$, $n$, and 
$\|\nabla\varphi\|_{L_\infty(\mathbb{R}^{n-1})}$. 
We can write the form $L_0(u,v)$ as

\begin{equation}\label{Lzero-uv}
\sum_{|\alpha| =|\beta| =m}\int_{\mathbb{R}^n_+}
\langle {A}_{\alpha\beta}(x)\,D^\beta u(x),\,D^\alpha v(x)\rangle\,dx
\end{equation}

\noindent where the coefficient matrices $A_{\alpha\beta}$ are given by

\begin{equation}\label{Aalbet}
A_{\alpha\beta}={\rm det}\,\lambda'\,\sum_{|\gamma|=|\tau|=m}
P_{\alpha\beta}^{\gamma\tau}(\varkappa'\circ\lambda)
({\mathfrak A}_{\gamma\tau}\circ\lambda),
\end{equation}

\noindent for some scalar homogeneous polynomials 
$P_{\alpha\beta}^{\gamma\tau}$ of the elements of the matrix 
$\varkappa'(\lambda(x))$ with ${\rm deg}\,P_{\alpha\beta}^{\gamma\tau}=2m$. 
In view of (\ref{1.33})-(\ref{1.36}),  

\begin{equation}\label{E44}
\sum_{|\alpha|=|\beta|=m}[A_{\alpha\beta}]_{{\rm BMO}(\mathbb{R}^n_+)} 
\leq c\Bigl(\kappa^{-1}[\nabla\varphi]_{{\rm BMO}(\mathbb{R}^{n-1})} 
+\sum_{|\alpha|=|\beta|=m}[{\mathfrak A}_{\alpha\beta}]_{{\rm BMO}(G)}\Bigr),
\end{equation}

\noindent where $c$ depends on $n$, $m$, and 
$\|\nabla\varphi\|_{L_\infty(\mathbb{R}^{n-1})}$.

By (\ref{E43}) 

\begin{equation}\label{L-L}
|L(u,u)-L_0(u,u)|\leq c\,\delta\|u\|_{V_2^{m,0}(\mathbb{R}^n_+)}^2
\end{equation}

\noindent and, therefore,

\begin{equation}\label{ReLzero}
\Re\,L_0(u,u)\geq\Re\,{\mathcal L}({\mathcal U},{\mathcal U})
-c\,\delta\|u\|^2_{V_2^{m,0}(\mathbb{R}^n_+)}.
\end{equation}

\noindent Using (\ref{B25}) and the equivalence

\begin{equation}\label{norm-U}
\|{\mathcal U}\|_{V_2^{m,0}(G)}\sim\|u\|_{V_2^{m,0}(\mathbb{R}^n_+)}
\end{equation}

\noindent (cf. the discussion in \S{5.3}), we arrive at (\ref{B5}). 
Thus, all conditions of Lemma~\ref{lem5} hold and the result follows. 
The improvement of (\ref{E8b}) for $m=1$ mentioned in the statement 
(i) holds because in this case $L=L_0$.
\hfill$\Box$
\vskip 0.08in

\section{Dirichlet problem in a bounded Lipschitz domain}
\setcounter{equation}{0}

\subsection{Preliminaries}

Let $\Omega$ be a {\it bounded Lipschitz domain} in $\RR^n$ which means
(cf. {\bf\cite{St}}, p.\,189)
that there exists a finite open covering $\{{\mathcal O}_j\}_{1\leq j\leq N}$
of $\partial\Omega$ with the property that, for every $j\in\{1,...,N\}$,
${\mathcal O}_j\cap\Omega$ coincides with the portion of ${\mathcal O}_j$ 
lying in the over-graph of a Lipschitz function $\varphi_j:\RR^{n-1}\to\RR$
(where $\RR^{n-1}\times\RR$ is a new system of coordinates obtained from 
the original one via a rigid motion). We then define the 
{\it Lipschitz constant} of a bounded Lipschitz domain $\Omega\subset\RR^n$ as

\begin{equation}\label{Lip-ct}
\inf\,\Bigl(\max\{\|\nabla\varphi_j\|_{L_\infty(\RR^{n-1})}:\,1\leq j\leq N\}
\Bigr)
\end{equation}

\noindent where the infimum is taken over all possible families 
$\{\varphi_j\}_{1\leq j\leq N}$ as above. 

It is a classical result that the surface measure $d\sigma$ is
well-defined and that there exists an outward pointing normal vector 
$\nu$ at almost every point on $\partial\Omega$. 
 
We denote by $\rho(X)$ the distance from $X\in\RR^n$ to $\partial\Omega$
and, for $p$, $a$ and $m$ as in (\ref{indices}), introduce the weighted 
Sobolev space $V^{m,a}_p(\Omega)$ naturally associated with the norm

\begin{equation}\label{normU2}
\|{\mathcal U}\|_{V_p^{m,a}(\Omega)}
:=\Bigl(\sum_{0\leq |\beta|\leq m}\int_{\Omega}
|\rho(X)^{|\beta|-m} D^\beta{\mathcal U}(X)|^p\,\rho(X)^{pa}\,dX\Bigr)^{1/p}. 
\end{equation}

\noindent One can check the equivalence of the norms

\begin{equation}\label{equiv-Nr2}
\|{\mathcal U}\|_{V_p^{m,a}(\Omega)}\sim 
\|\rho_{\rm reg}^a\,{\mathcal U}\|_{V_p^{m,0}(\Omega)}, 
\end{equation}

\noindent where $\rho_{\rm reg}(X)$ stands for the regularized distance 
from $X$ to $\partial\Omega$ (in the sense of Theorem~2, p.\,171 in 
{\bf\cite{St}}). It is also easily proved that 
$C_0^\infty(\Omega)$ is dense in $V^{m,a}_p(\Omega)$ and that 

\begin{equation}\label{sumDU}
\|{\mathcal U}\|_{V_p^{m,a}(\Omega)}\sim 
\Bigl(\sum_{|\beta|=m}\int_\Omega|D^\beta{\mathcal U}(X)|^p\,\rho(X)^{pa}\,dX
\Bigr)^{1/p}
\end{equation}

\noindent uniformly for ${\mathcal U}\in C_0^\infty(\Omega)$. 
As in (\ref{dual-VG}), we set 

\begin{equation}\label{dual-V}
V^{-m,a}_p(\Omega):=\Bigl(V^{m,-a}_{p'}(\Omega)\Bigr)^*.
\end{equation}

Let us fix a Cartesian coordinates system and consider 
the differential operator 

\begin{equation}\label{E444}
{\mathcal A}\,{\mathcal U}={\mathcal A}(X,D_X)\,{\mathcal U}
:=\sum_{|\alpha|=|\beta|=m}D^\alpha({\mathcal A}_{\alpha\beta}(X)
\,D^\beta{\mathcal U}),\qquad X\in\Omega,
\end{equation}

\noindent with measurable $l\times l$ matrix-valued coefficients. 
The corresponding sesquilinear form will be denoted by 
${\mathcal A}({\mathcal U},{\mathcal V})$. 
Similarly to (\ref{E4a}) and (\ref{B25}) we impose the conditions

\begin{equation}\label{E4b}
\sum_{|\alpha|=|\beta|=m}\|{\mathcal A}_{\alpha\beta}\|_{L_\infty(\Omega)} 
\leq \kappa^{-1}
\end{equation}

\noindent and

\begin{equation}\label{B25b}
\Re\,{\mathcal A}({\mathcal U},{\mathcal U})\geq\kappa\sum_{|\gamma|=m}
\|D^\gamma\,{\mathcal U}\|^2_{L_2(\Omega)}\quad\mbox{for all}\,\,\,
{\mathcal U}\in V_2^{m,0}(G).
\end{equation}

\subsection{Interior regularity of solutions}

\begin{lemma}\label{lem2}
Let $\Omega\subset\RR^n$ be a bounded Lipschitz domain. 
Pick two functions ${\mathcal H},{\mathcal Z}\in C^\infty_0(\Omega)$ 
such that ${\mathcal H}\,{\mathcal Z}={\mathcal H}$, and assume that 

\begin{equation}\label{E5}
\sum_{|\alpha|=|\beta|=m}[{\mathcal A}_{\alpha\beta}]_{{\rm BMO}(\Omega)}
\leq\delta
\end{equation}

\noindent where
 
\begin{equation}\label{E6}
\delta\leq \frac{c(m,n,\kappa)}{p\,p'}
\end{equation}

\noindent with a sufficiently small constant $c(m,n,\kappa)>0$.

If ${\mathcal U}\in W_q^m(\Omega,loc)$ for a certain $q<p$ and 
${\mathcal A}\,{\mathcal U}\in W_p^{-m}(\Omega,loc)$, then 
${\mathcal U}\in W_p^m(\Omega,loc)$ and

\begin{equation}\label{E3}
\|{\mathcal H}\,{\mathcal U}\|_{W_p^m(\Omega)}
\leq C\,(\|{\mathcal H}\,{\mathcal A}(\cdot,D)\,
{\mathcal U}\|_{W_p^{-m}(\Omega)}
+\|{\mathcal Z}\,{\mathcal U}\|_{W_q^m(\Omega)}).
\end{equation}
\end{lemma}

\noindent{\bf Proof.} We start with a trick applied in {\bf\cite{CFL1}} under 
slightly  different circumstances. We shall use the notation ${\mathcal A}_Y$ 
for the operator ${\mathcal A}(Y, D_X)$, where $Y\in\Omega$ and the notation 
$\Phi_Y$ for  a fundamental solution of ${\mathcal A}_Y$ in $\mathbb{R}^n$. 
Then, with star denoting the convolution product, 

\begin{equation}\label{HUPhi}
{\mathcal H}\,{\mathcal U}
+\Phi_Y\ast({\mathcal A}-{\mathcal A}_Y)({\mathcal H}\,{\mathcal U}) 
=\Phi_Y\ast({\mathcal H}\,{\mathcal A}{\mathcal U})
+\Phi_Y\ast ([{\mathcal A},{\mathcal H}]({\mathcal Z}{\mathcal U}))
\end{equation}

\noindent and, consequently, for each multi-index $\gamma$, $|\gamma|=m$, 

\begin{eqnarray}\label{DHU}
&& D^\gamma({\mathcal H}\,{\mathcal U})+\sum_{|\alpha|=|\beta|=m} 
D^{\alpha+\gamma}\Phi_Y\ast\bigl(({\mathcal A}_{\alpha\beta}
-{\mathcal A}_{\alpha\beta}(Y)) D^\beta({\mathcal H}\,{\mathcal U})\bigr)
\nonumber\\[6pt]
&& \qquad\quad
=D^\gamma\Phi_Y\ast({\mathcal H}\,{\mathcal A}{\mathcal U})
+D^\gamma\Phi_Y\ast([{\mathcal A},{\mathcal H}]({\mathcal Z}{\mathcal U})).
\end{eqnarray}

\noindent Writing this equation at the point $Y$ and using (\ref{b1}), 
we obtain

\begin{eqnarray}\label{E3a}
&& (1-C\,pp'\delta)
\sum_{|\gamma|=m}\|D^\gamma({\mathcal H}\,{\mathcal U})\|_{L_p(\Omega)}
\nonumber\\[6pt]
&& \qquad\quad
\leq C(p,\kappa)
(\|{\mathcal H}\,{\mathcal A}{\mathcal U}\|_{W_p^{-m}(\Omega)} 
+\|[{\mathcal A},{\mathcal H}]({\mathcal Z}{\mathcal U})\|_{W_p^{-m}(\Omega)}).
\end{eqnarray}

\noindent Let $p'<n$. We have for every ${\mathcal V}\in \ring W^m_p(\Omega)$

\begin{eqnarray}\label{AHZUV}
&& \Bigl|\int_\Omega\langle[{\mathcal A},{\mathcal H}]
({\mathcal Z}{\mathcal U}),{\mathcal V}\rangle\,dX\Bigr|
=|{\mathcal A}({\mathcal H}{\mathcal Z}{\mathcal U},{\mathcal V})
-{\mathcal A}({\mathcal Z}{\mathcal U},{\mathcal H}{\mathcal V})| 
\nonumber\\[6pt]
&& \qquad\qquad\quad
\leq c(\|{\mathcal Z}{\mathcal U}\|_{W_p^{m-1}(\Omega)} 
\|{\mathcal V}\|_{W_{p'}^m(\Omega)}
+\|{\mathcal Z}{\mathcal U}\|_{W_{\frac{pn}{n+p}}^{m} 
(\Omega)}\|{\mathcal V}\|_{W_{\frac{p'n}{n-p'}}^{m-1}(\Omega)}).
\end{eqnarray}

\noindent By Sobolev's theorem

\begin{equation}\label{ZZU}
\|{\mathcal Z}{\mathcal U}\|_{W_p^{m-1}(\Omega)}
\leq c\,\|{\mathcal Z}{\mathcal U}\|_{W_{\frac{pn}{n+p}}^{m}(\Omega)} 
\end{equation}

\noindent and

\begin{equation}\label{VWpm}
\|{\mathcal V}\|_{W_{\frac{p'n}{n-p'}}^{m-1}(\Omega)} 
\leq c\,\|{\mathcal V}\|_{W_{p'}^m(\Omega)}.
\end{equation}

\noindent Therefore,

\begin{equation}\label{AHZ}
\Bigl|\int_\Omega\langle
[{\mathcal A},{\mathcal H}]({\mathcal Z}{\mathcal U}),{\mathcal V}\rangle 
\,dX\Bigr|
\leq c\,\|{\mathcal Z}{\mathcal U}\|_{W_{\frac{pn}{n+p}}^{m}(\Omega)} 
\|{\mathcal V}\|_{W_{p'}^m(\Omega)}
\end{equation}

\noindent which is equivalent to the inequality

\begin{equation}\label{AHZ3}
\|[{\mathcal A},{\mathcal H}]({\mathcal Z}{\mathcal U})\|_{W_p^{-m}(\Omega)}
\leq c\,\|{\mathcal Z}{\mathcal U}\|_{W_{\frac{pn}{n+p}}^{m} (\Omega)}. 
\end{equation}

\noindent In the case $p'\geq n$, the same argument leads to a 
similar inequality, where $pn/(n+p)$ is replaced by $1+\varepsilon$ with 
an arbitrary $\varepsilon>0$ for $p'>n$ and $\varepsilon=0$ for $p'=n$.

Now, (\ref{E3}) follows from (\ref{E3a}) if $p'\geq n$ and $p'<n$, 
$q\geq pn/(n+p)$. In the remaining case the goal is achieved by iterating 
this argument finitely many times. 
\hfill$\Box$
\vskip 0.08in

\begin{corollary}\label{cor3}
Let $p\geq 2$ and suppose that {\rm (\ref{E5})} and {\rm (\ref{E6})} hold. 
If ${\mathcal U}\in W_2^m(\Omega,loc)$ and 
${\mathcal A}\,{\mathcal U}\in W_p^{-m}(\Omega,loc)$, then 
${\mathcal U}\in W_p^m(\Omega,loc)$ and

\begin{equation}\label{E3-aaa}
\|{\mathcal H}\,{\mathcal U}\|_{W_p^m(\Omega)}
\leq C\,(\|{\mathcal Z}\,{\mathcal A}(\cdot,D)\, 
{\mathcal U}\|_{W_p^{-m}(\Omega)}
+\|{\mathcal Z}\,{\mathcal U}\|_{W_2^{m-1}(\Omega)}).
\end{equation}
\end{corollary}

\noindent{\bf Proof.} Let ${\mathcal Z}_0$ denote a real-valued function 
in $C^\infty_0(\Omega)$ such that ${\mathcal H}{\mathcal Z}_0={\mathcal H}$ 
and ${\mathcal Z}_0{\mathcal Z}={\mathcal Z}_0$. By (\ref{E3})

\begin{equation}\label{E3b}
\|{\mathcal H}\,{\mathcal U}\|_{W_p^m(\Omega)}
\leq C\,(\|{\mathcal H}\,{\mathcal A}(\cdot,D)\, 
{\mathcal U}\|_{W_p^{-m}(\Omega)}
+\|{\mathcal Z}_0\,{\mathcal U}\|_{W_2^{m}(\Omega)})
\end{equation}

\noindent and it follows from (\ref{B25b}) that

\begin{equation}\label{ZUW}
\|{\mathcal Z}_0\,{\mathcal U}\|^2_{W_2^m(\Omega)}\leq c\kappa^{-1}
\Re\,{\mathcal A}({\mathcal Z}_0{\mathcal U},{\mathcal Z}_0{\mathcal U}).
\end{equation}

\noindent Furthermore,

\begin{equation}\label{AZ5}
|{\mathcal A}({\mathcal Z}_0{\mathcal U},{\mathcal Z}_0{\mathcal U})
-{\mathcal A}({\mathcal U},{\mathcal Z}_0^2{\mathcal U})|\leq c\kappa^{-1}
\|{\mathcal Z}{\mathcal U}\|_{W_2^{m-1}(\Omega)}\,
\|{\mathcal Z}_0{\mathcal U}\|_{W_2^{m}(\Omega)}.
\end{equation}

\noindent Hence

\begin{equation}\label{ZU6}
\|{\mathcal Z}_0\,{\mathcal U}\|^2_{W_2^m(\Omega)}\leq c\kappa^{-1}
(\|{\mathcal Z}\,{\mathcal A}{\mathcal U}\|_{W_2^{-m}(\Omega)}\, 
\|{\mathcal Z}_0^2\,{\mathcal U}\|_{W_2^m(\Omega)}
+\kappa^{-1}\|{\mathcal Z}{\mathcal U}\|_{W_2^{m-1}(\Omega)}\,
\|{\mathcal Z}_0{\mathcal U}\|_{W_2^{m}(\Omega)})
\end{equation}

\noindent and, therefore, 

\begin{equation}\label{Zu7}
\|{\mathcal Z}_0\,{\mathcal U}\|_{W_2^m(\Omega)}
\leq c\kappa^{-1}(\|{\mathcal Z}\,{\mathcal A} 
{\mathcal U}\|_{W_2^{-m}(\Omega)}\,
+\kappa^{-1}\|{\mathcal Z}{\mathcal U}\|_{W_2^{m-1}(\Omega)}). 
\end{equation}

\noindent Combining this inequality with (\ref{E3b}) we arrive 
at (\ref{E3-aaa}). 
\hfill$\Box$
\vskip 0.08in

\subsection{Invertibility of 
${\mathcal A}:V_p^{m,a}(\Omega)\longrightarrow V_p^{-m,a}(\Omega)$}

Recall the infinitesimal mean oscillations as defined in (\ref{e60}).

\begin{theorem}\label{th1a}
Let  $1<p<\infty$, $0<s<1$, and $a=1-s-1/p$. Furthermore, let $\Omega$ be a 
bounded Lipschitz domain in $\mathbb{R}^n$. Suppose that the differential
operator ${\mathcal A}$ is as in {\rm \S{6.1}} and that, in addition, 

\begin{equation}\label{E16}
\sum_{|\alpha|=|\beta|=m}\{{\mathcal A}_{\alpha\beta}\}_{{\rm Osc}(\Omega)}
+\{\nu\}_{{\rm Osc}(\partial\Omega)}\leq\delta,
\end{equation}

\noindent where

\begin{equation}\label{E17}
\Bigl(pp'+\frac{1}{s(1-s)}\Bigr)\frac{\delta}{s(1-s)}\leq c
\end{equation}

\noindent for a sufficiently small constant $c>0$ independent of $p$ and $s$. 
In the case $m=1$ the factor $\delta/s(1-s)$ in {\rm (\ref{E17})} can 
be replaced by $\delta$.

Then the operator

\begin{equation}\label{cal-L}
{\mathcal A}:V_p^{m,a}(\Omega)\longrightarrow V_p^{-m,a}(\Omega)
\end{equation}

\noindent is an isomorphism. 
\end{theorem}

\noindent{\bf Proof.} We shall proceed in a series a steps starting with

\vskip 0.08in
(i) {\it The construction of the auxiliary domain $G$ and 
operator ${\mathcal L}$}.

\noindent Let $\varepsilon$ be small enough so that 

\begin{equation}\label{m1}
\sum_{|\alpha|=|\beta|=m}\meanint_{\!\!\!B_r\cap \Omega}
\meanint_{\!\!\!B_r\cap\Omega}
|{\mathcal A}_{\alpha\beta}(X)-{\mathcal A}_{\alpha\beta}(Y)|\,dXdY
\leq 2\delta
\end{equation}

\noindent for all balls in $\{B_r\}_\Omega$ with radii $r<\varepsilon$ and 

\begin{equation}\label{m2}
\meanint_{\!\!\!B_r\cap\partial\Omega}\meanint_{\!\!\!B_r\cap\partial\Omega}\,
\Bigl|\nu(X)-\nu(Y)\,\Bigr|\,d\sigma_Xd\sigma_Y\leq 2\delta
\end{equation}

\noindent for all balls in $\{B_r\}_{\partial\Omega}$ with radii 
$r<\varepsilon$. 

We fix a ball $B_\varepsilon$ in $\{B_\varepsilon\}_{\partial\Omega}$ 
and assume without loss of generality that, in a suitable system of 
Cartesian coordinates,  

\begin{equation}\label{newGGG}
\Omega\cap B_\varepsilon=\{X=(X',X_n)\in B_\varepsilon:\,X_n>\varphi(X')\}
\end{equation}

\noindent for some Lipschitz function $\varphi:\RR^{n-1}\to\RR$. 
Consider now the unique cube $Q(\varepsilon)$ (relative to this system of 
coordinates) which is inscribed in $B_\varepsilon$ and denote its projection 
onto $\mathbb{R}^{n-1}$ by $Q'(\varepsilon)$. 
Since $\nabla\varphi=-\nu'/\nu_n$, it follows from (\ref{m2}) that 

\begin{equation}\label{m3}
\meanint_{\!\!\!B'_r}\meanint_{\!\!\!B'_r}\,
\Bigl|\nabla\varphi(X')-\nabla\varphi(Y')\,\Bigr|\,dX'dY'\leq c(n)\,\delta,
\end{equation}

\noindent where $B'_r=B_r\cap \mathbb{R}^{n-1}$, $r<\varepsilon$. Let us 
retain the notation $\varphi$ for the mirror extension of the function 
$\varphi$ from $Q'(\varepsilon)$ onto $\mathbb{R}^{n-1}$.  

We extend ${\mathcal A}_{\alpha\beta}$ from $Q(\varepsilon)\cap\Omega$ onto 
$Q(\varepsilon)\backslash\Omega$ by setting

\begin{equation}\label{A=A}
{\mathcal A}_{\alpha\beta}(X)
:={\mathcal A}_{\alpha\beta}(X',-X_n+2\varphi(X')),
\qquad X\in Q(\varepsilon)\backslash\Omega, 
\end{equation}
 
\noindent and we shall use the notation ${\mathfrak A}_{\alpha\beta}$ for 
the periodic extension of ${\mathcal A}_{\alpha\beta}$ from $Q(\varepsilon)$ 
onto $\mathbb{R}^n$. 

Consistent with the earlier discussion in \S{5}, we shall denote the special 
Lipschitz domain $\{X=(X',X_n):\,X'\in\mathbb{R}^{n-1},
\,X_n>\varphi(X')\}$ by $G$. One can easily see that, owing to 
$2\varepsilon n^{-1/2}$-periodicity of $\varphi$ and 
${\mathcal A}_{\alpha\beta}$, 

\begin{equation}\label{SumA}
\sum_{|\alpha|=|\beta|=m}[{\mathcal A}_{\alpha\beta}]_{{\rm BMO}(G)}
+[\nabla\varphi]_{{\rm BMO}(\mathbb{R}^{n-1})}\leq c(n)\,\delta.
\end{equation}

\noindent Now, with the operator ${\mathcal A}(X,D_X)$ in $\Omega$, we 
associate an auxiliary operator ${\mathcal L}(X,D_X)$ in $G$ 
given by (\ref{E4}). 

\vskip 0.08in
(ii) {\it Uniqueness.} 

\noindent Assuming that ${\mathcal U}\in V_p^{m,a}(\Omega)$ satisfies 
${\mathcal L}\,{\mathcal U}=0$ in $\Omega$, we shall show that 
${\mathcal U}\in V_2^{m,0}(\Omega)$. This will imply that ${\mathcal U}=0$ 
which proves the injectivity of the operator (\ref{cal-L}). 

To this end, pick a function ${\mathcal H}\in C_0^\infty(Q(\varepsilon))$ 
and write 
${\mathcal L}({\mathcal H}\,{\mathcal U})
=[{\mathcal L},\,{\mathcal H}]\,{\mathcal U}$. 
Also, fix a small $\theta>0$ and select a smooth function $\Lambda$ on 
$\mathbb{R}^1_+$, which is identically $1$ on $[0,1]$ and which 
vanishes identically on $(2,\infty)$. Then by (ii) in Lemma~\ref{lem5a},

\begin{equation}\label{E19}
{\mathcal L}({\mathcal H}\,{\mathcal U})-[{\mathcal L},\,{\mathcal H}]\,
(\Lambda(\rho_{\rm reg}/\theta)\,{\mathcal U})
\in V_2^{-m,0}(G)\cap V_p^{-m,a}(G).
\end{equation}

\noindent Note that the operator 

\begin{equation}\label{LH-1}
[{\mathcal  L},\,{\mathcal H}]\rho_{\rm reg}^{-1}: V_p^{m,a}(G) 
\longrightarrow V_p^{-m,a}(G)
\end{equation}

\noindent is bounded and that the norm of the multiplier 
$\rho_{\rm reg}\,\Lambda(\rho_{\rm reg}/\theta)$ in $V_p^{m,a}(G)$ 
is $O(\theta)$. Moreover, the same is true for $p=2$ and $a=0$. 

The inclusion (\ref{E19}) can be written in the form

\begin{equation}\label{E21}
{\mathcal L}({\mathcal H}\,{\mathcal U})
+{\mathcal M}({\mathcal Z}\,{\mathcal U})\in V_p^{-m,a}(G)\cap V_2^{-m,0}(G),
\end{equation}

\noindent where ${\mathcal Z}\in C^\infty_0(\mathbb{R}^n)$, 
${\mathcal Z}\,{\mathcal H}={\mathcal H}$ and ${\mathcal M}$ is a 
linear operator mapping 

\begin{equation}\label{V2Vp}
V_p^{m,a}(G)\to V_p^{-m,a}(G)\quad {\rm and}\quad  
V_2^{m,0}(G)\to V_2^{-m,0}(G)
\end{equation}

\noindent with both norms of order $O(\theta)$. 

Select a finite covering of $\overline{\Omega}$ by cubes $Q_j(\varepsilon)$ 
and let $\{{\mathcal H}_j\}$ be a smooth partition of unity subordinate to 
$\{Q_j(\varepsilon)\}$. 
Also, let ${\mathcal Z}_j\in C_0^\infty(Q_j(\varepsilon))$ 
be such that ${\mathcal H}_j{\mathcal Z}_j={\mathcal H}_j$. 
By $G_j$ we denote the special Lipschitz domain generated by the 
cube $Q_j(\varepsilon)$ as in part (i) of the present proof. 
The corresponding operators ${\mathcal L}$ and ${\mathcal M}$ will 
be denoted by ${\mathcal L}_j$ and ${\mathcal M}_j$, respectively. 
It follows from (\ref{E21}) that 

\begin{equation}\label{Hu8}
{\mathcal H}_j\,{\mathcal U}
+\sum_k({\mathcal L}_j^{-1}\,{\mathcal M}_j\,{\mathcal Z}_j
\,{\mathcal Z}_k)({\mathcal H}_k\,{\mathcal U})
\in V_p^{m,a}(\Omega)\cap V_2^{m,0}(\Omega).
\end{equation}

\noindent Taking into account that the norms of the matrix operator 
${\mathcal L}_j\,{\mathcal M}_j\,{\mathcal Z}_j\,{\mathcal Z}_k$ 
in the spaces $V_p^{m,a}(\Omega)$ and $V_2^{m,0}(\Omega)$ are $O(\theta)$, 
we may take $\theta>0$ small enough and obtain 
${\mathcal H}_j\,{\mathcal U}\in V_2^{m,0}(\Omega)$, i.e. 
${\mathcal U}\in V_2^{m,0}(\Omega)$. Therefore, 
${\mathcal L}:V_p^{m,a}(\Omega)\to V_p^{-m,a}(\Omega)$ is injective.

\vskip 0.08in
(iii) {\it A priori estimate}. 

\noindent Let $p\geq 2$ and assume that ${\mathcal U}\in V_p^{m,a}(\Omega)$. 
Referring to Corollary~\ref{cor3} and arguing as in part (ii) of the present 
proof, we arrive at the equation

\begin{equation}\label{m32}
{\mathcal H}_j\,{\mathcal U}
+\sum_k({\mathcal L}_j^{-1}\,{\mathcal M}_j\,{\mathcal Z}_j\,
{\mathcal Z}_k)({\mathcal H}_k\,{\mathcal U})={\mathcal F},
\end{equation}

\noindent whose right-hand side satisfies

\begin{equation}\label{FVO}
\|{\mathcal F}\|_{V_p^{m,a}(\Omega)}\leq c
(\|{\mathcal A}\,{\mathcal U}\|_{V_p^{-m,a}(\Omega)}
+\|{\mathcal U}\|_{W_2^{m-1}(\omega)}),
\end{equation}

\noindent for some domain $\omega$ with ${\overline\omega}\subset\Omega$. 
Since the $V_p^{m,a}(\Omega)$-norm of the sum in (\ref{m32}) does not 
exceed $ C\theta\|{\mathcal U}\|_{V_p^{m,a}(\Omega)}$, we obtain the estimate 

\begin{equation}\label{m32-bis}
\|{\mathcal U}\|_{V_p^{m,a}(\Omega)}
\leq c\,(\|{\mathcal A}\,{\mathcal U}\|_{V_p^{-m,a}(\Omega)}
+\|{\mathcal U} \|_{W_2^{m-1}(\omega)}).
\end{equation}

\vskip 0.08in
(iv) {\it End of proof.} 

\noindent Let $p\geq 2$. The range of the operator 
${\mathcal A}:V_p^{m,a}(\Omega)\to V_p^{-m,a}(\Omega)$ is 
closed by (\ref{m32}) and the compactness of the restriction operator: 
$V_p^{m,a}(\Omega) \to W_2^{m-1}(\omega)$. Since the coefficients of the 
adjoint operator ${\mathcal L}^*$ satisfy the same conditions as those of 
${\mathcal L}$, the operator 
${\mathcal L}^*: V_{p'}^{m,a}(\Omega)\to V_{p'}^{-m,-a}(\Omega)$ 
is injective. Therefore, we conclude that 
${\mathcal L}:V_{p}^{m,a}(\Omega)\to V_{p}^{-m,-a}(\Omega)$ is surjective. 
Being also injective, ${\mathcal L}$ is isomorphic if $p\geq 2$. 
Hence ${\mathcal L}^*$ is isomorphic for $p'\leq 2$. 
This means that ${\mathcal L}$ is isomorphic for $p\leq 2$. The result follows.
\hfill$\Box$
\vskip 0.08in

\subsection{Traces and extensions}

Let $\Omega\subset\RR^n$ be a bounded Lipschitz domain and, for $m\in\NN$ 
as well as $1<p<\infty$ and $-1/p<a<1-1/p$, consider a new space, 
$W_p^{m,a}(\Omega)$, consisting of functions ${\mathcal U}\in L_p(\Omega,loc)$ 
with the property that $\rho^{a}D^\alpha{\mathcal U}\in L_p(\Omega)$ for all 
multi-indices $\alpha$ with $|\alpha|=m$. We equip $W_p^{m,a}(\Omega)$ 
with the norm 

\begin{equation}\label{newW}
\|{\mathcal U}\|_{W_p^{m,a}(\Omega)}
:=\sum_{|\alpha|=m}\|D^\alpha{\mathcal U}\|_{L_p(\Omega,\,\rho(X)^{ap}\,dX)} 
+\|{\mathcal U}\|_{L_p(\omega)},
\end{equation}

\noindent where $\omega$ is an open non-empty domain, 
$\overline{\omega}\subset\Omega$. An equivalent norm is given by
the expression in (\ref{W-Nr}). We omit the standard proof of the fact that 

\begin{equation}\label{dense}
C^\infty({\overline\Omega})\hookrightarrow W_p^{m,a}(\Omega)
\quad\mbox{densely}. 
\end{equation}

Recall that for $p\in(1,\infty)$ and $s\in(0,1)$ the Besov space 
$B_p^s(\partial\Omega)$ is then defined via the requirement (\ref{Bes-xxx}). 
If we introduce the $L_p$-modulus of continuity

\begin{equation}\label{omega-p}
\omega_p(f,t):=\Bigl(\int\!\!\!\!\!\!\int
\limits_{{|X-Y|<t}\atop{X,Y\in\partial\Omega}}
|f(X)-f(Y)|^p\,d\sigma_Xd\sigma_Y\Bigr)^{1/p},
\end{equation}

\noindent then 

\begin{equation}\label{Eqv}
\|f\|_{B_p^s(\partial\Omega)}\sim\|f\|_{L_p(\partial\Omega)}
+\left(\int_0^\infty\frac{\omega_p(f,t)^p}{t^{n-1+ps}}dt\right)^{1/p},
\end{equation}

\noindent uniformly for $f\in B^s_p(\partial\Omega)$. 

The nature of our problem requires that we work with Besov spaces
(defined on Lipschitz boundaries) which exhibit a higher order
of smoothness. In accordance with {\bf\cite{JW}}, we now make the following
definition. 

\begin{definition}\label{def1} 
For $p\in(1,\infty)$, $m\in\NN$ and $s\in(0,1)$, define the (higher order) 
Besov space $\dot{B}^{m-1+s}_p(\partial\Omega)$ as the collection of all 
finite families $\dot{f}=\{f_\alpha\}_{|\alpha|\leq m-1}$ of functions defined 
on $\partial\Omega$ with the following property. For each multi-index
$\alpha$ of length $\leq m-1$ let 

\begin{equation}\label{reminder}
R_\alpha(X,Y):=f_\alpha(X)-\sum_{|\beta|\leq m-1-|\alpha|}\frac{1}{\beta!}\,
f_{\alpha+\beta}(Y)\,(X-Y)^\beta,\qquad X,Y\in\partial\Omega,
\end{equation}

\noindent and consider the $L_p$-modulus of continuity 

\begin{equation}\label{rem-Rr}
r_\alpha(t):=\Bigl(\int\!\!\!\!\!\!\int
\limits_{{|X-Y|<t}\atop{X,Y\in\partial\Omega}}
|R_\alpha(X,Y)|^p\,d\sigma_Xd\sigma_Y\Bigr)^{1/p}. 
\end{equation}

\noindent Then 

\begin{equation}\label{Bes-Nr}
\|\dot{f}\|_{\dot{B}^{m-1+s}_p(\partial\Omega)}
:=\sum_{|\alpha|\leq m-1}\|f_\alpha\|_{L_p(\partial\Omega)}
+\sum_{|\alpha|\leq m-1}\Bigl(\int_0^\infty
\frac{r_\alpha(t)^p}{t^{p(m-1+s-|\alpha|)+n-1}}\,dt\Bigr)^{1/p}<\infty.
\end{equation}
\end{definition}

For further reference we note here that for each fixed $\kappa>0$, 
an equivalent norm is obtained by replacing $r_\alpha(t)$ by 
$r_\alpha(\kappa\,t)$ in (\ref{Bes-Nr}). 
Also, when $m=1$, the above definition agrees with (\ref{Bes-xxx}), 
thanks to (\ref{Eqv}).

A few notational conventions which make the exposition more transparent are
as follows. Given a family of functions $\{f_\alpha\}_{|\alpha|\leq m-1}$
on $\partial\Omega$ and $X\in\Omega$, $Y,Z\in\partial\Omega$, set 

\begin{equation}\label{PPP}
\begin{array}{l}
{\displaystyle{
P_\alpha(X,Y):=\sum_{|\beta|\leq m-1-|\alpha|}\frac{1}{\beta!}\,
f_{\alpha+\beta}(Y)\,(X-Y)^\beta,\qquad \forall\,\alpha\,:\,|\alpha|\leq m-1,}}
\\[25pt]
P(X,Y):=P_{(0,...,0)}(X,Y),
\end{array}
\end{equation}

\noindent so that 

\begin{equation}\label{PR-0}
R_\alpha(Y,Z)=f_\alpha(Y)-P_\alpha(Y,Z),
\qquad\forall\,\alpha\,:\,|\alpha|\leq m-1,
\end{equation}

\noindent and the following elementary identities hold for each 
multi-index $\alpha$ of length $\leq m-1$:

\begin{eqnarray}\label{PR}
D^{\beta}_XP_{\alpha}(X,Y) & = & P_{\alpha+\beta}(X,Y),
\qquad|\beta|\leq m-1-|\alpha|,
\nonumber\\[6pt]
P_{\alpha}(X,Y) -P_{\alpha}(X,Z) & =& 
\sum_{|\beta|\leq m-1-|\alpha|}R_{\alpha+\beta}(Y,Z)\frac{(X-Y)^\beta}{\beta!}.
\end{eqnarray}

\noindent See, e.g., p.\,177 in {\bf\cite{St}} for the last formula. 

\begin{lemma}\label{trace-1}
For each $1<p<\infty$, $-1/p<a<1-1/p$ and $s=1-a-1/p$, the trace operator

\begin{equation}\label{TR-1}
{\rm Tr}:W^{1,a}_p(\Omega)\longrightarrow B^s_p(\partial\Omega)
\end{equation}

\noindent is well-defined, linear, bounded, onto and has $V^{1,a}_p(\Omega)$ 
as its null-space. Furthermore, there exists a linear, continuous mapping

\begin{equation}\label{Extension}
{\mathcal E}:B^s_p(\partial\Omega)\longrightarrow W^{1,a}_p(\Omega),
\end{equation}

\noindent called extension operator, such that ${\rm Tr}\circ{\mathcal E}=I$
(i.e., the operator (\ref{TR-1}) has a bounded, linear right-inverse). 
\end{lemma}

\noindent{\bf Proof.} By a standard argument involving a smooth partition 
of unity it suffices to deal with the case when $\Omega$ is the domain 
lying above the graph of a Lipschitz function $\varphi:\RR^{n-1}\to\RR$. 
Composing with the bi-Lipschitz homeomorphism 
$\RR^n_+\ni(X',X_n)\mapsto(X',\varphi(X')+X_n)\in\Omega$ further reduces
matters to the case when $\Omega=\RR^n_+$, in which situation the claims 
in the lemma have been proved in {\bf\cite{Usp}}. 
\hfill$\Box$
\vskip 0.08in

We need to establish an analogue of Lemma~\ref{trace-1} for 
higher smoothness spaces. While for $\Omega=\RR^n_+$ this has been done by 
Uspenski\u{\i} in {\bf\cite{Usp}}, the flattening argument used in 
Lemma~\ref{trace-1} is no longer effective in this context. 
Let us also mention here that a result similar in spirit, valid for any 
Lipschitz domain $\Omega$ but with $B^{m-1+s+1/p}(\Omega)$ in place of
$W^{m,a}_p(\Omega)$ (cf. (\ref{incls}) for the relationship between 
these spaces) has been proved by A.\,Jonsson and H.\,Wallin in {\bf\cite{JW}} 
(in fact, in this latter context, these authors have dealt with much more
general sets than Lipschitz domains). The result which serves 
our purposes is as follows. 

\begin{proposition}\label{trace-2}
Let $1<p<\infty$, $-1/p<a<1-1/p$, $s=1-a-1/p\in(0,1)$ and $m\in\NN$. 
Define the {\rm higher} {\rm order} trace operator

\begin{equation}\label{TR-11}
{\rm tr}_{m-1}:W^{m,a}_p(\Omega)\longrightarrow
\dot{B}^{m-1+s}_p(\partial\Omega)
\end{equation}

\noindent by setting 

\begin{equation}\label{Tr-DDD}
{\rm tr}_{m-1}\,\,{\mathcal U}
:=\Bigl\{i^{|\alpha|}\,{\rm Tr}\,[D^\alpha\,{\mathcal U}]\Bigr\}
_{|\alpha|\leq m-1},
\end{equation}

\noindent where the traces in the right-hand side are taken in the sense
of Lemma~\ref{trace-1}. Then (\ref{TR-11})-(\ref{Tr-DDD}) is a
a well-defined, linear, bounded operator, which is onto and has 
$V^{m,a}_p(\Omega)$ as its null-space. Moreover, it has a bounded, linear 
right-inverse, i.e. there exists a linear, continuous operator 

\begin{equation}\label{Ext-222}
{\mathcal E}:\dot{B}^{m-1+s}_p(\partial\Omega)
\longrightarrow W^{m,a}_p(\Omega)
\end{equation}

\noindent such that 

\begin{equation}\label{Ext-333}
\dot{f}=\{f_\alpha\}_{|\alpha|\leq m-1}\in\dot{B}^{m-1+s}(\partial\Omega)
\Rightarrow i^{|\alpha|}\,{\rm Tr}\,[D^\alpha({\mathcal E}\,\dot{f})]=f_\alpha,
\quad\forall\,\alpha\,:\,|\alpha|\leq m-1. 
\end{equation}
\end{proposition}

In order to facilitate the exposition, we isolate a couple of preliminary
results prior to the proof of Proposition~\ref{trace-2}. 

\begin{lemma}\label{Lemma-R}
Assume that $\varphi:\RR^{n-1}\to\RR$ is a Lipschitz function 
and define $\Phi:\RR^{n-1}\to\partial\Omega\hookrightarrow\RR^n$ by setting 
$\Phi(X'):=(X',\varphi(X'))$ at each $X'\in\RR^{n-1}$. Define the Lipschitz 
domain $\Omega$ as $\{X=(X',X_n)\in\RR^n:\,X_n>\varphi(X')\}$
and, for some fixed $m\in\NN$, $p\in(1,\infty)$ and $s\in(0,1)$  
consider a system of functions $f_\alpha\in B^s_p(\partial\Omega)$, 
$\alpha\in\NN_0^n$, $|\alpha|\leq m-1$, with the property that 

\begin{equation}\label{B-CC}
\frac{\partial}{\partial X_k}[f_\alpha(\Phi(X'))]
=\sum_{j=1}^{n}f_{\alpha+e_j}(\Phi(X'))\partial_k\Phi_j(X'),
\qquad 1\leq k\leq n-1,
\end{equation}

\noindent for each multi-index $\alpha$ of length $\leq m-2$, 
where $\{e_j\}_j$ is the canonical orthonormal basis in $\RR^n$. 
Finally, for each $l\in\{1,...,m-1\}$ introduce 
$\Delta_l:=\{(t_1,...,t_{l}):\,0\leq t_{l}\leq\cdots\leq t_1\leq 1\}$, 
and define $R_\alpha(X,Y)$ as in (\ref{reminder}). Then if $\alpha$ is an 
arbitrary multi-index of length $\leq m-2$ and $r:=m-1-|\alpha|$, 
the following identity holds: 

\begin{eqnarray}\label{RRR=id}
&& R_\alpha(\Phi(X'),\Phi(Y'))
\nonumber\\[6pt]
&&\quad
=\sum_{(j_1,...,j_{r})\in\{1,...,n\}^{r}}\Bigl\{\int_{\Delta_{r}}\Bigl[
f_{\alpha+e_{j_1}+\cdots+e_{j_r}}(\Phi(Y'+t_{r}(X'-Y')))
-f_{\alpha+e_{j_1}+\cdots+e_{j_r}}(\Phi(Y'))\Bigr]
\nonumber\\[6pt]
&&
\qquad\times\prod_{k=1}^{r}\,\nabla\Phi_{j_k}(Y'+t_{k}(X'-Y'))\cdot(X'-Y')
\,dt_{r}\cdots dt_1\Bigr\},\qquad X',\,Y'\in\RR^{n-1}.
\end{eqnarray}
\end{lemma}

\noindent{\bf Proof.} We shall show that for any system of functions
$\{f_\alpha\}_{|\alpha|\leq m-1}$ which satisfies (\ref{B-CC}), any 
multi-index $\alpha\in\NN_0^n$ with $|\alpha|\leq m-2$ and any $l\in\NN$
with $l\leq r:=m-1-|\alpha|$, there holds

\begin{eqnarray}\label{FF=id}
&& f_\alpha(\Phi(X'))-\sum_{|\beta|\leq l}\frac{1}{\beta!}
f_{\alpha+\beta}(\Phi(Y'))(\Phi(X')-\Phi(Y'))^\beta 
\nonumber\\[6pt]
&&\quad=\sum_{(j_1,...,j_{l})\in\{1,...,n\}^{l}}\Bigl\{\int_{\Delta_l}
\Bigl[f_{\alpha+e_{j_1}+\cdots+e_{j_l}}(\Phi(Y'+t_{l}(X'-Y')))
-f_{\alpha+e_{j_1}+\cdots+e_{j_l}}(\Phi(Y'))\Bigr]
\nonumber\\[6pt]
&&\qquad\qquad\qquad\qquad\quad
\times\prod_{k=1}^{l}\,\nabla\Phi_{j_k}(Y'+t_{k}(X'-Y'))\cdot(X'-Y')
\,dt_{l}\cdots dt_1\Bigr\}.
\end{eqnarray}

\noindent Clearly, (\ref{RRR=id}) follows from (\ref{reminder}) 
and (\ref{FF=id}) by taking $l:=r$. 

In order to justify (\ref{FF=id}) we proceed by induction on $l$. 
Concretely, when $l=1$ we may write, based on (\ref{B-CC}) and 
the Fundamental Theorem of Calculus,  

\begin{eqnarray}\label{RRR=0}
&& f_{\alpha}(\Phi(X'))-f_{\alpha}(\Phi(Y'))
-\sum_{j=1}^n f_{\alpha+e_j}(\Phi(Y'))(\Phi_j(X')-\Phi_j(Y'))
\nonumber\\[6pt]
&&\quad =\int_0^1\frac{d}{dt}\Bigl[f_\alpha(\Phi(Y'+t(X'-Y')))\Bigr]\,dt
\nonumber\\[6pt]
&& \qquad\qquad\quad-\sum_{j=1}^n f_{\alpha+e_j}(\Phi(Y'))
\int_0^1\frac{d}{dt}\Bigl[\Phi_j(Y'+t(X'-Y'))\Bigr]\,dt
\nonumber\\[6pt]
&&\quad = \sum_{j=1}^n\Bigl\{\int_0^1\Bigl[f_{\alpha+e_j}(Y'+t(X'-Y'))
-f_{\alpha+e_j}(\Phi(Y'))\Bigr]
\nonumber\\[6pt]
&&\qquad\qquad\times\nabla\Phi_j(Y'+t(X'-Y'))\cdot(X'-Y')\,dt\Bigr\},
\end{eqnarray}

\noindent as wanted. To prove the version of (\ref{FF=id}) when $l$ 
is replaced by $l+1$ we split the sum in the left-hand side of (\ref{FF=id}),
written for $l+1$ in place of $l$, according to whether $|\beta|\leq l$
or $|\beta|=l+1$ and denote the expressions created in this fashion 
by $S_1$ and $S_2$, respectively. Next, based on 
(\ref{B-CC}) and the Fundamental Theorem of Calculus, we write  

\begin{eqnarray}\label{FTC}
&& f_{\alpha+e_{j_1}+\cdots+e_{j_l}}(\Phi(Y'+t_{l}(X'-Y')))
-f_{\alpha+e_{j_1}+\cdots+e_{j_l}}(\Phi(Y'))
\\[6pt]
&&
=\sum_{i=1}^n\int_0^{t_{l}}f_{\alpha+e_{j_1}+\cdots+e_{j_l}+e_i}
(\Phi(Y'+t_{l+1}(X'-Y')))
\nabla\Phi_{i}(Y'+t_{l+1}(X'-Y'))\cdot(X'-Y')\,dt_{l+1}
\nonumber
\end{eqnarray}

\noindent and use the induction hypothesis to conclude that 

\begin{eqnarray}\label{FF=id-2}
&& S_1=\sum_{(j_1,...,j_{l+1})\in\{1,...,n\}^{l+1}}
\Bigl\{\int_{\Delta_{l+1}}
f_{\alpha+e_{j_1}+\cdots+e_{j_{l+1}}}(\Phi(Y'+t_{l+1}(X'-Y')))
\nonumber\\[6pt]
&&\qquad\qquad
\times\prod_{k=1}^{l+1}\,\nabla\Phi_{j_k}(Y'+t_{k}(X'-Y'))\cdot(X'-Y')
\,dt_{l+1}\cdots dt_1\Bigr\}.
\end{eqnarray}

\noindent Thus, if 

\begin{equation}\label{F=Psi}
F_j(t):=\Phi_j(Y'+t(X'-Y'))-\Phi_j(Y'),\qquad 1\leq j\leq n,
\end{equation}

\noindent we may express $S_1$ in the form 

\begin{eqnarray}\label{another}
&& S_1=\sum_{(j_1,...,j_{l+1})\in\{1,...,n\}^{l+1}}\Bigl\{\int_{\Delta_{l+1}}
\Bigl[f_{\alpha+e_{j_1}+\cdots+e_{j_{l+1}}}(\Phi(Y'+t_{l+1}(X'-Y')))
-f_{\alpha+e_{j_1}+\cdots+e_{j_{l+1}}}(\Phi(Y'))\Bigr]
\nonumber\\[6pt]
&&\qquad\qquad\qquad\qquad
\times\prod_{k=1}^{l+1}\,\nabla\Phi_{j_k}(Y'+t_{k}(X'-Y'))\cdot(X'-Y')
\,dt_{l+1}\cdots dt_1\Bigr\}
\\[6pt]
&&\qquad\quad +\sum_{(j_1,...,j_{l+1})\in\{1,...,n\}^{l+1}}
f_{\alpha+e_{j_1}+\cdots+e_{j_{l+1}}}(\Phi(Y'))\int_{\Delta_{l+1}}\,
\prod_{k=1}^{l+1}\,F'_{j_k}(t_{k})\,dt_{l+1}\cdots dt_1.
\nonumber
\end{eqnarray}

Note that the first double sum above corresponds precisely to 
the expression in the right-hand side of (\ref{FF=id}) written with $l$
replaced by $l+1$. Our proof of (\ref{FF=id}) by induction is therefore 
complete as soon as we show that for each multi-index $\beta$ of length $l+1$,

\begin{equation}\label{S2}
\sum_{{(j_1,...,j_{l+1})\in\{1,...,n\}^{l+1}}\atop
{e_{j_1}+\cdots+e_{j_{l+1}}=\beta}}
\int_{\Delta_{l+1}}\,
\prod_{k=1}^{l+1}\,F'_{j_k}(t_{k})\,dt_{l+1}\cdots dt_1
=\frac{1}{\beta!}(\Phi(X')-\Phi(Y'))^\beta.
\end{equation}

\noindent In turn, this is going to be a consequence of a general identity,
to the effect that 

\begin{equation}\label{Fprime}
\sum_{{(j_1,...,j_{l})\in\{1,...,n\}^{l}}
\atop{e_{j_1}+\cdots+e_{j_{l}}=\beta}}
\int_0^{t_0}\int_0^{t_1}\cdots\int_0^{t_{l-1}}
\prod_{k=1}^{l}\,F'_{j_k}(t_{k})\,dt_{l}\cdots dt_1
=\frac{1}{\beta!}F(t_0)^\beta,
\end{equation}

\noindent for any Lipschitz function $F=(F_1,...,F_n):[0,1]\to\CC^n$
with $F(0)=0$, any point $t_0\in[0,1]$ any $l\in\NN$ and any multi-index 
$\beta$ of length $l$. Of course, the case most relevant for our purposes 
is when the $F_j$'s are as in (\ref{F=Psi}), $t_0=1$ and when $l$ is replaced 
by $l+1$, but the above formulation is best suited for proving (\ref{Fprime}) 
via induction on $l$. Indeed, the case $l=1$ is immediate from the 
Fundamental Theorem of Calculus and to pass from $l$ to $l+1$ it suffices 
to show that the two sides of (\ref{Fprime}) have the same derivative with 
respect to $t_0$. The important observation in carrying out the latter 
step is that the derivative of the left-hand side of (\ref{Fprime})
with respect to $t_0$ is an expression to which the current induction 
hypothesis is readily applicable. 
This justifies (\ref{Fprime}) and completes the proof of (\ref{RRR=id}). 
\hfill$\Box$
\vskip 0.08in

\begin{corollary}\label{Cor-R}
Under the assumptions of Lemma~\ref{Lemma-R}, 
for each multi-index of length $\leq m-2$ the following estimate holds

\begin{equation}\label{RRR=est}
\Bigl(\int_0^\infty
\frac{r_\alpha(t)^p}{t^{p(m-1+s-|\alpha|)+n-1}}\,dt\Bigr)^{1/p}
\leq C\sum_{|\gamma|=m-1}\|f_\gamma\|_{B^s_p(\partial\Omega)},
\end{equation}

\noindent where the constant $C$ depends only on $n$, $p$, 
$s$ and $\|\nabla\varphi\|_{L_\infty(\RR^{n-1})}$. 
\end{corollary}

\noindent{\bf Proof.} The identity (\ref{RRR=id}) gives 

\begin{eqnarray}\label{pw-est-R}
&& |R_\alpha(\Phi(X'),\Phi(Y'))|
\\[6pt]
&&\qquad\qquad
\leq C|X'-Y'|^{m-1-|\alpha|}
\sum_{|\gamma|=m-1}\int_0^1
|f_{\gamma}(\Phi(Y'+\tau(X'-Y')))-f_{\gamma}(\Phi(Y'))|\,d\tau
\nonumber
\end{eqnarray}

\noindent for each $X',Y'\in\RR^{n-1}$, where the constant $C$ 
depends only on $n$ and $\|\nabla\Phi\|_{L_\infty}$ which, in turn, is 
controlled in terms of $\|\nabla \varphi\|_{L_\infty}$. 
Given an arbitrary $t>0$ we now integrate the $p$-th power of 
both sides in (\ref{pw-est-R}) for 
$X',Y'\in\partial\RR^{n-1}$ subject to $|\Phi(X')-\Phi(Y')|<t$.
Using Fubini's Theorem and making the change of variables
$Z':=Y'+\tau(X'-Y')$ we obtain, after noticing that $|Z'-Y'|\leq \tau t$, 

\begin{eqnarray}\label{r-vs-om}
r_\alpha(t)^p & \leq & C\,t^{p(m-1-|\alpha|)}\sum_{|\gamma|=m-1}
\int\!\!\!\!\!\!\int\limits_{{X',Y'\in\RR^{n-1}}\atop{|X'-Y'|< c\,t}}
\int_0^1|f_{\gamma}(\Phi(Y'+\tau(X'-Y')))-f_{\gamma}(\Phi(Y'))|^p
\,d\tau\,dX'dY'
\nonumber\\[6pt]
&\leq & C\,t^{p(m-1-|\alpha|)}\sum_{|\gamma|=m-1}\int_0^1
\int\!\!\!\!\!\!\int\limits_{{Z',Y'\in\RR^{n-1}}\atop{|Z'-Y'|<c\,\tau t}}
|f_{\gamma}(\Phi(Z'))-f_{\gamma}(\Phi(Y'))|^p\,dZ'dY'd\tau
\nonumber\\[6pt]
&\leq & C\,t^{p(m-1+s-|\alpha|)+n-1}
\sum_{|\gamma|=m-1}\int_0^1
\frac{\omega_p(f_{\gamma},\,c\,\tau t)^p}{\tau^{n-1}t^{ps+n-1}}\,d\tau.
\end{eqnarray}

\noindent Consequently, 

\begin{eqnarray}\label{r-om-bis}
\int_0^\infty
\frac{r_\alpha(t)^p}{t^{p(m-1+s-|\alpha|)+n-1}}\,dt
& \leq & C\sum_{|\gamma|=m-1}\int_0^\infty \int_0^1
\frac{\omega_p(f_{\gamma},\,c\,\tau t)^p}{\tau^{n-1}t^{ps+n-1}}\,d\tau dt
\nonumber\\[6pt]
& \leq & C\sum_{|\gamma|=m-1}\Bigl(\int_0^\infty
\frac{\omega_p(f_{\gamma},r)^p}{r^{ps+n-1}}\,dr\Bigr)
\Bigl(\int_0^1\frac{1}{\tau^{1-sp}}\,d\tau\Bigr)
\nonumber\\[6pt]
& \leq & C\sum_{|\gamma|=m-1}\int_0^\infty
\frac{\omega_p(f_{\gamma},t)^p}{t^{ps+n-1}}\,dt,
\end{eqnarray}

\noindent after making the change of variables $r:=c\,\tau t$ in the 
second step. With this in hand, the estimate (\ref{RRR=est}) follows 
by virtue of (\ref{Eqv}). 
\hfill$\Box$
\vskip 0.08in

After this preamble, we are in a position to present the 

\vskip 0.08in
\noindent{\bf Proof of Proposition~\ref{trace-2}.} We divide the proof 
into a series of steps, starting with 

\vskip 0.08in
\noindent{\it Step I: The well-definiteness of trace.}  
Let ${\mathcal U}$ be an arbitrary function in $W^{m,a}_p(\Omega)$ and set 

\begin{equation}\label{ff-aa}
f_\alpha:=i^{|\alpha|}\,{\rm Tr}\,[D^\alpha\,{\mathcal U}],\qquad
\forall\,\alpha\,:\,|\alpha|\leq m-1.
\end{equation}
 
\noindent It follows from Lemma~\ref{trace-1} that these trace functions 
are well-defined and, in fact, 

\begin{equation}\label{falpha}
\sum_{|\alpha|\leq m-1}\|f_\alpha\|_{B^s_p(\partial\Omega)}
\leq C\|\,{\mathcal U}\|_{W^{m,a}_p(\Omega)}.
\end{equation}

In order to prove that $\dot{f}:=\{f_\alpha\}_{|\alpha|\leq m-1}$ belongs to 
$\dot{B}^{m-1+s}_p(\partial\Omega)$, let $R_\alpha(X,Y)$ and $r_\alpha(t)$ 
be as in (\ref{reminder})-(\ref{rem-Rr}). Our goal is to show that for 
every multi-index $\alpha$ with $|\alpha|\leq m-1$, 

\begin{equation}\label{B-alpha}
\Bigl(\int_0^\infty
\frac{r_\alpha(t)^p}{t^{p(m-1+s-|\alpha|)+n-1}}\,dt\Bigr)^{1/p}
\leq C\|\,{\mathcal U}\|_{W^{m,a}_p(\Omega)}. 
\end{equation}

\noindent To this end, we first observe that if $|\alpha|=m-1$ then the 
expression in the left-hand side of (\ref{B-alpha}) is majorized by 
$C\Bigl(\int_0^\infty\omega_p(f_\alpha,t)^p/t^{ps+n-1}\,dt\Bigr)^{1/p}$
which, by (\ref{Eqv}) and (\ref{falpha}), is indeed 
$\leq C\|\,{\mathcal U}\|_{W^{m,a}_p(\Omega)}$. To treat the case 
when $|\alpha|<m-1$ we assume that $\Omega$ is locally
represented as $\{X:\,X_n>\varphi(X')\}$ for some Lipschitz function 
$\varphi:\RR^{n-1}\to\RR$ and, as before, set $\Phi(X'):=(X',\varphi(X'))$, 
$X'\in\RR^{n-1}$. Then (\ref{B-CC}) holds, thanks to (\ref{ff-aa}),
for every multi-index $\alpha$ of length $\leq m-2$. 
Consequently, Corollary~\ref{Cor-R} applies and, in concert with 
(\ref{falpha}), yields (\ref{B-alpha}). This proves that the operator 
(\ref{TR-11})-(\ref{Tr-DDD}) is well-defined and bounded. 

\vskip 0.08in
\noindent{\it Step II: The extension operator.}  
We introduce a co-boundary operator 
${\mathcal E}$ which acts on $\dot{f}=\{f_\alpha\}_{|\alpha|\leq m-1}\in
\dot{B}^{m-1+s}_p(\partial\Omega)$ according to 

\begin{equation}\label{def-Ee}
({\mathcal E}\dot{f})(X)
=\int_{\partial\Omega}{\mathcal K}(X,Y)\,P(X,Y)\,d\sigma_Y,
\qquad X\in\Omega,
\end{equation}

\noindent where $P(X,Y)$ is the polynomial associated with $\dot{f}$ as 
in (\ref{PPP}). The integral kernel ${\mathcal K}$ is assumed to satisfy

\begin{eqnarray}\label{ker-prp}
&& \int_{\partial\Omega}{\mathcal K}(X,Y)\,d\sigma_Y=1
\qquad\mbox{for all }\,X\in\Omega,
\\[6pt]
&&
|D^\alpha_X{\mathcal K}(X,Y)|\leq c_\alpha\,\rho(X)^{1-n-|\alpha|},
\quad\forall\,X\in\Omega,\,\,\forall\,Y\in\partial\Omega,
\label{more-Kp}
\end{eqnarray}

\noindent where $\alpha$ is an arbitrary multi-index, and 

\begin{equation}\label{last-Kp}
{\mathcal K}(X,Y)=0\quad\mbox{if }\,\,|X-Y|\geq 2\rho(X).
\end{equation}

\noindent One can take, for instance, the kernel

\begin{equation}\label{K-def}
{\mathcal K}(X,Y):=\eta\left(\frac{X-Y}{\varkappa\rho_{\rm reg}(X)}\right)
\left(\int_{\partial\Omega}\eta\left(\frac{X-Z}{\varkappa\rho_{\rm reg}(X)} 
\right)d\sigma_Z\right)^{-1},
\end{equation}

\noindent where $\eta\in C^\infty_0(B_2)$, $\eta=1$ on $B_1$, $\eta\geq 0$ and 
$\varkappa$ is a positive constant depending on the Lipschitz constant 
of $\partial\Omega$. Here, as before, $\rho_{\rm reg}(X)$ stands for 
the regularized distance from $X$ to $\partial\Omega$.  

For each $X\in\Omega$ and $Z\in\partial\Omega$ and for every multi-index 
$\gamma$ with $|\gamma|=m$ we then obtain

\begin{equation}\label{N0}
D^\gamma{\mathcal E}\dot{f}(X)
=\sum_{{\alpha+\beta=\gamma}\atop{|\alpha|\geq 1}} 
\frac{\gamma!}{\alpha!\beta!}\int_{\partial\Omega} 
D^\alpha_X{\mathcal K}(X,Y)\,(P_\beta(X,Y)-P_\beta (X,Z))\,d\sigma_Y.
\end{equation}

\noindent If for a fixed $\mu>1$ and for each $X\in\Omega$ and $t>0$ 
we set 

\begin{equation}\label{def-Gamma}
\Gamma_t:=\{Y\in\partial\Omega:\,|X-Y|<\mu t\}
\end{equation}

\noindent we may then estimate 

\begin{eqnarray}\label{N1}
|D^\gamma{\mathcal E}\dot{f}(X)|^p &\leq & C
\sum_{{\alpha+\beta=\gamma}\atop{|\alpha|\geq 1}} 
\rho(X)^{-p|\alpha|}\meanint_{\Gamma_{\rho(X)}} 
|P_\beta(X,Y)-P_\beta(X,Z)|^p\,d\sigma_Y
\nonumber\\[6pt]
& \leq & C\sum_{{\alpha+\beta=\gamma}\atop{|\alpha|\geq 1}}  
\sum_{|\beta|+|\delta|\leq m-1}\rho(X)^{-p|\alpha|}
\meanint_{\Gamma_{\rho(X)}} 
|R_{\delta+\beta}(Y,Z)|^p\,|X-Y|^{p|\delta|}\,d\sigma_Y,
\nonumber\\[6pt]
& \leq & C\sum_{|\tau|\leq m-1}\rho(X)^{p(|\tau|-m)}
\meanint_{\Gamma_{\rho(X)}}|R_{\tau}(Y,Z)|^p\,d\sigma_Y,
\end{eqnarray} 

\noindent where we have used H\"older's inequality and (\ref{PR}).
Averaging the extreme terms in (\ref{N1}) for $Z$ in $\Gamma_{\rho(X)}$, 
we arrive at

\begin{equation}\label{N2}
|D^\gamma{\mathcal E}\dot{f}(X)|^p \leq C 
\sum_{|\tau|\leq m-1}\rho(X)^{p(|\tau|-m)-2(n-1)} 
\int_{\Gamma_{\rho(X)}}\!\int_{\Gamma_{\rho(X)}}
|R_\tau(Y,Z)|^p\,d\sigma_Yd\sigma_Z.
\end{equation}

Consider now a Whitney decomposition of $\Omega$ into a family of 
dyadic cubes, $\{Q_i\}_{i\in {\mathcal I}}$. In particular, 
$l_i:={\rm diam}\,Q_i\sim{\rm dist}\,(Q_i,\partial\Omega)$ 
uniformly for $i\in{\mathcal I}$. Thus, if $X\in Q_{i}$ for some 
${i}\in I_{j}:=\{i\in{\mathcal I}:\,l_i=2^{-j}\}$, $j\in\ZZ$, 
the estimate (\ref{N2}) yields 

\begin{equation}\label{EST-E}
|D^\gamma{\mathcal E}\dot{f}(X)|
\leq C\sum_{|\tau|\leq m-1}2^{-j(|\tau|-m)}
\Bigl(2^{2j(n-1)}\int\!\!\!\!\!\!\!\!\!\!\!\!\int
\limits_{{Y,Z\in\partial\Omega\cap\,\varkappa\,Q_{i}}\atop
{|Y-Z|<\varkappa\,2^{-j}}}|R_\tau(Y,Z)|^p\,d\sigma_Yd\sigma_Z\Bigr)^{1/p}
\end{equation}

\noindent for some $\varkappa=\varkappa(\partial\Omega)>1$.
In fact, by choosing the constant $\mu$ in (\ref{def-Gamma}) sufficiently 
close to $1$, matters can be arranged so that the family 
$\{\varkappa Q_i\}_{i\in{\mathcal I}}$ has finite overlap.  
Keeping this in mind and availing ourselves 
of the fact that $\rho(X)\sim l_i$ uniformly for 
$X\in Q_i$, $i\in{\mathcal I}$, for each multi-index $\gamma$ of length $m$ 
we may then estimate:

\begin{eqnarray}\label{bigstep}
&& \int_{\Omega}|D^\gamma{\mathcal E}\dot{f}(X)|^p\rho(X)^{p(1-s)-1}\,dX
\nonumber\\[6pt]
&& 
\quad
\leq C\sum_{j\in\ZZ}\sum_{i\in I_j} 2^{-jp(1-s)-j}
\int_{Q_i}|D^\gamma{\mathcal E}\dot{f}(X)|^p\,dX
\nonumber\\[6pt]
&& 
\quad
\leq C\sum_{j\in\ZZ}\sum_{i\in I_j}\sum_{|\tau|\leq m-1} 
2^{jp(m-1+s-|\tau|)+j(n-1)}
\int\!\!\!\!\!\!\!\!\!\!\!\!\int
\limits_{{Y,Z\in\partial\Omega\cap\,\varkappa\,Q_i}\atop
{|Y-Z|<\varkappa\,2^{-j}}}|R_\tau(Y,Z)|^p\,d\sigma_Yd\sigma_Z
\nonumber\\[6pt]
&& 
\quad
\leq C\sum_{j\in\ZZ}^\infty\sum_{|\tau|\leq m-1} 
2^{jp(m-1+s-|\tau|)+j(n-1)}
\int\!\!\!\!\!\!\!\!\!\!\!\!\int
\limits_{{Y,Z\in\partial\Omega}\atop
{|Y-Z|<\varkappa\,2^{-j}}}|R_\tau(Y,Z)|^p\,d\sigma_Yd\sigma_Z
\nonumber\\[6pt]
&& 
\quad
\leq C\sum_{|\tau|\leq m-1} 
\int_0^\infty\frac{r_\tau(t)^p}{t^{p(m-1+s-|\tau|)+n-1}}\,dt
\nonumber\\[6pt]
&& 
\quad
\leq C\|\dot{f}\|^p_{\dot{B}^{m-1+s}_p(\partial\Omega)},
\end{eqnarray}

\noindent where in the last step we have used (\ref{Bes-Nr}). 
This proves that the operator (\ref{Extension}) is well-defined and bounded.

\vskip 0.08in
\noindent{\it Step III: The right-invertibility property.} 
We shall now show that the operator (\ref{def-Ee}) is a right-inverse for 
the trace operator (\ref{TR-11}), i.e., whenever 
$\dot{f}=\{f_\gamma\}_{|\gamma|\leq m-1}\in
\dot{B}^{m-1+s}_p(\partial\Omega)$, there holds 

\begin{equation}\label{N3a}
f_\gamma=i^{|\gamma|}\,{\rm Tr}[D^\gamma{\mathcal E}\dot{f}]
\end{equation}

\noindent for every multi-index $\gamma$ of length $\leq m-1$. 
To this end, for $|\gamma|\leq m-1$ we write 

\begin{equation}\label{N4} 
D^\gamma{\mathcal E}\dot{f}(X)-{\mathcal E}_\gamma\dot{f}(X) 
=\sum_{{\alpha+\beta=\gamma}\atop{|\alpha|\geq 1}}
\frac{\gamma!}{\alpha!\beta!}\int_{\partial\Omega} 
D^\alpha_X {\mathcal K}(X,Y)(P_\beta(X,Y)-P_\beta(X,Z))\,d\sigma_Y,
\end{equation}

\noindent where

\begin{equation}\label{N3b}
{\mathcal E}_\gamma\dot{f}(X)
:=\int_{\partial\Omega}{\mathcal K}(X,Y)\, P_\gamma(X,Y)\,d\sigma_Y,
\qquad X\in\Omega. 
\end{equation}

\noindent Estimating the right-hand side in (\ref{N4}) in the same way 
as we did with the right-hand side of (\ref{N0}), we obtain

\begin{eqnarray}\label{Dgamma-E}
\int_{\partial\Omega}
|D^\gamma{\mathcal E}\dot{f}(X)-{\mathcal E}_\gamma\dot{f}(X)|^p 
\rho(X)^{-ps-1}\,dX
&\leq & C\sum_{|\tau|\leq m-1}\int_0^\infty\frac{r_\tau(t)^p}
{t^{p(|\gamma|+s-|\tau|)+n-1}}\,dt 
\nonumber\\[6pt]
& \leq & C\,\|\dot{f}\|^p_{\dot{B}^{m-1+s}_p(\partial\Omega)}.
\end{eqnarray}

\noindent In a similar fashion, we check that

\begin{eqnarray}\label{sim-fash}
&& \int_{\partial\Omega}|\nabla (D^\gamma{\mathcal E}\dot{f}(X)
-{\mathcal E}_\gamma\dot{f}(X)) |^p \rho(X)^{p-ps-1}\,dX
\nonumber\\[6pt]
&&\qquad\qquad 
\leq C\sum_{|\tau|\leq m-1}\int_0^\infty\frac{r_\tau(t)^p}
{t^{p(|\gamma|+s-|\tau|)+n-1}}\,dt
\leq C\,\|\dot{f}\|^p_{\dot{B}^{m-1+s}_p(\partial\Omega)}.
\end{eqnarray}

\noindent The two last inequalities imply 
$D^\gamma{\mathcal E}\dot{f}-{\mathcal E}\dot{f}\in V^{1,a}_p(\Omega)$ 
and, therefore,

\begin{equation}\label{N4a}
{\rm Tr}\,(D^\gamma{\mathcal E}\dot{f}-{\mathcal E}_\gamma\dot{f})=0.
\end{equation}

Going further, let us set

\begin{equation}\label{EgX}
Eg(X):=\int_{\partial\Omega}{\mathcal K}(X,Y)\,g(Y)\,d\sigma_Y, 
\qquad X\in\Omega.
\end{equation}

\noindent A simpler version of the reasoning in Step~II yields 
that $E$ maps $B^s_p(\partial\Omega)$ boundedly into $W^{1,a}_p(\Omega)$.
Also, a standard argument based on the Poisson kernel-like behavior of 
${\mathcal K}(X,Y)$ shows that ${\rm Tr}\,Eg=g$ for each 
$g\in B^s_p(\partial\Omega)$. 

Based on the definition (\ref{PPP}) and (\ref{N3b}) we have

\begin{eqnarray}\label{E-E}
&& |{\mathcal E}_\gamma\dot{f}(X)-Ef_\gamma(X)|^p
+\rho(X)^p|\nabla({\mathcal E}_\gamma\dot{f}(X)-Ef_\gamma(X))|^p 
\nonumber\\[6pt]
&& \qquad\qquad
\leq C\sum_{{|\beta|\leq m-1-|\gamma|}\atop{|\beta|\geq 1}} 
\rho(X)^{p|\beta|}\meanint_{\Gamma_{\rho(X)}} 
|f_{\gamma+\beta}(Y)|^p\,d\sigma_Y.
\end{eqnarray}

\noindent Consequently, for an arbitrary Whitney cube $Q_i$ we have

\begin{eqnarray}\label{n1}
&& \int_{Q_i}|{\mathcal E}_\gamma\dot{f}(X)-Ef_\gamma(X)|^p\rho(X)^{-ps-1}\,dX
+\int_{B_\delta}
|\nabla({\mathcal E}_\gamma\dot{f}(X)-Ef_\gamma(X))|^p\rho(X)^{p-ps-1}\,dX
\nonumber\\[6pt]
&& \qquad\qquad\qquad\qquad\qquad\qquad
\leq C\sum_{{|\beta|\leq m-1-|\gamma|}\atop{|\beta|\geq 1}}
l_i^{p(|\beta|-s)}\int_{\partial\Omega\cap \varkappa Q_i}
|f_{\gamma+\beta}(Y)|^p\,d\sigma_Y.
\end{eqnarray}

\noindent Summing over all Whitney cubes we find

\begin{equation}\label{Sm-wb}
\|{\mathcal E}_\gamma\dot{f}-Ef_\gamma\|_{V_p^{1,a}(\Omega)} 
\leq C\sum_{|\alpha|\leq m-1}\|f_\alpha\|_{L_p(\partial\Omega)}
\end{equation}

\noindent which implies

\begin{equation}\label{N5}
{\rm Tr}\,({\mathcal E}_\gamma\dot{f}-Ef_\gamma) =0.
\end{equation}

\noindent Finally, combining (\ref{N5}), (\ref{N4a}), and 
${\rm Tr}\,Ef_\gamma=f_\gamma$, we arrive at (\ref{N3a}). 

\vskip 0.08in
\noindent{\it Step IV: The kernel of the trace.} 
We now turn to the task of identifying the null-space of the trace operator 
(\ref{TR-11})-(\ref{Tr-DDD}). For each $k\in\NN_0$ we 
denote by ${\mathcal P}_k$ the collection of all vector-valued, complex 
coefficient polynomials of degree $\leq k$ (and agree that 
${\mathcal P}_k= 0$ whenever $k$ is a negative integer). 
The claim we make at this stage is that the null-space of the operator

\begin{equation}\label{TRW}
W^{m,a}_p(\Omega)\ni {\mathcal W}\mapsto 
\Bigl\{{\rm Tr}\,[D^\gamma{\mathcal W}]\Bigr\}_{|\gamma|=m-1}
\in B^s_p(\partial\Omega)
\end{equation}

\noindent is given by 

\begin{equation}\label{Null-tr}
{\mathcal P}_{m-2}+V^{m,a}_p(\Omega).
\end{equation}

\noindent The fact that the null-space of the trace operator
(\ref{TR-11})-(\ref{Tr-DDD}) is $V^{m,a}_p(\Omega)$ follows readily from this.

That (\ref{Null-tr}) is included in the null-space of the operator (\ref{TRW})
is obvious. The opposite inclusion amounts to showing that if 
${\mathcal W}\in W_p^{m,a}(\Omega)$ is such that 
${\rm Tr}\,[D^\gamma{\mathcal W}]=0$ for all multi-indices $\gamma$ with 
$|\gamma|=m-1$, then there exists $P_{m-2}\in {\mathcal P}_{m-2}$ with the 
property that ${\mathcal W}-P_{m-2}\in V_p^{m,a}(\Omega)$.  
To this end, we note that the case $m=1$ is a consequence of (\ref{Hardy}) 
and consider next the case $m=2$, i.e. when 

\begin{equation}\label{WWT}
{\mathcal W}\in W_p^{2,a}(\Omega),\qquad{\rm Tr}\,[\nabla{\mathcal W}]=0
\quad\mbox{on}\,\,\partial\Omega.
\end{equation}

\noindent Assume that $\{{\mathcal W}_j\}_{j\geq 1}$ is a 
sequence of smooth in $\overline{\Omega}$ (even polynomial) vector-valued 
functions approximating ${\mathcal W}$ in $W_p^{2,a}(\Omega)$. In particular, 
 
\begin{equation}\label{w1}
{\rm Tr}\,[\nabla{\mathcal W}_j]\to 0\,\,\,\mbox{ in }\,\,\,L_p(\partial\Omega)
\,\,\,\mbox{ as }\,\,j\to\infty. 
\end{equation}

\noindent If in a neighborhood of a point on $\partial\Omega$ the domain
$\Omega$ is given by $\{X:\,X_n>\varphi(X')\}$ for some Lipschitz function 
$\varphi$, the following chain rule holds for the gradient of the function 
$w_j:B'\ni X'\mapsto{\mathcal W}_j(X',\varphi(X'))$, where $B'$ is a 
$(n-1)$-dimensional ball: 
 
\begin{equation}\label{w2}
\nabla w_j(X')=
\Bigl(\nabla_{Y'}{\mathcal W}_j(Y',\varphi(X'))\Bigr)\Bigl|_{Y'= X'} 
+\Bigl(\frac{\partial}{\partial Y_n}{\mathcal W}_j(X',Y_n)\Bigr)
\Bigl|_{Y_n=\varphi(X')}\,\nabla\varphi(X').
\end{equation}

\noindent Since the sequence $\{w_j\}_{j\geq 1}$ is bounded in $L_p(B')$ and 
$\nabla w_j\to 0$ in $L_p(B')$, it follows that there exists a subsequence 
$\{j_i\}_i$ such that $w_{j_i}\to const$ in $L_p(B')$ 
(see Theorem~1.1.12/2 in {\bf\cite{Maz1}}). 
Hence, ${\rm Tr}\,{\mathcal W}=P_0=const$ on $\partial\Omega$. 
In view of ${\rm Tr}\,[{\mathcal W}-P_0]=0$ and 
${\rm Tr}\,[\nabla{\mathcal W}]=0$, we may conclude that 
${\mathcal W}-P_0\in V_p^{2,a}(\Omega)$ by Hardy's inequality. 

The general case follows in an inductive fashion, by reasoning as before 
with $D^\alpha{\mathcal W}$ with $|\alpha|=m-2$ in place of ${\mathcal W}$.
\hfill$\Box$
\vskip 0.08in

We now present a short proof of (\ref{Bes-X}), based on 
Proposition~\ref{trace-2}. 

\begin{proposition}\label{B-EQ}
Assume that $1<p<\infty$, $s\in(0,1)$ and $m\in\NN$. Then 

\begin{equation}\label{eQ-11}
\|\dot{f}\|_{\dot{B}^{m-1+s}_p(\partial\Omega)}\sim \sum_{|\alpha|\leq m-1}
\|f_\alpha\|_{B^{s}_p(\partial\Omega)},
\end{equation}

\noindent uniformly for 
$\dot{f}=\{f_\alpha\}_{|\alpha|\leq m-1}\in\dot{B}^{m-1+s}_p(\partial\Omega)$.
As a consequence, (\ref{Bes-X}) holds. 
\end{proposition}

\noindent{\bf Proof.} The left-pointing inequality in (\ref{eQ-11})
is implicit in (\ref{RRR=est}). As for the opposite one, let  
$\dot{f}=\{f_\alpha\}_{|\alpha|\leq m-1}\in\dot{B}^{m-1+s}_p(\partial\Omega)$
and, with $a:=1-s-1/p$, consider 
${\mathcal U}:={\mathcal E}(\dot{f})\in W^{m,a}_p(\Omega)$. 
Then Lemma~\ref{trace-1} implies that, for each multi-index 
$\alpha$ of length $\leq m-1$, the function 
$f_\alpha=i^{|\alpha|}\,{\rm Tr}\,[D^\alpha{\mathcal U}]$ belongs 
to $B^s_p(\partial\Omega)$, plus a naturally accompanying norm estimate. 
This concludes the proof of (\ref{eQ-11}). Finally, the last claim in 
the proposition is a consequence of (\ref{eQ-11}), (\ref{dense}) and 
the fact that the operator (\ref{TR-11})-(\ref{Tr-DDD}) is onto. 
\hfill$\Box$
\vskip 0.08in

We include one more equivalent characterization of the space 
$\dot{B}^{m-1+s}_p(\partial\Omega)$, in the spirit of work in 
{\bf\cite{AP}}, {\bf\cite{PV}}, {\bf\cite{Ve2}}. 
To state it, recall that $\{e_j\}_j$ is the canonical 
orthonormal basis in $\RR^n$. 

\begin{proposition}\label{CC-Aray}
Assume that $1<p<\infty$, $s\in(0,1)$ and $m\in\NN$. Then 

\begin{equation}\label{eQ0}
\{f_\alpha\}_{|\alpha|\leq m-1}\in\dot{B}^{m-1+s}_p(\partial\Omega)
\Longleftrightarrow
\left\{
\begin{array}{l}
f_\alpha\in B^{s}_p(\partial\Omega),\quad \forall\,\alpha\,:\,|\alpha|\leq m-1
\\[10pt]
\qquad\qquad\qquad
\mbox{ and }
\\[10pt]
(\nu_j\partial_k-\nu_k\partial_j)f_{\alpha}=
\nu_jf_{\alpha+e_k}-\nu_kf_{\alpha+e_j}
\\[6pt]
\forall\,\alpha\,:\,|\alpha|\leq m-2,\quad\forall\,j,k\in\{1,...,n\}.
\end{array}
\right.
\end{equation}
\end{proposition}

\noindent{\bf Proof.} The left-to-right implication is a consequence 
of (\ref{eQ-11}) and of the fact that (\ref{ff-aa}) holds for 
some ${\mathcal U}\in W^{m,a}_p(\Omega)$ (cf. Proposition~\ref{trace-2}). 
As for the opposite implication, we proceed as in the proof of 
Proposition~\ref{trace-2} and estimate (\ref{rem-Rr}) based on 
the identities (\ref{B-CC}) and knowledge that $f_\alpha$ belongs to 
$B^{s}_p(\partial\Omega)$ for each $\alpha$ of length $\leq m-1$. 
\hfill$\Box$
\vskip 0.08in

We close this section with two remarks on the nature of the space 
$\dot{B}^{m-1+s}_p(\partial\Omega)$. First, we claim that the assignment 

\begin{equation}\label{ass}
\dot{B}^{m-1+s}_p(\partial\Omega)\ni\dot{f}=\{f_\alpha\}_{|\alpha|\leq m-1}
\mapsto 
\Bigl\{i^k\sum_{|\alpha|=k}\frac{k!}{\alpha!}\,\nu^\alpha\,f_\alpha\Bigr\}
_{0\leq k\leq m-1}\in L_p(\partial\Omega)
\end{equation}

\noindent is one-to-one. This is readily justified with the help of the 
identity 

\begin{equation}\label{m-xxx}
D^\alpha=i^{-|\alpha|}\,\nu^\alpha\,
\frac{\partial^{|\alpha|}}{\partial\nu^{|\alpha|}}
+\sum_{|\beta|=|\alpha|-1}\sum_{j,k=1}^n p_{\alpha,\beta,j,k}(\nu)
\frac{\partial}{\partial\tau_{jk}}D^\beta
\end{equation}

\noindent where $\partial/\partial\tau_{jk}:=\nu_j\partial/\partial x_k
-\nu_k\partial/\partial x_j$ and the $p_{\alpha,\beta,j,k}$'s are polynomial 
functions. Indeed, let $\dot{f}\in \dot{B}^{m-1+s}_p(\partial\Omega)$ be 
mapped to zero by the assignment (\ref{ass}) and consider 
${\mathcal U}:={\mathcal E}(\dot{f})\in W^{m,a}_p(\Omega)$. Then 
$f_\alpha=i^{|\alpha|}\,{\rm Tr}\,\,[D^\alpha\,{\mathcal U}]$ 
on $\partial\Omega$ for each $\alpha$ with $|\alpha|\leq m-1$ and, 
granted the current hypotheses, $\partial^k{\mathcal U}/\partial\nu^k=0$ 
for $k=0,1,...,m-1$. Consequently, (\ref{m-xxx}) and induction on $|\alpha|$
yield that ${\rm Tr}\,\,[D^\alpha\,{\mathcal U}]=0$ on $\partial\Omega$ 
for each $\alpha$ with $|\alpha|\leq m-1$. Thus, $f_\alpha=0$ for each 
$\alpha$ with $|\alpha|\leq m-1$, as desired. 

The elementary identity (\ref{m-xxx}) can be proved by writing 

\begin{eqnarray}\label{m-xxx2}
i^{|\alpha|}\,D^\alpha & = & \prod_{j=1}^n\Bigl(\frac{\partial}{\partial x_j}
\Bigr)^{\alpha_j}
\\[6pt]
& = & \prod_{j=1}^n\Bigl[\sum_{k=1}^n\xi_k\Bigl(\xi_k
\frac{\partial}{\partial x_j}-\xi_j\frac{\partial}{\partial x_k}\Bigr)
+\sum_{k=1}^n\xi_j\xi_k\frac{\partial}{\partial x_k}\Bigr]^{\alpha_j}
\Bigl|_{\xi=\nu}
\nonumber\\[6pt]
& = & \prod_{j=1}^n\Bigl[\sum_{l=0}^{\alpha_j}\frac{\alpha_j!}{l!(\alpha_j-l)!}
\Bigl(\sum_{k=1}^n\xi_k\Bigl(\xi_k
\frac{\partial}{\partial x_j}-\xi_j\frac{\partial}{\partial x_k}\Bigr)\Bigr)
^{\alpha_j-l}\nu_j^{l}\frac{\partial^l}{\partial\nu^l}\Bigr]\Bigl|_{\xi=\nu}
\nonumber\\[6pt]
& = & \prod_{j=1}^n\Bigl[
\nu_j^{\alpha_j}\frac{\partial^{\alpha_j}}{\partial\nu^{\alpha_j}}+
\sum_{l=0}^{\alpha_j-1}\frac{\alpha_j!}{l!(\alpha_j-l)!}
\Bigl(\sum_{k=1}^n\xi_k\Bigl(\xi_k
\frac{\partial}{\partial x_j}-\xi_j\frac{\partial}{\partial x_k}\Bigr)\Bigr)
^{\alpha_j-l}\nu_j^{l}\frac{\partial^l}{\partial\nu^l}\Bigr]\Bigl|_{\xi=\nu}
\nonumber
\end{eqnarray}

\noindent and noticing that 
$\prod_{j=1}^n\nu_j^{\alpha_j}\partial^{\alpha_j}/\partial\nu^{\alpha_j}
=\nu^\alpha\partial^{|\alpha|}/\partial\nu^{|\alpha|}$, 
whereas $(\xi_k\partial/\partial x_j-\xi_j\partial/\partial x_k)|_{\xi=\nu}
=-\partial/\partial\tau_{jk}$. 

Our second remark concerns the image of the mapping (\ref{ass}) in the 
case when $\partial\Omega$ is sufficiently smooth. More precisely, 
assume that $\partial\Omega\in C^{m-1,1}$ and, for $0\leq k\leq m-1$,  
the space $B^{m-1-k+s}_p(\partial\Omega)$ is defined starting from 
$B^{m-1-k+s}_p(\RR^{n-1})$ and then transporting this space to
$\partial\Omega$ via a smooth partition of unity argument and locally 
flattening the boundary (alternatively, $B^{m-1-k+s}_p(\partial\Omega)$ is 
the image of the trace operator acting from $B^{m-1-k+s+1/p}_p(\RR^{n})$). 
We claim that 

\begin{equation}\label{image}
\partial\Omega\in C^{m-1,1}\Longrightarrow
\mbox{the image of the mapping (\ref{ass}) is 
$\oplus_{k=0}^{m-1}B^{m-1-k+s}_p(\partial\Omega)$}.
\end{equation}

\noindent Indeed, granted that $\partial\Omega\in C^{m-1,1}$, 
it follows from (\ref{eQ0}) that 
$f_\alpha\in B^{m-1-|\alpha|+s}_p(\partial\Omega)$ for each $\alpha$ with 
$|\alpha|\leq m-1$ and, hence, 
$g_k:=i^k\sum_{|\alpha|=k}\frac{k!}{\alpha!}\,\nu^\alpha\,f_\alpha
\in B^{m-1-k+s}_p(\partial\Omega)$ for each $k\in\{0,...,m-1\}$. 

Conversely, given a family 
$\{g_k\}_{0\leq k\leq m-1}\in\oplus_{k=0}^{m-1}B^{m-1-k+s}_p(\partial\Omega)$,
we claim that there exists
$\dot{f}=\{f_\alpha\}_{|\alpha|\leq m-1}\in\dot{B}^{m-1+s}_p(\partial\Omega)$
such that $g_k=i^k\sum_{|\alpha|=k}\frac{k!}{\alpha!}\,\nu^\alpha\,f_\alpha$ 
for each $k\in\{0,...,m-1\}$. 
One way to see this is to start with ${\mathcal U}\in B^{m-1+s+1/p}_p(\Omega)$
solution of $\Delta^m{\mathcal U}=0$ in $\Omega$, 
$\partial^k{\mathcal U}/{\partial\nu^k}=g_k
\in B^{m-1-k+s}_p(\partial\Omega)$, $0\leq k\leq m-1$ 
(a system which satisfies the Shapiro-Lopatinskij condition) and then define 
the $f_\alpha$'s as in (\ref{ff-aa}).

\subsection{Proof of the main result and further comments}}

Theorem~\ref{Theorem} is a particular case of the next theorem concerning 
the unique solvability of the Dirichlet problem in $W_p^{m,a}(\Omega)$. 

\begin{theorem}\label{Theorem1}
Let all assumptions of Theorem~\ref{th1a} be satisfied. Also let 
${\mathcal F}\in V_p^{-m,a}(\Omega)$. Then the Dirichlet problem
 
\begin{equation}\label{m9}
\left\{
\begin{array}{l}
{\mathcal A}(X,D_X)\,{\mathcal U}={\mathcal F}
\qquad\mbox{in}\,\,\Omega,
\\[15pt] 
{\displaystyle{\frac{\partial^k{\mathcal U}}{\partial\nu^k}}}
=g_k\quad\,\,\mbox{on}\,\,\partial\Omega,\,\,\,\,\,0\leq k\leq m-1,
\end{array}
\right.
\end{equation}
 
\noindent has a solution ${\mathcal U}\in W_p^{m,a}(\Omega)$ if and only if 
(\ref{data-B}) is satisfied. In this latter case, the solution is unique
and satisfies

\begin{equation}\label{estUU}
\|{\mathcal U}\|_{W_p^{m,a}(\Omega)}
\leq C\sum_{|\alpha|\leq m-1}\|f_\alpha\|_{B^{s}_p(\partial\Omega)}
+C\|{\mathcal F}\|_{V_p^{-m,a}(\Omega)}. 
\end{equation}
\end{theorem}

\noindent{\bf Proof.} It is clear from definitions that the operator 

\begin{equation}\label{ADXW}
{\mathcal A}(X,D_X):W_p^{m,a}(\Omega)\longrightarrow V_p^{-m,a}(\Omega)
\end{equation}

\noindent is well-defined and bounded. Thus, granted that we seek solutions
for (\ref{m9}) in the space $W_p^{m,a}(\Omega)$, the membership of 
${\mathcal F}$ to $V_p^{-m,a}(\Omega)$, as well as the fact that the $g_k$'s 
satisfy (\ref{data-B}), are necessary conditions for the solvability 
of (\ref{m9}).

Conversely, let 
$\dot{f}=\{f_\alpha\}_{|\alpha|\leq m-1}\in\dot{B}^{m-1+s}_p(\partial\Omega)$
be such that (\ref{data-B}) holds and, with ${\mathcal E}$ denoting the 
extension operator from Proposition~\ref{trace-2}, seek a solution for
(\ref{m9}) in the form ${\mathcal U}={\mathcal E}(\dot{f})+{\mathcal W}$. 
where ${\mathcal W}\in V^{m,a}_p(\Omega)$ solves

\begin{equation}\label{m10}
\left\{
\begin{array}{l}
{\mathcal A}(X,D_X){\mathcal W}
={\mathcal F}-{\mathcal A}(X,D_X)({\mathcal E}(\dot{f}))
\qquad\mbox{in}\,\,\Omega,
\\[15pt] 
{\rm Tr}\,\,[D^\gamma\,{\mathcal W}]=0\quad\mbox{on}\,\,\partial\Omega,
\,\,\,\forall\,\gamma\,:\,|\gamma|\leq m-1.
\end{array}
\right.
\end{equation}

\noindent Since the boundary conditions in (\ref{m10}) are automatically
satisfied if ${\mathcal W}\in V_p^{m,a}(\Omega)$, the solvability of 
(\ref{m10}) is a direct consequence of Theorem~\ref{th1a}. As for uniqueness, 
assume that ${\mathcal U}\in W_p^{m,a}(\Omega)$ solves (\ref{m9}) with 
${\mathcal F}=0$ and $g_k=0$, $0\leq k\leq m-1$. From the fact that 
(\ref{ass}) is one-to-one, we infer that 
${\rm Tr}\,\,[D^\gamma\,{\mathcal U}]=0$ on $\partial\Omega$ 
for all $\gamma$ with $|\gamma|\leq m-1$. Then, by Proposition~\ref{trace-2}, 
${\mathcal U}\in V_p^{m,a}(\Omega)$ is a null-solution of 
${\mathcal A}(X,D_X)$. In turn, Theorem~\ref{th1a} gives that 
${\mathcal U}=0$, proving uniqueness for (\ref{m9}). Finally, (\ref{estUU})
is a consequence of the results in \S{6.4}.
\hfill$\Box$
\vskip 0.08in

We conclude this section with a couple of comments, the first of which 
regards the effect of the presence of lower order terms. 
More specifically, assume that 

\begin{equation}\label{E444-bis}
{\mathcal A}(X,D_X)\,{\mathcal U}
:=\sum_{0\leq |\alpha|,|\beta|\leq m}D^\alpha({\mathcal A}_{\alpha\beta}(X)
\,D^\beta{\mathcal U}),\qquad X\in\Omega,
\end{equation}

\noindent where the top part of ${\mathcal A}(X,D_X)$ satisfies the 
hypotheses made in Theorem~\ref{Theorem} and the lower order terms are bounded.
Then the Dirichlet problem (\ref{m9}) is Fredholm solvable, of index zero, 
in the sense that a solution ${\mathcal U}\in W_p^{m,a}(\Omega)$ exists 
if and only if the data ${\mathcal F}$, $\{g_k\}_{0\leq k\leq m-1}$ 
satisfy finitely many linear conditions, whose number matches the 
dimension of the space of null-solutions for (\ref{m9}). 
Furthermore, the estimate 

\begin{equation}\label{estUU-bis}
\|{\mathcal U}\|_{W_p^{m,a}(\Omega)}
\leq C\, \Bigl(\, \sum_{|\alpha|\leq m-1}\|{\mathcal F}\|_{V_p^{-m,a}(\Omega)} 
+\|f_\alpha\|_{B^{s}_p(\partial\Omega)}+\|{\mathcal U}\|_{L_p(\Omega)}\Bigr)
\end{equation}
 
\noindent holds for any solution ${\mathcal U}\in W_p^{m,a}(\Omega)$
of (\ref{m9}). 

Indeed, the operator

\begin{equation}\label{cal-AAA}
{\mathcal A}:V_p^{m,a}(\Omega)\longrightarrow V_p^{-m,a}(\Omega)
\end{equation}

\noindent is Fredholm with index zero, as can be seen by decomposing 
${\mathcal A}=\ring{\mathcal A}+({\mathcal A}-\ring{\mathcal A})$ where
$\ring{\mathcal A}:=\sum_{|\alpha|=|\beta|=m} 
D^\alpha {\mathcal A}_{\alpha\beta}\,D^\beta$, and then invoking 
Theorem~\ref{th1a}. 
Now, it can be shown that the problem (\ref{m9}) is solvable if and only if 
${\mathcal F}-{\mathcal A}(X,D_X){\mathcal E}\dot{f}
\in {\rm Im}\,{\mathcal A}$, the image of the operator (\ref{cal-AAA}). 
Thus, if $T({\mathcal F},\{g_k\}_{0\leq k\leq m-1}):=
{\mathcal F}-{\mathcal A}(X,D_X){\mathcal E}\dot{f}$, this membership entails 
$({\mathcal F},\{g_k\}_{0\leq k\leq m-1})\in T^{-1}
\Bigl({\rm Im}\,{\mathcal A}\Bigr)$. Note that $T$ maps the space of data 
onto $V_p^{-m,a}(\Omega)$, hence the number of linearly independent 
compatibility conditions the data should satisfy is 

\begin{equation}\label{comp-cond-X}
{\rm codim}\,T^{-1}\Bigl({\rm Im}\,{\mathcal A}\Bigr)
={\rm codim}\,({\rm Im}\,{\mathcal A}).
\end{equation}

\noindent On the other hand, from by Proposition~\ref{trace-2} and the 
fact that (\ref{ass}) is one-to-one we infer that the space of null-solutions 
for (\ref{m9}) is precisely ${\rm ker}\,{\mathcal A}$, the kernel of the
operator (\ref{cal-AAA}). Since, as already pointed out,
this operator has index zero, it follows that the problem (\ref{m9})
has index zero. Finally, (\ref{estUU-bis}) follows from what we have proved so 
far via a standard reasoning as in {\bf\cite{Ho}}. 

Our last comment regards the statement of the Dirichlet
problem (\ref{e0}) with data 

\begin{equation}\label{newdata}
\partial^k{\cal U}/\partial\nu^k=g_k\in B_p^{m-1-k+s}(\partial\Omega),
\qquad k=0,1,\ldots,m-1,
\end{equation}

\noindent where $B_p^{m-1-k+s}(\partial\Omega)$ is defined here as the
{\it range of ${\rm Tr}$ acting from $B_p^{m-1-k+s+1/p}(\RR^n)$}. 
If $\partial\Omega$ is smooth ($C^{1,1}$ will do) this problem is, certainly, 
well-posed. Let us illustrate some features of this particular formulation 
as the smoothness of $\partial\Omega$ deteriorates. 

Suppose we are looking for the solution ${\cal U}\in W_2^2(\Omega)$ of the
Dirichlet problem for the biharmonic operator

\begin{equation}\label{m9b}
\left\{
\begin{array}{l}
\Delta ^2\,{\cal U}=0\qquad\mbox{in}\,\,\Omega,
\\[10pt]
{\rm Tr}\,{\cal U}=g_1\qquad\mbox{on}\,\,\partial\Omega,
\\[10pt]
\langle\nu,{\rm Tr}\,[\nabla{\cal U}]\rangle 
=g_2\quad\mbox{on}\,\,\partial\Omega.
\end{array}
\right.
\end{equation}

\noindent The simplest class of data $(g_1, g_2)$  would be, of course,
$B_2^{3/2}(\partial\Omega)\times B_2^{1/2}(\partial\Omega)$, where
$B_2^{3/2}(\partial\Omega)$ and $B_2^{1/2}(\partial\Omega)$ are the spaces
of traces on $\partial\Omega$ for functions in $ W_2^2(\Omega)$ and $
W_2^1(\Omega)$, respectively. However, this formulation has several serious
drawbacks.

The first one is that the mapping

\begin{equation}\label{w22}
W_2^2(\Omega)\ni {\cal U}\to\langle\nu,{\rm Tr}\,[\nabla{\mathcal U}]\rangle 
\in B_2^{1/2}(\partial\Omega)
\end{equation}

\noindent is generally unbounded. In fact, by choosing ${\mathcal U}$ to be a
linear function we see that the continuity of (\ref{w22}) implies 
$\nu\in B_2^{1/2}(\partial\Omega)$ which is not necessarily the case 
for a Lipschitz domain, even for such a simple one as the square $S=[0,1]^2$.

The same problem fails to have a solution in the class in $W_2^2(\Omega)$
when when $(g_1, g_2)$ is an arbitrary pair in
$B_2^{3/2}(\partial\Omega)\times B_2^{1/2}(\partial\Omega)$. Indeed,
consider the problem (\ref{m9b}) for $\Omega=S$ and the data $g_1 =0$
and $g_2 =1$. It is standard (see Theorem 7.2.4 in {\bf\cite{KMR1}} 
and Sect.\,7.1 in {\bf\cite{KMR2}}) that the main term of the asymptotics 
near the origin of any solution ${\mathcal U}$ in $W_2^1(S)$ is given 
in polar coordinates $(r,\omega)$ by

\begin{equation}\label{asymp}
\frac{2r}{\pi+2}\left( (\omega -\frac{\pi}{2})
\sin\omega - \omega\cos\omega\right).
\end{equation}

\noindent Since this function does not belong to $W_2^2(S)$, there 
is no solution of problem (\ref{m9b}) in this space.

\vskip 0.10in
\noindent --------------------------------------
\vskip 0.20in

\noindent {\tt Vladimir Maz'ya}

\noindent Department of Mathematics

\noindent Ohio State University 

\noindent Columbus, OH 43210, USA

\noindent {\tt e-mail}: {\it vlmaz\@@math.ohio-state.edu}

and 

\noindent Department of Mathematical Sciences

\noindent University of Liverpool

\noindent Liverpool L69 3BX, UK

\vskip 0.15in

\noindent {\tt Marius Mitrea}

\noindent Department of Mathematics

\noindent University of Missouri at Columbia

\noindent Columbia, MO 65211, USA

\noindent {\tt e-mail}: {\it marius\@@math.missouri.edu}

\vskip 0.15in

\noindent {\tt Tatyana Shaposhnikova}

\noindent Department of Mathematics

\noindent Ohio State University 

\noindent Columbus, OH 43210, USA

\noindent {\tt e-mail}: {\it tasha@math.ohio-state.edu}

and

\noindent Department of Mathematics

\noindent Link\"oping University

\noindent Link\"oping SE-581 83, Sweden

\noindent {\tt e-mail}: {\it tasha@mai.liu.se

\end{document}